\def\E{\end{document}}
\documentclass[10pt]{article}
\usepackage{cite}
\usepackage[tbtags]{amsmath}
\usepackage{enumerate}
\usepackage{amsmath,amssymb,amsthm,dsfont}
\usepackage{csquotes}
\usepackage{appendix}
\usepackage{enumitem}
\usepackage{titlesec,hyperref}
\usepackage{color}
\usepackage{verbatim}
\usepackage{fancyhdr}
\usepackage[margin=2.5cm]{geometry}
\usepackage{indentfirst} % 首行缩进

\numberwithin{equation}{section}
\newcommand{\eqrow}[2]{\ensuremath{\text{\eqref{#1}}_{#2}}}
\newtheorem{theorem}{Theorem}[section]
\newtheorem{proposition}[theorem]{Proposition}
\newtheorem{lemma}[theorem]{Lemma}
\newtheorem{corollary}[theorem]{Corollary}
\theoremstyle{definition}
\newtheorem{definition}[theorem]{Definition}

\theoremstyle{remark}
\newtheorem{remark}[theorem]{Remark}
\newcommand{\Om}{\Omega}
\newcommand{\OO}{\mathcal O}

\newcommand{\dd}{\,\mathrm d}
\newcommand{\R}{\mathbb R}
\newcommand{\eps}{\varepsilon}
\newcommand{\Div}{\operatorname{div}}
\newcommand{\curl}{\operatorname{curl}}
\newcommand{\supp}{\operatorname{supp}}
\newcommand{\norm}[2][]{\left\|#2\right\|\if\relax\detokenize{#1}\relax\else_{#1}\fi}

\newcommand{\X}{\mathcal X}
\newcommand{\Y}{\mathcal Y}
\newcommand{\F}{\mathcal F}
\newcommand{\Uc}{\mathcal U}
\newcommand{\Hh}{\mathcal H}
\newcommand{\Vv}{\mathcal V}

\setlist{nosep}

\newcommand{\NN}{\mathcal N_\beta}

\newcommand{\tr}{\mathrm{tan}}

\newcommand{\nb}{\overline n}
\newcommand{\cb}{\overline c}
\newcommand{\ub}{\overline u}
\newcommand{\pb}{\overline p}
\newcommand{\an}{\overline{\alpha}_n}
\newcommand{\ac}{\overline{\alpha}_c}

\begin{document}
\title{Small-time global boundary controllability of chemotaxis--fluid systems under nonnegativity constraints}

	\author{Qiang Tao  \thanks {School of Mathematical Sciences,
			Shenzhen University, Shenzhen 518060, China. E-mail: taoq@szu.edu.cn; taoq060@126.com} \quad and  \quad Ying Yang\thanks {Corresponding author. School of Mathematical Sciences,
			Shenzhen University, Shenzhen 518060, China. E-mail: yingyang@szu.edu.cn}
	}

\date{}
\maketitle

{\bf Abstract.}\ \ {\small In this paper, we study the global boundary controllability of a coupled Keller-Segel-Navier-Stokes system in dimensions two and three. The boundary controls only act on a part of the boundary which intersects each connected component. On the remaining boundary, the biological variables satisfy homogeneous Neumann conditions and the velocity satisfies a Navier slip--with--friction condition with a smooth, time-independent symmetric friction matrix.
We establish global approximate and exact controllability to admissible positive, possibly spatially nonconstant, reference trajectories in any positive time. The initial data need not be close to the reference, and either biological component may initially vanish. Both biological components remain nonnegative throughout the control process.
The approach connects a positive preparation, the return method, appropriate asymptotic boundary layer expansions and a local controllability result.
A common $L^4\times W^{1,4}\times L^4$ framework  closes the nonlinear stability argument in both dimensions.

}



\section{Introduction}
Chemotaxis is a biological process in which cells move toward a chemically more favorable environment. This movement can be toward or away from a higher concentration of a chemoattractant.
Keller and Segel \cite{KS1970, KS1971} first established a chemotaxis model, and since then many modified models have been developed to describe bacterial aggregation, slime mold formation, fish pigmentation patterning, tumor angiogenesis, primitive streak formation, wound healing, etc. In fact, bacteria and microorganisms often live in fluids, where both the chemical concentration (e.g.
oxygen) and the cell density are transported and diffuse through the fluid. When cells are suspended in a fluid, their motion is affected by diffusion, chemotactic transport, and fluid advection, while the cell density contributes to the force in the momentum equation. Tuval et al. \cite{TCDWKG2005} proposed a coupled system of the chemotaxis model and a viscous incompressible fluid.

The controllability of such a system concerns the possibility of steering the coupled evolution towards a prescribed trajectory by an external action. Since the cell density and chemical concentration represent nonnegative quantities, a control construction should preserve their signs throughout the evolution.
The purpose of this paper is to establish global boundary controllability to positive, possibly spatially nonconstant, trajectories of a coupled Keller–Segel–Navier–Stokes system with homogeneous Neumann and Navier slip--with--friction boundary conditions in 2D or 3D domain. The initial data need not be close to the reference trajectories. 
\subsection{Setting of the problem}
Let $d\in\{2,3\}$ and $\Om\subset\R^d$ be a bounded connected domain with smooth boundary. Denote $\Gamma_c$ to be a nonempty open subset of $\partial\Om$ meeting every connected component  and  $\Gamma_0=\partial\Om\setminus\Gamma_c$. Let $Q_T^\Omega=(0,T)\times\Omega$, $\Sigma_T^0=(0,T)\times\Gamma_0$ and $\Sigma_T^c=(0,T)\times\Gamma_c$. $e_2=(0,1)\in\mathbb R^2$ and $e_3=(0, 0, 1)\in\mathbb R^3$. Let $s$ denote physical time and $\nu$ the outward unit normal. We consider the chemotaxis--fluid system with consumption term
\begin{equation}\label{1.1}
\left\{\begin{aligned}
 &n_s+u\cdot\nabla n-D_n\Delta n=-\chi\Div(n\nabla c), &&\text{in }Q_T^\Omega,\\[2pt]
 &c_s+u\cdot\nabla c-D_c\Delta c=-nc, &&\text{in }Q_T^\Omega,\\[2pt]
 &u_s+(u\cdot\nabla)u-\mu\Delta u+\nabla p=ne_d,  \quad
 \Div u=0, &&\text{in }Q_T^\Omega,\\[2pt]
 &\partial_\nu n=\partial_\nu c=0,\quad u\cdot\nu=0,\quad\NN u=0, &&\text{on }\Sigma_T^0,\\[2pt]
 &\partial_\nu n=h_n^*,\quad\partial_\nu c=h_c^*,  \quad u\cdot\nu=h_\nu^*,\quad\NN u=h_\tau^*, &&\text{on }\Sigma_T^c,\\[2pt]
 &(n,c,u)|_{s=0}=(n_0,c_0,u_0), &&\text{in }\Omega.
\end{aligned}\right.
\end{equation}
where $D_n,D_c,\mu,\chi>0$ are fixed and $n_0,c_0\ge0$. The functions $n,c,u,p$ denote cell density, chemical concentration, velocity, and pressure, respectively. The two biological variables are required to satisfy
\begin{equation*}
 n(s,\cdot)\ge0,\qquad c(s,\cdot)\ge0\qquad(0\le s\le T).
\end{equation*}
Define
\begin{equation*}
 D(u)=\frac{\nabla u+\nabla u^{\mathsf T}}2,\qquad
 [z]_{\tr}=z-(z\cdot\nu)\nu,\qquad
 \NN u=[2\mu D(u)\nu+\beta(x)u]_{\tr},
\end{equation*}
where $\beta\in C^\infty(\partial\Om;\mathbb R^{d\times d})$ denotes a symmetric $d\times d$ matrix-valued function regarded as a measure of the friction near the boundary. Set the tangential projection
\begin{equation*}
 P_\tau=I-\nu\otimes\nu,\qquad
 \beta_\tau=P_\tau\beta P_\tau.
\end{equation*}
For tangent $u$, $[\beta u]_{\tr}=\beta_\tau u$.
The controls $h_n^*,h_c^*,h_\nu^*,h_\tau^*$ act on $\Gamma_c$. The incompressibility condition requires
\begin{equation}\label{1.2}
 \int_{\Gamma_c}h_\nu^*(s)\dd S=0.
\end{equation}
In the 2D case, set $\operatorname{curl}u=\partial_1u_2-\partial_2u_1$ and $\tau=(-\nu_2,\nu_1)$,
$b_\tau=\tau^{\mathsf T}\beta\tau$ and $b_\nu=\tau^{\mathsf T}\beta\nu$. Then
\begin{equation*}
 \NN u\cdot\tau
 =2\mu D(u)\nu\cdot\tau+b_\tau(u\cdot\tau)+b_\nu(u\cdot\nu).
\end{equation*}
On $\Gamma_0$ the last term is zero, while on $\Gamma_c$ it is
$b_\nu h_\nu^*$.

The target reference trajectory $(\nb,\cb,\ub,\pb)$ satisfies the same interior equations and homogeneous boundary conditions on $\Gamma_0$:
\begin{equation*}
\left\{\begin{aligned}
 &\nb_s+\ub\cdot\nabla\nb-D_n\Delta\nb=-\chi\Div(\nb\nabla\cb), &&\text{in }Q_T^\Omega,\\[2pt]
 &\cb_s+\ub\cdot\nabla\cb-D_c\Delta\cb=-\nb\cb, &&\text{in }Q_T^\Omega,\\[2pt]
 &\ub_s+(\ub\cdot\nabla)\ub-\mu\Delta\ub+\nabla\pb=\nb e_d,  \quad \Div\ub=0, &&\text{in }Q_T^\Omega,\\[2pt]
 &\partial_\nu\nb=\partial_\nu\cb=0,\quad\ub\cdot\nu=0,\quad\NN\ub=0, &&\text{on }\Sigma_T^0,
 \\[2pt]
 &(\nb,\cb,\ub)|_{s=0}=(\bar n_0,\bar c_0,\bar u_0), &&\text{in }\Omega.
\end{aligned}\right.
\end{equation*}
Our goal is to find controls such that
\begin{equation*}
 (n,c,u)(T,\cdot)=(\nb,\cb,\ub)(T,\cdot).
\end{equation*}

\subsection{Previous results and main difficulties}

In the past decade, with the deepening study of the well-posedness of biological chemotaxis models, an increasing number of works have focused on the mathematical control theory of chemotaxis models.
There are some important results in the study of the controllability of the related  models. Guo and Zhang \cite{GZ2014} first considered local null controllability for a parabolic-elliptic chemotaxis model. Then, Chaves-Silva and Guerrero \cite{CG2015} established the local controllability of a parabolic–parabolic Keller–Segel system near a constant equilibrium, with the control entering through the chemical equation, uniformly with respect to the parameter governing the parabolic–elliptic limit. In \cite{GZ2016}, the authors obtained local exact controllability toward positive trajectories by means of a single bounded internal control, while preserving positivity.
For logarithmic sensitivity, Tao and Zhang \cite{TZ2025} used a Hopf-Cole transformation to establish local exact controllability to nonnegative trajectories for the transformed system. However, these results are all concerned with local controllability based on Carleman
estimates for the adjoint system.

The controllability of incompressible fluid equations
has been extensively discussed by various methods,
including the duality method and the return method.
For the incompressible Navier--Stokes equations,
local exact controllability to trajectories with
internal controls was established by Imanuvilov
\cite{I2001} and Fern\'andez-Cara et al.
\cite{FGIP2004}. In \cite{G2006}, Guerrero treated nonlinear Navier-slip boundary conditions. In recent years, there has been growing interest in global controllability of Navier-Stokes equations.  
Coron, Marbach and Sueur \cite{CMS2020} extended the well-prepared dissipation method, originally developed for the global null controllability of the Burgers equation, combined with the return method to prove small-time global exact controllability under Navier slip-with-friction conditions. 
In \cite{CMSZ2019}, Coron et al. considered two-dimensional incompressible Navier–Stokes equations in a rectangle with no-slip conditions. They obtained that any finite energy initial data can be driven to rest by boundary controls together with a small force in any Sobolev norm.
Liao, Sueur and Zhang \cite{LSZ2022} constructed smooth controlled solutions for the corresponding three-dimensional null controllability problem and derived a small-time global approximate Lagrangian controllability result as well.
Recently, related global boundary controllability results were developed for Boussinesq and micropolar fluids in \cite{CFKMS2023,TYY2025}. Parallel developments include reducing the number of control components \cite{CarG2013,CL2014}, as well as global approximate controllability via degenerate or spatially localized distributed controls on the two-dimensional torus \cite{AS2006,NRB2025, NRNS2025}.

The controllability of coupled chemotaxis-fluid
systems has also attracted attention, as these systems keep both the characteristics of the chemotactic coupling and the nonlinearity of incompressible fluid flow.
In this setting,  Chaves-Silva and Guerrero \cite{CG2017} proved local exact controllability of a Keller-Segel-Navier-Stokes system around trajectories with constant cell density, a spatially homogeneous chemical concentration and zero velocity.
Their controls act on the chemical equation and, in dimension three, on one more velocity component.
The local controllability results for chemotaxis and
chemotaxis--fluid systems above do not provide
a global transfer from arbitrarily large nonnegative
initial data to a prescribed positive trajectory.
This motivates the question of whether such a transfer can be achieved while preserving both biological components throughout the evolution.
In this paper, we aim to establish global exact boundary controllability
in dimensions two and three, using two scalar Neumann
controls and velocity boundary controls on $\Gamma_c$ without requiring the initial state to be
close to the reference.

Indeed, preserving the nonnegativity of the state
variables presents additional challenges in controllability
problems, particularly when global controllability
is considered. Constraints imposed on the state or on the control can lead to a strictly positive minimal time for controllability, even in the case of the linear heat equation 
\cite{LTZ2017}. Related waiting time phenomena and large-time controllability for semilinear heat equations have also been investigated in \cite{PZ2018}.
In our problem,  the states must
satisfy $n,c\ge0$ along the entire controlled trajectory.

Compared with the coupled fluid systems in \cite{CFKMS2023,TYY2025}, the present problem involves second order chemotactic cross-diffusion between two scalar components that must both remain nonnegative. For global controllability, these pointwise constraints must be preserved throughout the control process, even when the initial state is far from the target and intermediate states may be large. Moreover, the interaction of chemotactic cross diffusion with three-dimensional convection and Navier boundary layers further complicates the constrained control of the coupled dynamics.
Further difficulties arise when the reference trajectory
depends on space. If we set
$(r,q,v)=(n-\bar n,c-\bar c,u-\bar u)$, the linearized
system contains
$
\chi\Div(\bar n\nabla q), 
v\cdot\nabla\bar n, 
v\cdot\nabla\bar c.
$
The first term cannot be considered as a lower order perturbation since it has the same differential order as
the scalar diffusion. The last two coupling terms, in turn, feed the velocity error back into the two scalar equations, while the density error enters the fluid momentum equation. Therefore, even at the linearized level, the two scalar components remain strongly coupled to the fluid component.
The difficulty therefore lies not only in achieving exact terminal matching, but also in doing so from large initial states without violating either nonnegativity constraint during the transfer.

\subsection{Notations and reference trajectories}
Fix a bounded connected smooth extended domain $\OO$ with
\begin{equation*}
 \Om\subset\OO,\qquad  \Gamma=\partial\OO,\qquad\Gamma_c\subset\OO,\qquad\Gamma_0\subset\partial\OO,
 \qquad E=\OO\setminus\overline\Om.
\end{equation*}
Reducing the active boundary part if necessary, we take the extension geometry as in Section 2.1 of \cite{CMS2020}, which implies that the associated potential return flow carries each particle into the exterior. We can also choose nonempty open sets $\omega_0\subset\subset\omega\subset\subset E$. For $L>0$, let $Q_L=(0,L)\times\OO$ and $\Sigma_L=(0,L)\times\partial\OO$. The extended friction matrix remains smooth, symmetric and time independent and coincides with the original coefficient on the common boundary.

\paragraph{Notation and state spaces.}
For brevity, we write, for example,
\[
 \|f\|_4=\left(\int_{\OO}|f|^4\dd x\right)^{1/4},\qquad
 \|q\|_{W^{1,4}}=\big(\|q\|_4^4+\|\nabla q\|_4^4\big)^{1/4}
\]
to denote the ordinary Lebesgue and Sobolev norms. 
$\|f\|_{p,\partial\OO}$ is the boundary $L^p$ norm.
For an interval $I$, the shortcut $L^2H^1$ represents $L^2(I;H^1(\OO))$, and similarly for other spaces. Unless specified otherwise, time norms are taken over the interval under consideration, and spatial norms and unmarked integrals are taken over $\OO$. 
We use $(\nabla u)_{ij}=\partial_j u_i$. For a smooth bounded domain $D$ and $1<p<\infty$, put $p'=p/(p-1)$ and define
\begin{equation*}
 \mathcal H_p(D)=\left\{v\in L^p(D)^d:
       \int_D v\cdot\nabla\varphi=0
       \text{ for every }\varphi\in W^{1,p'}(D)\right\}.
\end{equation*}
Let 
$\Hh=\mathcal H_2(\OO)$, $\Vv=H^1(\OO)^d\cap\Hh$, 
$H_N^2(\OO)=\{z\in H^2(\OO):\partial_\nu z=0\}$ and 
$\mathbb P_p$ be the bounded Helmholtz projection onto $\mathcal H_p(\OO)$:
$\mathbb P_pf=f-\nabla\pi_f$, where $\pi_f\in W^{1,p}(\OO)/\mathbb R$
satisfies $\int_\OO\nabla\pi_f\cdot\nabla\varphi
=\int_\OO f\cdot\nabla\varphi$ for every $\varphi\in W^{1,p'}(\OO)$.
We suppress the subscript when its value is determined by the norm, such as
$\mathbb P=\mathbb P_2$ on $L^2$. We denote by $\nabla_\Gamma=P_\tau\nabla$ and $\Div_\Gamma$ the tangential
gradient and surface divergence, respectively. The Stokes operator with Navier boundary conditions is defined by $A_S=-\mu\mathbb P\Delta$ with
\begin{equation*} 
 D(A_S)=\{v\in H^2(\OO)^d\cap\Vv:\NN v=0\}.
\end{equation*}

Let $L^2_\gamma(\Omega)$ be the closure in
$L^2(\Omega)^d$ of smooth divergence-free vector
fields whose normal components vanish on
$\Gamma_0$, 
$
V_\gamma=H^1(\Omega)^d\cap L^2_\gamma(\Omega)
$
and
\begin{equation}\label{1.3}
 \X_0=L^2\times H^1\times\Hh,\qquad
 \Y=H^1\times H_N^2\times\Vv,\qquad
 \F=L^2\times H^1\times L^2(\OO)^d.
\end{equation}
The control functions will be taken in the space
\begin{equation}\label{1.4}
 \Uc=H^1(\OO)\times H^1(\OO)\times L^2(\OO)^d.
\end{equation}
For an interval $I$, the strong state space $\X(I)$ is the product of
\begin{align}
 \X_n(I)&=H^1(I;L^2)\cap L^2(I;H_N^2)\cap C(\overline I;H^1),\nonumber\\
 \X_c(I)&=H^1(I;H^1)\cap L^2(I;H^3)\cap C(\overline I;H_N^2),\label{1.5}\\
 \X_u(I)&=H^1(I;\Hh)\cap L^2(I;D(A_S))\cap C(\overline I;\Vv).\nonumber
\end{align}

Next, we will consider the following reference trajectories.
\begin{definition}\label{definition 1.1}
An admissible reference is the restriction to $\Omega$ of $(\nb,\cb,\ub,\pb)$ on $[0,T]\times\OO$ satisfying
\[
 \nb\in C([0,T];L^4(\OO)),\qquad
 \cb\in C([0,T];W^{1,4}(\OO)),\qquad
 \ub\in C([0,T];L^4(\OO)^d),
\]
\[
 (\nb,\cb,\ub)\in\X(\delta,T),\qquad
 \pb\in L^2(\delta,T;H^1(\OO)/\mathbb R)
 \quad\text{for every }\delta>0,
\]
with $\X$ defined in \eqref{1.5}.
The initial data satisfy
\[
 \bar n_*\in L^4(\OO),\qquad
 \bar c_*\in W^{1,4}(\OO),\qquad
 \bar u_*\in L^4(\OO)^d,\quad
 \Div\bar u_*=0,\quad\bar u_*\cdot\nu=0.
\]
The physical and extended initial data satisfy
\begin{equation*}
\left\{\begin{aligned}
 &(\nb,\cb,\ub)|_{s=0}=(\bar n_*,\bar c_*,\bar u_*), &&\text{in }\OO,\\[2pt]
 &(\bar n_*,\bar c_*,\bar u_*)=(\bar n_0,\bar c_0,\bar u_0), &&\text{in }\Omega,\\[2pt]
 &\bar n_*,\bar c_*\ge0, &&\text{in }\OO.
\end{aligned}\right.
\end{equation*}
The corresponding states satisfy
\begin{equation*}
\left\{\begin{aligned}
 &\Div\ub=0, &&\text{in }Q_T,\\[2pt]
 &\partial_\nu\nb=\partial_\nu\cb=0,\quad\ub\cdot\nu=0,\quad\NN\ub=0, &&\text{on }\Sigma_T.
\end{aligned}\right.
\end{equation*}
Furthermore, for every $\tau\in(0,T)$, there is $m_\tau>0$ such that
\begin{equation}\label{1.6}
 \nb,\cb\ge m_\tau
 \quad\text{a.e. on }[\tau,T]\times\OO.
\end{equation}
Their equation defects
\begin{align*}
 \bar f_n&=\nb_s+\ub\cdot\nabla\nb-D_n\Delta\nb+\chi\Div(\nb\nabla\cb),\\
 \bar f_c&=\cb_s+\ub\cdot\nabla\cb-D_c\Delta\cb+\nb\cb,\\
 \bar f_u&=\ub_s+(\ub\cdot\nabla)\ub-\mu\Delta\ub+\nabla\pb-\nb e_d
\end{align*}
belong to $L^2(0,T;\F)$, vanish in $\Omega$, and are smooth on
$[\tau,T]\times\overline\OO$ for every $\tau>0$. 
\end{definition}

\begin{remark}\label{remark 1.2}
For $N_{\rm r},\mathcal{C}_{\rm r}>0$, we can choose
\[
 \nb=N_{\rm r},\qquad \cb=\mathcal{C}_{\rm r}e^{-N_{\rm r}s},\qquad \ub=0,\qquad \pb=N_{\rm r}x_d.
\]
Then $\bar f_n=\bar f_c=\bar f_u=0$ and
$
 \min_{[0,T]\times\overline\OO}\min\{\nb,\cb\}
 =\min\{N_{\rm r},\mathcal{C}_{\rm r}e^{-N_{\rm r}T}\}>0.
$
Thus Definition~\ref{definition 1.1} is nonempty for every $T>0$.
\end{remark}

Indeed, applying Lemma~\ref{lemma 2.4} to the reference and its sources, and choosing a smooth pressure normalization, we have
\[
 (\nb,\cb,\ub,\pb)\in C^\infty([\tau,T]\times\overline\OO)
 \qquad(0<\tau<T).
\]
In addition, in order to retain the non-negativity of the solution of the biological equations, we need to define the exterior reference gains. That is, for $s>0$,
\begin{equation}\label{1.7}
 \an=\frac{\bar f_n}{\nb} ~~\text{and~~} \ac=\frac{\bar f_c}{\cb}.
\end{equation}

\subsection{Main results}
\begin{theorem} \label{thm 1.3}
Let $d\in\{2,3\}$, $T>0$ and 
$n_0\in L^4(\Omega), c_0\in W^{1,4}(\Omega), 
 u_0\in\mathcal H_4(\Omega)$ with $n_0,c_0\ge0$ be the given initial data.
Assume that the
reference satisfies Definition~\ref{definition 1.1} on $[0,T]$.
There are constants $C,\varepsilon_0>0$ such that, for every
$0<\varepsilon<\varepsilon_0$, there exists a controlled solution of
\eqref{1.1} with $n,c\ge0$ and
\begin{equation}\label{1.8}
 \|n(T)-\nb(T)\|_{L^4(\Omega)}
 +\|c(T)-\cb(T)\|_{W^{1,4}(\Omega)}
 +\|u(T)-\ub(T)\|_{L^4(\Omega)}\le C\varepsilon^{1/8},
\end{equation}
where the constants depend only on the initial data, $T$,  $\beta$
and the reference trajectories.
The state satisfies
\begin{equation}\label{1.9}
\begin{aligned}
 n&\in C([0,T];L^4(\Omega))\cap L^2(0,T;H^1(\Omega)), 
 \qquad  c\in C([0,T];W^{1,4}(\Omega))\cap L^2(0,T;H^2(\Omega)),\\
 u&\in C([0,T];L^4(\Omega)^d)\cap L^2(0,T;H^1(\Omega)^d).
\end{aligned}
\end{equation}
For every $0<\delta<T$, the state also satisfies
\[
\begin{aligned}
 n&\in H^1(\delta,T;L^2(\Omega))\cap L^2(\delta,T;H^2(\Omega)), \qquad
 c\in H^1(\delta,T;H^1(\Omega))\cap L^2(\delta,T;H^3(\Omega)),\\
 u&\in H^1(\delta,T;L^2(\Omega)^d)\cap L^2(\delta,T;H^2(\Omega)^d).
\end{aligned}
\]
The controls satisfy
\begin{equation}\label{1.10}
\begin{aligned}
 h_n^*,h_c^*&\in L^2(\delta,T;H^{1/2}(\Gamma_c)), \qquad
 h_\tau^*&\in L^2(\delta,T;H^{1/2}(\Gamma_c)^d),\qquad
 h_\nu^*\in L^2(\delta,T;H^{3/2}(\Gamma_c)),
\end{aligned}
\end{equation}
with $h_\tau^*\cdot\nu=0$ and \eqref{1.2}.
\end{theorem}

\begin{theorem}\label{thm 1.4}
Under the assumptions of Theorem~\ref{thm 1.3}, there are boundary
controls and a solution of \eqref{1.1} with $n,c\ge0$, with the same state and control regularity as in Theorem~\ref{thm 1.3}, such that, for some
$s_3\in(0,T)$,
\begin{equation}\label{1.11}
 (n,c,u)(s,\cdot)=(\nb,\cb,\ub)(s,\cdot)
 \quad\text{in }\Omega,\qquad s_3\le s\le T.
\end{equation}
\end{theorem}

\begin{remark} \label{remark 1.5}
For a regular solution of \eqref{1.1} with $n,c\ge0$, comparison on
interior positive time cylinders and the strong maximum principle give
\[
 z_0\not\equiv0\quad\Longrightarrow\quad
 z(T,x)>0\quad(x\in\Omega),\qquad z\in\{n,c\}.
\]
Here the regularity is only needed on interior cylinders.
Thus a nontrivial nonnegative component cannot
be driven exactly to zero in finite time, irrespective of the signs of
its boundary control.  
\end{remark}

\begin{remark} Our results include bounded
$n_0,u_0$ and Lipschitz $c_0$, subject to $\Div u_0=0$ and the full zero
normal trace. They also include $n_0\in H^1$, $c_0\in H^2$ and
$u_0\in H^1\cap\mathcal H_4$, without imposing $\partial_\nu c_0=0$.
By Corollary~\ref{corollary 2.6}, admissible references can also be generated locally from the same kinds of initial data:
\begin{equation}\label{1.12}
 \bar n_0\in L^4(\Omega),\qquad \bar c_0\in W^{1,4}(\Omega),\qquad
 \bar u_0\in\mathcal H_4(\Omega),\qquad \bar n_0,\bar c_0\ge0.
\end{equation}
\end{remark}

\begin{remark} 
\label{remark 1.7}
The space $L^4\times W^{1,4}\times L^4$ is chosen
to treat dimensions two and three within a common
framework. This framework is
 compatible with the weak boundary-layer
estimates and the strong continuation argument.
However, for the two-dimensional case, the global controllability
can also be obtained through an $L^2$-based
remainder argument, followed by regularization
into $\mathcal Y$ and local exact control.
\end{remark}

\begin{remark} The global controllability results hold for every given time
$T>0$. Indeed, all stages of the construction,
including the initial regularization, positive
preparation, global approximation and local exact
correction, can be completed within an arbitrarily
small positive time. The parameter $\varepsilon$ can be chosen sufficiently
small, depending on $T$, the initial data and the
reference trajectory, such that the 
approximation error satisfies both the regularization
and local controllability thresholds.
\end{remark}

\paragraph{Main ideas of the proof.}
In the proofs of Theorems~\ref{thm 1.3} and
\ref{thm 1.4}, we combine a global approximation
with a local exact correction.
The initial data can be far from the reference,
and both biological components  may vanish initially.
The main aim here is to construct a controlled
state that respects the state constraints and
leads to an error sufficiently small for local
controllability. The argument proceeds through
positive preparation, global approximation,
regularization of the error and local exact control.

First, the regularity of the initial data is improved
by short-time regularization, and then a positivity
preparation is carried out such that both biological
scalar components become strictly positive.
Subsequently, along the return flow, we construct
the scalar profiles by interpolating between the
prepared data and the reference state, allowing
the exterior scalar sources to be chosen in
multiplicative form. Thus, nonnegativity can be
preserved by the actual nonlinear equations.
The joint remainder estimates in the
$L^4\times W^{1,4}\times L^4$ framework justify
the above large-amplitude transfer construction
and, together with a strong continuation criterion,
yield a solution throughout the control interval
and a terminal error of order $\varepsilon^{1/8}$.
These estimates treat the chemoattractant,
three-dimensional convection and weak boundary-layer
defects within a common framework, while retaining
the Navier pressure and the inhomogeneous boundary
stress.

Next, we introduce a short regularization stage
with the reference gains since the above approximation
does not directly satisfy the smallness condition
in $\mathcal Y$ required for local exact controllability.
This regularization preserves nonnegativity and
gives an error sufficiently small in $\mathcal Y$
at a suitable time.

Finally, local exact controllability near spatially
nonconstant references is obtained through a joint
scalar--Navier Carleman estimate and the fixed-point
method.
% Added
An additional time weight controls the second-order
chemotactic coupling, while weighted Stokes regularity
handles the pressure terms. The nonlinear argument
also retains the control state products required
by the multiplicative formulation.
Combining the global approximation and local exact
correction yields the desired controllability while
preserving both biological state constraints.

\section{The extended system}\label{section 2}
On the extended domain $\OO$, recall that $
 Q_T=(0,T)\times\mathcal O, 
 \Sigma_T=(0,T)\times\partial\mathcal O.
$ Now, we consider
\begin{equation}\label{2.1}
\left\{\begin{aligned}
 &n_s+u\cdot\nabla n-D_n\Delta n=-\chi\Div(n\nabla c)+f_n, &&\text{in }Q_T,\\[2pt]
 &c_s+u\cdot\nabla c-D_c\Delta c=-nc+f_c, &&\text{in }Q_T,\\[2pt]
 &u_s+(u\cdot\nabla)u-\mu\Delta u+\nabla p=ne_d+f_u, 
 ~~\Div u=\sigma, &&\text{in }Q_T,\\[2pt]
 &\partial_\nu n=\partial_\nu c=0,\quad u\cdot\nu=0,\quad\NN u=0, &&\text{on }\Sigma_T,\\[2pt]
 &(n,c,u)|_{s=0}=(n_*,c_*,u_*), &&\text{in }\OO.
\end{aligned}\right.
\end{equation}
All four sources $(f_n, f_c, f_u, \sigma)$ are internal controls, which vanish in $\Om$, and $\int_\OO\sigma=0$.  

 \subsection{Nonnegative extension and preservation}\label{subsection 2.1}

We first extend the physical initial data without imposing initial Neumann
or Navier stress compatibility.

\begin{proposition} \label{prop 2.1}
For any initial data satisfying $n_0\in L^4(\Omega), c_0\in W^{1,4}(\Omega), 
 u_0\in\mathcal H_4(\Omega)$ with $n_0,c_0\ge0$, there exists
an extension such that
\begin{equation}\label{2.2}
\left\{\begin{aligned}
 &(n_*,c_*,u_*)=(n_0,c_0,u_0),&&\text{in }\Omega,\\
 &n_*\in L^4(\mathcal O),\quad
 c_*\in W^{1,4}(\mathcal O),\quad u_*\in\mathcal H_4(\mathcal O),\\
 &n_*,c_*\geq0,&&\text{a.e. in }\mathcal O,
\end{aligned}\right.
\end{equation}
and
\begin{equation}\label{2.3}
 \|n_*\|_4+\|c_*\|_{W^{1,4}}+\|u_*\|_4
 \leq C\bigl(\|n_0\|_{L^4(\Omega)}
       +\|c_0\|_{W^{1,4}(\Omega)}+\|u_0\|_{L^4(\Omega)}\bigr).
\end{equation}
Moreover, the total masses of the density and chemical
extensions can also be made strictly positive by adding a fixed nonzero function
$0\le \zeta\in C_c^\infty(E)$ to both extensions.
In this case, the right-hand side of \eqref{2.3} is increased by
$\|\zeta\|_4+\|\zeta\|_{W^{1,4}}$.
\end{proposition}
\begin{proof}
Let $E_1:W^{1,4}(\Omega)\to W^{1,4}(\mathbb R^d)$ be a bounded Sobolev
extension operator. Set
\[
 n_*=\mathbf1_\Omega n_0,\qquad
 c_*=(E_1c_0)^+|_{\mathcal O},\qquad
 u_*=\mathbf1_\Omega u_0.
\]
It holds that
$(E_1c_0)^+|_\Omega=c_0$, and $\nabla(E_1c_0)^+=\mathbf1_{\{E_1c_0>0\}}\nabla E_1c_0,$
which implies that $c_*\geq0$ and
$\|c_*\|_{W^{1,4}}\leq C\|c_0\|_{W^{1,4}(\Omega)}$.
Moreover,
\[
 \|n_*\|_4=\|n_0\|_{L^4(\Omega)},\qquad
 \|u_*\|_4=\|u_0\|_{L^4(\Omega)}.
\]
For every $\varphi\in W^{1,4/3}(\mathcal O)$, by the definition of
$\mathcal H_4(\Omega)$, one can see
\[
 \int_{\mathcal O}u_*\cdot\nabla\varphi\,dx
 =\int_\Omega u_0\cdot\nabla(\varphi|_\Omega)\,dx=0.
\]
Thus $u_*\in\mathcal H_4(\mathcal O)$, which implies
\eqref{2.2}-\eqref{2.3}.
For $0\leq\zeta\in C_c^\infty(E)$ with $\int_E\zeta>0$, replace $n_*,c_*$ by
$n_*+\zeta,c_*+\zeta$. Their restrictions and signs are unchanged, and the
extended masses are positive. The new bound is the right-hand side of \eqref{2.3} plus $\|\zeta\|_4+\|\zeta\|_{W^{1,4}}$.
\end{proof}

We now show that suitable exterior controls preserve the nonnegativity throughout the extended domain.
 
\begin{lemma}\label{lemma 2.2}
Let $(n,c,u)$ be a strong solution of \eqref{2.1} on $[0,L]$ satisfying
$\operatorname{div}u=\sigma\in L^1(0,L;L^\infty)$ and the homogeneous outer
boundary conditions, with $n_*, c_*\ge0$. Assume
\begin{equation*}
 f_n=\alpha_n n+g_n,\qquad f_c=\alpha_c c+g_c,\qquad g_n,g_c\ge0,
\end{equation*}
where $\alpha_n,\alpha_c\in L^2(0,L;L^3)$ and $\nabla c\in L^4(0,L;L^6)$. Then, the biological states satisfy 
\[
n\geq 0, ~~c\geq 0, ~~~~a.e. ~x\in \OO, ~~0\leq s \leq L.
\]
\end{lemma}
\begin{proof}
Let $z=n^-:=\max\{-n,0\}$. Multiplying the density equation in \eqref{2.1} by $-z$ and using integration by parts, one has
\begin{equation*}
 \frac12\frac{\dd}{\dd s}\|z\|_2^2+D_n\|\nabla z\|_2^2
 =\frac12\int\sigma z^2+\int \alpha_nz^2
       +\chi\int z\nabla c\cdot\nabla z-\int g_nz.
\end{equation*}
By interpolation, we get
\[
 \|z\|_3\le \|z\|_2^{1/2}\|z\|_6^{1/2}
 \le C\|z\|_2^{1/2}(\|\nabla z\|_2+\|z\|_2)^{1/2}.
\]
For  $\eta>0$, we have
\begin{align*}
 \int|\alpha_n|z^2
 &\le\|\alpha_n\|_3\|z\|_3^2
 \le C\|\alpha_n\|_3\|z\|_2(\|\nabla z\|_2+\|z\|_2) 
 \le\eta\|\nabla z\|_2^2+C_\eta(1+\|\alpha_n\|_3^2)\|z\|_2^2,\\
 \chi\int|z||\nabla c||\nabla z|
 &\le C\|\nabla c\|_6\|z\|_2^{1/2}\|\nabla z\|_2^{3/2}
       +C\|\nabla c\|_6\|z\|_2\|\nabla z\|_2 
 \le\eta\|\nabla z\|_2^2+C_\eta(1+\|\nabla c\|_6^4)\|z\|_2^2.
\end{align*}
Taking $2\eta<D_n$ and using $g_nz\ge0$, we obtain
\[
\frac{\dd}{\dd s}\|z\|_2^2\le C\bigl(1+\|\sigma\|_\infty+\|\alpha_n\|_3^2+\|\nabla c\|_6^4\bigr)\|z\|_2^2,
\]
which, together with Gronwall's inequality, implies $n\ge0$. Similarly, we can also obtain $c\ge0$ by using $-\int_{\mathcal O}n(c^-)^2\leq0$.
\end{proof}

\subsection{Regularization and positive preparation}\label{subsection 2.2}

In this subsection, we first regularize the extended data and then make the 
controlled scalar states strictly positive. For this purpose, we need to use semigroup estimates.

For $v,w\in\mathcal V$, define the symmetric form 
\begin{equation*}
 a_\beta(v,w)=2\mu\int_{\mathcal O}D(v):D(w)\,dx
       +\int_{\partial\mathcal O}(\beta v_{\tan})\cdot w_{\tan}\,dS.
\end{equation*}
The trace and Korn inequalities give positive constants $\gamma_{\rm E}$ and $C$ such that
\begin{equation}\label{2.4}
 a_\beta(v,v)\geq \gamma_{\rm E}\|v\|_{H^1}^2-C\|v\|_2^2,
 \qquad v\in\mathcal V.
\end{equation}
Positive lower bounds in subsequent energy estimates may be decreased as needed. Thus $a_\beta+\lambda(\cdot,\cdot)$ is coercive on $\mathcal V$ for all
sufficiently large $\lambda>0$. Then, let $\mathbb P_p$ with $1<p<\infty$ be the Helmholtz projection and
\begin{equation*}
\begin{aligned}
 A_{p,\beta}v=-\mu\mathbb P_p\Delta v, \quad
 D(A_{p,\beta})=\{v\in W^{2,p}(\mathcal O)^d\cap\mathcal H_p(\mathcal O):
                         \mathcal N_\beta v=0\text{ on }\partial\mathcal O\}.
\end{aligned}
\end{equation*}

  By virtue of Kunstmann's results on the bounded $H^\infty$-calculus,
after a spectral shift, and on fractional domains of the Stokes
operator with Robin and Navier boundary conditions
\cite{K2025}, together with Sobolev interpolation and the standard
Neumann heat semigroup estimates, we have the bounds below.
For nonnegative scalar friction, see \cite{AEG2018}.

\begin{proposition}\label{prop 2.3}
For each $1<p<\infty$, the operator $-A_{p,\beta}$ generates an
analytic $C_0$-semigroup on $\mathcal H_p(\mathcal O)$, 
which is denoted by $S_u(t)=S_\beta(t)$.
Let
\[
 S_n(t)=e^{D_nt\Delta_N},\qquad S_c(t)=e^{D_ct\Delta_N},
\]
where $\Delta_N$ is the Laplacian with homogeneous Neumann boundary conditions.
For $1<r\leq q<\infty$ and 
$
 f,g\in L^r(\mathcal O), v\in\mathcal H_r(\mathcal O), 
 F\in L^r(\mathcal O)^d,  G\in L^r(\mathcal O)^{d\times d},
$
and $0<t\leq1$, one has
\begin{equation}\label{2.5}
\begin{aligned}
 & \|S_u(t)v\|_q+t^{1/2}\|\nabla S_u(t)v\|_q
\leq Ct^{-\gamma_{r,q}}\|v\|_r,\\
 & \|S_j(t)f\|_q+t^{1/2}\|\nabla S_j(t)f\|_q
\leq Ct^{-\gamma_{r,q}}\|f\|_r,\qquad j\in\{n,c\},\\
 &\|S_n(t)\operatorname{div}_N F\|_q
\leq Ct^{-1/2-\gamma_{r,q}}\|F\|_r,\\
& \|S_u(t)\mathbb P\operatorname{div}G\|_q
\leq Ct^{-1/2-\gamma_{r,q}}\|G\|_r,\\
& \|S_c(t)g\|_{W^{1,q}}
\leq Ct^{-1/2-\gamma_{r,q}}\|g\|_r,
\end{aligned}
\end{equation}
where  
$\gamma_{r,q}=\frac d2(\frac1r-\frac1q)$ and the constant $C$ depends on $\mathcal O,d,\mu,D_n,D_c,\beta,r,q$,
but not on $t$. Here the weak divergence operators are defined by
\[
 \langle \operatorname{div}_N F,\varphi\rangle
 =-\int_{\mathcal O}F\cdot\nabla\varphi\,dx,
 \qquad \varphi\in W^{1,r'}(\mathcal O),
\]
and
\[
 \langle \mathbb P\operatorname{div}G,\psi\rangle
 =-\int_{\mathcal O}G:\nabla\psi\,dx,
\]
for every $\psi\in W^{1,r'}(\mathcal O)^d$ satisfying
$\operatorname{div}\psi=0$ in $\mathcal O$ and
$\psi\cdot\nu=0$ on $\partial\mathcal O$,
where $r'=r/(r-1)$.
 In particular, for $r=2$, $q=4$ and
$\theta_d=\frac12+\frac d8<1$ with $d\in\{2,3\}$,
\[
\begin{aligned}
 \|S_n(t)\operatorname{div}_N F\|_4
   +\|S_u(t)\mathbb P\operatorname{div}G\|_4
      \leq Ct^{-\theta_d}(\|F\|_2+\|G\|_2), ~~~~
 \|S_c(t)g\|_{W^{1,4}}
     \leq Ct^{-\theta_d}\|g\|_2.
\end{aligned}
\]
Moreover, $S_n,S_c,S_u$ are strongly continuous on
$L^4(\mathcal O)$, $W^{1,4}(\mathcal O)$ and
$\mathcal H_4(\mathcal O)$, respectively, and
\[
 \sup_{0\leq t\leq1}
 \bigl(\|S_n(t)f\|_4+\|S_c(t)g\|_{W^{1,4}}+\|S_u(t)v\|_4\bigr)
 \leq C\bigl(\|f\|_4+\|g\|_{W^{1,4}}+\|v\|_4\bigr)
\]
for $f\in L^4(\mathcal O)$, $g\in W^{1,4}(\mathcal O)$ and
$v\in\mathcal H_4(\mathcal O)$.
\end{proposition}

\begin{proof}
Let $M_w$ be the shape operator of $\partial\mathcal O$ with 
$M_w\eta=\nabla_\eta\nu$ for tangent vectors $\eta$ and $M_w\nu=0$.
Thus, the Navier
boundary condition is equivalent to
\[
 v\cdot\nu=0,\qquad
 [ (\nabla v)^T-\nabla v ]\nu=B_\beta v,\qquad
 B_\beta=\mu^{-1}\beta_\tau-2M_w,
\]
where the matrix $B_\beta$ is smooth and symmetric.
As $\mathcal O$ is smooth and bounded, the Helmholtz decomposition
holds for every $1<p<\infty$. By Theorem~5.5 and Corollary~5.6 in \cite{K2025},  
after a sufficiently large spectral shift, we obtain an analytic
semigroup and the square-root domain
$
 D((\lambda+A_{p,\beta})^{1/2})
   =W^{1,p}(\mathcal O)^d\cap\mathcal H_p(\mathcal O).
$
Therefore, $\|\nabla S_u(t)w\|_p\leq Ct^{-1/2}\|w\|_p$
for $w\in\mathcal H_p(\mathcal O)$ and $0<t\leq1$.
According to analyticity, elliptic regularity and Sobolev interpolation, we obtain
$L^r$-$L^q$ smoothing for nearby exponents and the semigroup property
extends it to all $1<r\leq q<\infty$. The first bound in
\eqref{2.5} follows by applying
 the $L^q$ gradient estimate to
$S_u(t/2)S_u(t/2)v$.
The other bounds follow from the Neumann heat semigroup
estimates. Since $t\leq1$, we also have the $W^{1,q}$ bound.

The symmetry of $\beta$ makes the Stokes semigroups at conjugate
exponents dual to each other, which implies for
$\varphi\in L^{q'}(\mathcal O)$ and
$\psi\in\mathcal H_{q'}(\mathcal O)$,
\[
\begin{aligned}
 |\langle S_n(t)\operatorname{div}_N F,\varphi\rangle|
  &=\left|\int_{\mathcal O}F\cdot\nabla S_n(t)\varphi\,dx\right|
    \leq Ct^{-1/2-\gamma_{r,q}}\|F\|_r\|\varphi\|_{q'},\\
 |\langle S_u(t)\mathbb P\operatorname{div}G,\psi\rangle|
  &=\left|\int_{\mathcal O}G:\nabla S_u(t)\psi\,dx\right|
    \leq Ct^{-1/2-\gamma_{r,q}}\|G\|_r\|\psi\|_{q'},
\end{aligned}
\]
provided $q'\leq r'$ and $\gamma_{q',r'}=\gamma_{r,q}$.
The two divergence estimates follow from duality and the boundedness of the Helmholtz projection.

Finally,
$D((I-D_c\Delta_N)^{1/2})=W^{1,4}(\mathcal O)$ on $L^4(\mathcal O)$ implies $S_c$ is strongly continuous and bounded for
$0\leq t\leq1$.
$S_n$ and $S_u$ are strongly continuous on
$L^4(\mathcal O)$ and $\mathcal H_4(\mathcal O)$, 
and are therefore uniformly bounded.
Since $e^{\lambda t}$ is bounded on $0\leq t\leq1$,  
the factor introduced by the spectral shift can be absorbed
into the constant $C$. The proof of this proposition is completed.
\end{proof}

Next, we establish local existence and regularization.

\begin{lemma} \label{lemma 2.4}
Let $T>0$ and
$
 n_*\in L^4(\mathcal O), c_*\in W^{1,4}(\mathcal O), 
 u_*\in\mathcal H_4(\mathcal O), n_*,c_*\geq0.
$
There is $0<T_r<T/16$ such that \eqref{2.1}, with
$f_n=f_c=f_u=\sigma=0$, has a unique mild solution satisfying
\begin{equation}\label{2.6}
\begin{aligned}
 n&\in C([0,2T_r];L^4)\cap L^2(0,2T_r;H^1), 
 ~~c\in C([0,2T_r];W^{1,4})\cap L^2(0,2T_r;H^2),\\
 u&\in C([0,2T_r];\mathcal H_4)\cap L^2(0,2T_r;H^1)
\end{aligned}
\end{equation}
with $n, c\geq0$. For any \(0<\delta<2T_r\), the solution is smooth and belongs to \(\X(\delta,2T_r)\) defined in \eqref{1.5}. 

Moreover, the same regularization property holds for any solution with \(\sigma=0\), homogeneous boundary conditions and smooth prescribed sources. More precisely, if \((n,c,u)\in\X(\delta,L)\) for every \(0<\delta<L\), then, for every integer \(m\geq1\),
\begin{equation*}
\sup_{\delta\leq s\leq L}
\left(
\|n(s)\|_{H^m}
+\|c(s)\|_{H^{m+1}}
+\|u(s)\|_{H^m}
\right)<\infty.
\end{equation*}
In particular, the unforced solution at \(T_r\) is smooth and satisfies all boundary conditions.
\end{lemma}
\begin{proof}
Using Proposition~\ref{prop 2.3} and setting 
$\theta=1-\theta_d=\frac12-\frac d8>0$, for $0\leq t\leq\ell\leq1$, one has 
\begin{equation}\label{2.7}
\left\{\begin{aligned}
 n(t)&=S_n(t)n_*
       -\int_0^t S_n(t-s)\operatorname{div}_N(un+\chi n\nabla c)(s)\,ds:=\Phi_n,\\
 c(t)&=S_c(t)c_*
       -\int_0^t S_c(t-s)(u\cdot\nabla c+nc)(s)\,ds:=\Phi_c,\\
 u(t)&=S_u(t)u_*
       -\int_0^t S_u(t-s)\mathbb P\operatorname{div}(u\otimes u)(s)\,ds
     +\int_0^t S_u(t-s)\mathbb P(ne_d)(s)\,ds:=\Phi_u.
\end{aligned}\right.
\end{equation}
Since $W^{1,4}\hookrightarrow L^\infty$ for $d\leq3$, we have
\[
\begin{aligned}
 \|un\|_2+\|n\nabla c\|_2
    &\leq\|n\|_4(\|u\|_4+\|\nabla c\|_4),\\
 \|u\cdot\nabla c\|_2+\|nc\|_2
    &\leq\|u\|_4\|\nabla c\|_4+C\|n\|_4\|c\|_{W^{1,4}},\\
 \|u\otimes u\|_2&\leq\|u\|_4^2.
\end{aligned}
\]
Moreover, 
$
 \int_0^t(t-s)^{-\theta_d}\,ds=\frac{t^\theta}{\theta}, 
 \theta\geq\frac18.
$
If we take
$\sup_{0\leq t\leq\ell}(\|n(t)\|_4+\|c(t)\|_{W^{1,4}}+\|u(t)\|_4)\leq R$, then
\begin{equation*}
\begin{aligned}
 \sup_{0\leq t\leq\ell}
   (\|\Phi_n(t)\|_4+\|\Phi_c(t)\|_{W^{1,4}}+\|\Phi_u(t)\|_4)
 &\leq C(\|n_*\|_4+\|c_*\|_{W^{1,4}}+\|u_*\|_4) +C\ell^\theta(R^2+R).
\end{aligned}
\end{equation*}
For any two different elements in this ball, one can find
\begin{equation*}
\begin{aligned}
 &\sup_{0\leq t\leq\ell}
 (\|\Phi_n-\widetilde\Phi_n\|_4
 +\|\Phi_c-\widetilde\Phi_c\|_{W^{1,4}}
 +\|\Phi_u-\widetilde\Phi_u\|_4)(t)\\
 &\qquad\leq C\ell^\theta(2R+1)
 \sup_{0\leq t\leq\ell}
 (\|n-\widetilde n\|_4+\|c-\widetilde c\|_{W^{1,4}}
                              +\|u-\widetilde u\|_4)(t).
\end{aligned}
\end{equation*}
Thus, we choose
\[
 R=2C(\|n_*\|_4+\|c_*\|_{W^{1,4}}+\|u_*\|_4)+1,
 \qquad C\ell^\theta(2R+1)<\frac14.
\]
The integral map $\Phi$ maps the closed ball
$\mathcal B_R$ into itself and is a strict contraction. Hence, the Banach fixed point
theorem yields a unique mild solution.

For the regularity, restart \eqref{2.7} on nested time intervals
away from zero. With $L^{p/2}$ products and $L^q$ states, the singularity in
\eqref{2.5} is
\[
 (t-s)^{-\kappa(p,q)},\qquad \kappa(p,q)=\frac12+\frac d2\left(\frac2p-\frac1q\right).
\]
Thus, let us select $(p,q)=(4,5),(5,8),(8,16)$ to obtain
$
 \kappa(4,5)\leq\frac{19}{20}, 
 \kappa(5,8)\leq\frac{73}{80}, 
 \kappa(8,16)\leq\frac{25}{32},
$
which implies
\[
 n,u,\nabla c\in L^\infty L^5,~~
  n,u,\nabla c\in L^\infty L^8,~~
 n,u,\nabla c\in L^\infty L^{16} \text{~~and~~} un, ~n\nabla c, ~u\otimes u \in L^\infty L^{8} \subset L^{8}L^{8}. 
\]
Then, by appealing to parabolic regularity and embedding theorems, we obtain for some $\alpha\in(0,1)$
\[
 n,u,\nabla c\in C^{\alpha,\alpha/2},\qquad
 c\in C^{2+\alpha,1+\alpha/2},\qquad
 n,u\in C^{1+\alpha,(1+\alpha)/2}.
\]
Using the equation and boundary regularity, we finally obtain arbitrary finite order Sobolev regularity.

To establish nonnegativity, one may first approximate the initial data by means of the Neumann and Stokes semigroups. The resulting approximations are smooth, uniformly bounded in the corresponding initial spaces, and converge to the original initial data. Moreover, the scalar approximations remain nonnegative. The contraction estimates above give a common existence interval and a uniform $L^4\times W^{1,4}\times L^4$ bound for these approximations. Local strong solvability and the continuation criterion in Lemma~\ref{lemma 3.3} show that their strong solutions exist on this common interval and coincide with the mild solutions by uniqueness. Lemma~\ref{lemma 2.2} then gives nonnegativity for every approximation. Continuous dependence in $C([0,\ell];L^4\times W^{1,4}\times L^4)$ yields $n,c\geq0$ in the limit.

For the smooth approximations, standard energy estimates, together with the uniform \(L^4\times W^{1,4}\times L^4\) bound and \eqref{2.4}, yield
\[
\int_0^\ell
\left(
\|\nabla n\|_2^2+
\|c\|_{H^2}^2+
\|u\|_{H^1}^2
\right)\,dt\leq C,
\]
which gives the \(L^2\)-in-time regularity in the lemma.
\end{proof}

Then, for constructing regular multiplicative controls later, we use nonnegative exterior sources to prepare strictly positive scalar states of the controlled system.

\begin{lemma} \label{lemma 2.5}
Let $(n_*,c_*,u_*)$ be smooth on $\overline{\mathcal O}$, with $n_*,c_*\geq0$,
$\operatorname{div}u_*=0$, and the boundary conditions in
\eqref{2.1}. For every $T>0$, there exist $0<T_p<T/8$ and a strong
exterior-controlled solution on $[0,T_p]$ such that scalar components remain
nonnegative and satisfy
\begin{equation*}
 \min_{\overline{\mathcal O}}n(T_p)>0,
 \qquad \min_{\overline{\mathcal O}}c(T_p)>0.
\end{equation*}
\end{lemma}
\begin{proof}
Choose $0\leq\zeta_0\in C_c^\infty(\omega)$ with
$I_0:=\int_{\mathcal O}\zeta_0>0$, and
$\theta\in C_c^\infty((1/8,3/8))$ with $0\leq\theta\leq1$ and
$\theta=1$ on $[3/16,5/16]$. Take
\[
 f_n=f_c=f=\theta(s/T_p)\zeta_0,\qquad f_u=\sigma=0.
\]
The local construction of Lemma~\ref{lemma 3.3} applies with $\varepsilon=1$, zero
profiles, $R_n=R_c=f$ and $F_u=h=0$. For $L\leq1$,
\[
 \|f\|_{L^2(0,L;L^2)}+\|f\|_{L^2(0,L;H^1)}
 \leq L^{1/2}(\|\zeta_0\|_2+\|\zeta_0\|_{H^1}).
\]
By the local strong solvability argument in Lemma~\ref{lemma 3.3}, there exists \(L_*>0\), such that the solution exists on \([0,L_*]\). In particular, \(L_*\) can be chosen independently of \(T_p\), since no time derivative of \(\theta(s/T_p)\) is involved. Taking \(T_p=\min\{T/16,L_*/4\}\) ensures existence on \([0,2T_p]\).
By Lemmas~\ref{lemma 2.2} and~\ref{lemma 2.4}, $n,c\geq0$, and the solution is smooth for positive times.

Let $t_0=T_p/2$, $M_z(s)=\int_{\mathcal O}z(s)$ for $z=n,c$ and
$K(s)=\int_0^s\|n(t)\|_\infty\,dt<\infty$. By integration of the scalar
equations, one has
\[
 M_n'=I_0\theta(s/T_p),\qquad
 M_c'+\|n(s)\|_\infty M_c\geq I_0\theta(s/T_p).
\]
Thus, one can see the mass
\[
 M_n(t_0)\geq\frac{I_0T_p}{8}>0,\qquad
 M_c(t_0)\geq e^{-K(t_0)}\frac{I_0T_p}{8}>0.
\]
On $[t_0,T_p]$, $f_n=f_c=0$, and each $z=n,c$ satisfies
\begin{equation}\label{2.8}
\left\{\begin{aligned}
 &z_s-\kappa\Delta z+b\cdot\nabla z+qz=0,
             &&\text{in }(t_0,T_p)\times\mathcal O,\\
 &\partial_\nu z=0,
             &&\text{on }(t_0,T_p)\times\partial\mathcal O,
\end{aligned}\right.
\end{equation}
where
\[
 (\kappa,b,q)=
 \begin{cases}
  (D_n,u+\chi\nabla c,\chi\Delta c),&z=n,\\
  (D_c,u,n),&z=c.
 \end{cases}
\]
Then, replacing $z$ by
$e^{-\lambda(s-t_0)}z$, with $\lambda\geq\|q^-\|_\infty$, we can make the zero-order
coefficient to be nonnegative.
By the strong maximum principle and the parabolic boundary point lemma
\cite{ProtterWeinberger1984}, we have
\[
 M_z(t_0)>0\quad\Longrightarrow\quad z(T_p,x)>0
 \quad(x\in\mathcal O).
\]
If $z(T_p,x_0)=0$ for $x_0\in\partial\mathcal O$, the boundary point lemma
gives $\partial_\nu z(T_p,x_0)<0$, contrary to
\eqref{2.8}. Hence, by compactness,
\[
 \min_{\overline{\mathcal O}}z(T_p)>0,\qquad z=n,c.
\] 
\end{proof}

\subsection{Reference trajectories and exterior gains}\label{subsection 2.3}

Now, let us turn to the reference trajectory. Different from the controlled state above, the construction below realizes the reference class in Definition~\ref{definition 1.1} from physical initial data.

\begin{corollary} \label{corollary 2.6}
Assume \eqref{1.12} hold.  Then there exists a nonnegative extension \((\bar n_*,\bar c_*,\bar u_*)\) satisfying
\[
\int_{\mathcal O}\bar n_*\,dx>0,\qquad
\int_{\mathcal O}\bar c_*\,dx>0.
\]
Let $T_{\max}>0$ be the maximal existence time of the corresponding
unforced mild solution in $L^4\times W^{1,4}\times\mathcal H_4$.
The restriction of this solution to
$\Omega$ satisfies Definition \ref*{definition 1.1} on every interval
$[0,T]$ with $T<T_{\max}$. For $0<\tau<T<T_{\max}$, it holds that
\[
(\bar n,\bar c,\bar u,\bar p)\in
C^\infty([\tau,T]\times\overline{\mathcal O}), ~~\min_{[\tau,T]\times\overline{\mathcal O}}
\{\bar n,\bar c\}>0.
\]
\end{corollary}
\begin{proof}
It is clear that the required extension, local existence, nonnegativity and positive time regularity follow from Proposition~\ref{prop 2.1} and Lemma \ref{lemma 2.4}. The density mass is conserved, while
\[
\frac{d}{dt}\int_{\mathcal O}\bar c\,dx
\geq-\|\bar n(t)\|_\infty\int_{\mathcal O}\bar c\,dx.
\]
For any $t>0$, continuity at zero gives a time $0<\tau_0<t$ with
$\int\bar c(\tau_0)>0$. On $[\tau_0,t]$, the coefficient
$\|\bar n\|_\infty$ is bounded, and therefore
\[
 \int_\OO\bar c(t)\ge
 \exp\left(-\int_{\tau_0}^t\|\bar n(s)\|_\infty\,ds\right)
 \int_\OO\bar c(\tau_0)>0.
\]
The strong maximum principle and the boundary point argument in
Lemma~\ref{lemma 2.5} give strict positivity up to the boundary.
Compactness gives a positive lower bound on every $[\tau,T]$ with
$0<\tau<T<T_{\max}$.  
\end{proof}

For a general admissible reference, recalling 
\eqref{1.7}, one has
\[
 \bar\alpha_n=\frac{\bar f_n}{\bar n},\qquad
 \bar\alpha_c=\frac{\bar f_c}{\bar c}\qquad(s>0).
\]
Fix $s_*\in(0,T)$ and take the lower bound $m_{s_*}$ from
\eqref{1.6}. By Lemma~\ref{lemma 2.4}, the reference and its
sources are smooth on $[s_*,T]$ with
$\bar n,\bar c\geq m_{s_*}>0$. 
Therefore, for each integer $k\geq0$, we have
\begin{equation*}
 \|\bar\alpha_n\|_{C^k([s_*,T]\times\overline{\mathcal O})}
 +\|\bar\alpha_c\|_{C^k([s_*,T]\times\overline{\mathcal O})}
 \leq C_{k,s_*}.
\end{equation*}

\section{Global approximation in a special case without fluid boundary layer}\label{section 3}
In this section and the next section, \((n_{\rm in},c_{\rm in},u_{\rm in})\) denotes the smooth, compatible initial state for the approximation step on \(\overline{\mathcal O}\) with \(n_{\rm in},c_{\rm in}>0\). From now on, we set \(s_0\geq s_*>0\) and \(T_{\rm app}>0\) such that \(s_0+T_{\rm app}\leq T\). The reference trajectory is smooth and strictly positive on \([s_0,s_0+T_{\rm app}]\), as required in Definition~\ref{definition 1.1}.

\subsection{Time scaling and the return flow construction}
We introduce a positive small parameter $0<\varepsilon<\min\{1,T_{\rm app}/2\}$ and set
\begin{equation*}
 s=s_0+\eps t,\quad S_\eps=T_{\rm app}/\eps,\quad
 N=\eps n(s),\quad \mathcal{C}=\eps c(s),\quad U=\eps u(s),\quad P=\eps^2p(s).
\end{equation*}
The scaling on $Q^\varepsilon=(0,S_\varepsilon)\times\OO$, with boundary $\Sigma^\varepsilon=(0,S_\varepsilon)\times\partial\OO$ can be read as
\begin{equation}\label{3.1}
\left\{\begin{aligned}
 &N_t+U\cdot\nabla N-\eps D_n\Delta N=-\chi\Div(N\nabla \mathcal{C})+G_n, &&\text{in }Q^\varepsilon,\\[2pt]
 &\mathcal{C}_t+U\cdot\nabla \mathcal{C}-\eps D_c\Delta \mathcal{C}=-N\mathcal{C}+G_c, &&\text{in }Q^\varepsilon,\\[2pt]
 &U_t+(U\cdot\nabla)U-\eps\mu\Delta U+\nabla P=\eps Ne_d+G_u,  ~~\Div U=\sigma^\varepsilon, &&\text{in }Q^\varepsilon,\\[2pt]
 &\partial_\nu N=\partial_\nu \mathcal{C}=0,\quad U\cdot\nu=0,\quad\NN U=0, &&\text{on }\Sigma^\varepsilon,
\end{aligned}\right.
\end{equation}
where 
\begin{equation*}
 G_j(t,x)=\eps^2f_j(s_0+\eps t,x)\quad \text{with~~} j=n,c,u,\qquad
 \sigma^\varepsilon(t,x)=\eps\sigma(s_0+\eps t,x).
\end{equation*}
For the reference, we denote
\begin{equation*}
 \bar z^\eps(t,x)=\bar z(s_0+\eps t,x)
 \quad \text{with~~}  z=n,c,u,p,f_n,f_c,f_u.
\end{equation*}

Following Lemma 2, Lemma 3 and Appendix A in \cite{CMS2020}, we can construct the return flow as follows.
\begin{lemma}\label{lemma 3.1}
There are smooth exterior sources and a potential returning flow $U^0$, compactly supported in $0<t<1$, satisfying
\begin{equation}\label{3.2}
\left\{\begin{aligned}
 &U_t^0+(U^0\cdot\nabla)U^0+\nabla p^0=\xi^0, &&\text{in }Q_1=(0,1)\times \OO,\\[2pt]
 &\Div U^0=\sigma^0,\quad\curl U^0=0, &&\text{in }Q_1,\\[2pt]
 &U^0\cdot\nu=0, &&\text{on }\Sigma_1,\\[2pt]
 &U^0(0,\cdot)=U^0(1,\cdot)=0, &&\text{in }\OO.
\end{aligned}\right.
\end{equation}
The flow is defined by
\begin{equation}\label{3.3}
\left\{\begin{aligned}
 &\partial_t\Phi(t,x_0)=U^0(t,\Phi(t,x_0)), &&(t,x_0)\in(0,1)\times\overline\OO,\\[2pt]
 &\Phi(0,x_0)=x_0, &&x_0\in\overline\OO.
\end{aligned}\right.
\end{equation}
Every path of $\Phi$ visits $\overline\OO\setminus\overline\Om$. Given the fixed smooth initial velocity $u_{\rm in}$ and terminal values $\bar u^\eps(1)$, there are first order fields $z^\eps,p^{1,\eps},\xi^{1,\eps}$ on $[0,1]$ such that
\begin{equation}\label{3.4}
\left\{\begin{aligned}
 &z_t^\eps+(U^0\cdot\nabla)z^\eps+(z^\eps\cdot\nabla)U^0+\nabla p^{1,\eps}=\xi^{1,\eps}, &&\text{in }Q_1,\\[2pt]
 &\Div z^\eps=0, &&\text{in }Q_1,\\[2pt]
 &z^\eps\cdot\nu=0, &&\text{on }\Sigma_1,\\[2pt]
 &z^\eps(0,\cdot)=u_{\rm in},\quad z^\eps(1,\cdot)=\bar u^\eps(1,\cdot)=\ub(s_0+\varepsilon,\cdot), &&\text{in }\OO.
\end{aligned}\right.
\end{equation}
For every finite number of derivatives used below, their norms are bounded independently of small $\eps$.
\end{lemma}
\begin{proof}
As in \cite{CMS2020}, choose a potential flushing flow $W$ on $[0,1]$ and, for every smooth solenoidal tangent field $\lambda$, a solution $w_\lambda$ of
\[
 \partial_t w_\lambda+(W\cdot\nabla)w_\lambda+(w_\lambda\cdot\nabla)W+\nabla p_\lambda=\xi_\lambda,
 \qquad w_\lambda(0)=\lambda,\quad w_\lambda(1)=0,
\]
with $\Div w_\lambda=0$, $w_\lambda\cdot\nu=0$ and exterior-supported $\xi_\lambda$.
Compressing the construction into a strict subinterval, we may take $W=0$ near both endpoints and $w_\lambda=0$ near $t=1$. Set
\[
U^0(t)=
\begin{cases}
2W(2t),&0\le t<\frac12,\\
-2W(2-2t),&\frac12\le t\le1,
\end{cases}~~~~~~~~
z^\varepsilon(t)=
\begin{cases}
w_{u_{\rm in}}(2t),&0\le t<\frac12,\\
w_{\bar u^\varepsilon(1)}(2-2t),
&\frac12\le t\le1.
\end{cases}
\]
The first-order pressure and source are $2(p_{u_{\rm in}},\xi_{u_{\rm in}})(2t)$ on the first half and $-2(p_{\bar u^\varepsilon(1)},\xi_{\bar u^\varepsilon(1)})(2-2t)$ on the second half. Thus \eqref{3.4} holds, and the profiles match smoothly at $t=1/2$. The first half of the
construction already ensures the flushing property.
 The smooth construction in \cite{CMS2020}, applied at each finite Sobolev order, is linear in $\lambda$. Furthermore, the required uniform bounds follow from the bounded smooth family $\bar u(s_0+\varepsilon)$.
\end{proof}
 
Notice that we will first consider the special case
\begin{equation}\label{3.5}
\beta=2\mu M_w\qquad\text{on }\partial\mathcal O,
\end{equation}
where \(M_w\) is the shape operator defined in Section~\ref{subsection 2.2}. As in \cite{CMS2020}, the Navier boundary condition reduces to
\[
z\cdot\nu=0,\qquad
[(\nabla z)^T-\nabla z]\nu=0
\quad\text{on }\partial\mathcal O.
\]
The returning flow \(U^0\) satisfies \(U^0\cdot\nu=0\) and \(\operatorname{curl}U^0=0\), which imply \(\mathcal N_\beta U^0=0\). Therefore, no fluid boundary layer is needed in either dimension.

We introduce the following asymptotic expansion for the solution of \eqref{3.1}, for $0\le t <1$, 
\begin{equation*}
 V^\eps=U^0+\eps z^\eps,\quad
 P_V=p^0+\eps p^{1,\eps}+\eps\mu\sigma^0,\quad
 G_u=\xi^0+\eps\xi^{1,\eps},
\end{equation*}
while for $1\le t\le S_\eps$,  
\begin{equation*}
 V^\eps=\eps\bar u^\eps,\qquad P_V=\eps^2\bar p^\eps,\qquad
 G_u=\eps^2\bar f_u^\eps.
\end{equation*}
From \eqref{3.2} and \eqref{3.4}, the velocity profile $V^\varepsilon$ is continuous at $t=1$.
Let $\sigma^\varepsilon=\sigma^0$ on $(0,1)$ and $\sigma^\varepsilon=0$ on $(1,S_\varepsilon)$, which gives $\Div V^\varepsilon=\sigma^\varepsilon$ and
$
 -\eps\mu\Delta U^0+\nabla(\eps\mu\sigma^0)=0.
$
It follows from \eqref{3.2}-\eqref{3.4} that, on $(0,1)$,
\begin{equation}\label{3.6}
 V_t^\eps+(V^\eps\cdot\nabla)V^\eps-\eps\mu\Delta V^\eps+\nabla P_V-G_u
 =\eps^2\big((z^\eps\cdot\nabla)z^\eps-\mu\Delta z^\eps\big),
\end{equation}
and on $(1,S_\eps)$, 
\begin{equation}\label{3.7}
 V_t^\eps+(V^\eps\cdot\nabla)V^\eps-\eps\mu\Delta V^\eps+\nabla P_V-G_u
 =\eps^2\bar n^\eps e_d.
\end{equation}

\subsection{Transport toward a nonconstant reference}
For each $t$, the map $\Phi(t,\cdot)$ in \eqref{3.3} is a smooth diffeomorphism of $\overline\OO$. Let $X(x_0)=\Phi(1,x_0)$. We can choose
$\vartheta\in C^\infty([0,1]\times\overline\OO)$, nonnegative, compactly supported in time in $(0,1)$, with relative spatial support in $\overline\OO\setminus\overline\Omega$, such that
\begin{equation}\label{3.8}
 J(x_0)=\int_0^1\vartheta(t,\Phi(t,x_0))\dd t\ge \ell_{\min}>0,
 \qquad H(t,x_0)=J(x_0)^{-1}\int_0^t\vartheta(\tau,\Phi(\tau,x_0))\dd\tau.
\end{equation}
Using the flushing property and compactness, we can prove the uniform lower bound. For each $x_0\in\overline\OO$, the flushing property supplies $t_{x_0}\in(0,1)$ and a relative open neighborhood $B_{x_0}\subset\overline\OO\setminus\overline\Omega$ of $\Phi(t_{x_0},x_0)$. By continuity, there exist a neighborhood \(U_{x_0}\) of \(x_0\), a closed interval \(I_{x_0}\subset\subset(0,1)\) of positive length, and a relative neighborhood \(B'_{x_0}\subset\subset B_{x_0}\) such that
\[
\Phi(t,y_0)\in B'_{x_0},
\qquad \text{for all }(t,y_0)\in I_{x_0}\times U_{x_0}.
\]
We choose smooth nonnegative cutoff functions \(\theta_{x_0}\) and \(\psi_{x_0}\) such that \(\theta_{x_0}=1\) on \(I_{x_0}\) and \(\psi_{x_0}=1\) on \(B'_{x_0}\), with \(\operatorname{supp}\theta_{x_0}\subset\subset(0,1)\) and \(\operatorname{supp}\psi_{x_0}\subset\subset B_{x_0}\) in the relative topology. By compactness of \(\overline{\mathcal O}\), we can select points \(x_1,\ldots,x_{N_{\rm cov}}\) whose neighborhoods \(U_{x_j}\) cover \(\overline{\mathcal O}\). Then, let us define
\[
\vartheta(t,x)=\sum_{j=1}^{N_{\rm cov}}\theta_{x_j}(t)\psi_{x_j}(x),
\qquad
\ell_{\min}=\min_{1\le j\le N_{\rm cov}}|I_{x_j}|>0.
\]
For any \(y_0\in U_{x_j}\), we have 
$
J(y_0)\geq |I_{x_j}|\geq \ell_{\min},
$
which yields \eqref{3.8}. Since trajectories starting on \(\partial\mathcal O\) remain on the boundary, the global exterior controls are allowed to intersect \(\partial\mathcal O\).
Thus,
\begin{equation*}
 0\le H\le1,\quad H_t=\vartheta(t,\Phi(t,x_0))/J(x_0),\quad
 H=0\text{ near }0,\quad H=1\text{ near }1.
\end{equation*}
With the superscript $\mathrm{tr}$ denoting transport along the return flow, define
\begin{align}
 n^{\rm tr}(t,\Phi(t,x_0))&=(1-H(t,x_0))n_{\rm in}(x_0)
                    +H(t,x_0)\bar n^\eps(t,X(x_0)),\label{3.9}\\
 c^{\rm tr}(t,\Phi(t,x_0))&=(1-H(t,x_0))c_{\rm in}(x_0)
                    +H(t,x_0)\bar c^\eps(t,X(x_0)).\label{3.10}
\end{align}
Here, we use \(x_0\) to denote the particle label and evaluate the target at its terminal position \(X(x_0)\). With this choice, the transport identities below involve only the time derivatives of the reference trajectory with order \(O(\varepsilon)\).

Accordingly, we take
\begin{align*}
 g_n^{\rm tr}(t,\Phi(t,x_0))&=H_t(t,x_0)
     \big(\bar n^\eps(t,X(x_0))-n_{\rm in}(x_0)\big),\\
 g_c^{\rm tr}(t,\Phi(t,x_0))&=H_t(t,x_0)
     \big(\bar c^\eps(t,X(x_0))-c_{\rm in}(x_0)\big).
\end{align*}
The sources are supported in the exterior. Thus, differentiating \eqref{3.9} along the flow, one has
\begin{equation*}
 \frac{\dd}{\dd t}n^{\rm tr}(t,\Phi(t,x_0))
 =n^{\rm tr}_t(t,\Phi(t,x_0))+\nabla n^{\rm tr}(t,\Phi(t,x_0))\cdot U^0(t,\Phi(t,x_0)).
\end{equation*}
Differentiating the right hand side of \eqref{3.9}, we get
$
 H_t\big(\bar n^\eps(t,X(x_0))-n_{\rm in}(x_0)\big)
 +\eps H(t,x_0)(\nb_s)(s_0+\eps t,X(x_0)).
$
Therefore
\begin{equation*}
\left\{\begin{aligned}
 &\begin{aligned}[b] &n^{\rm tr}_t+U^0\cdot\nabla n^{\rm tr}=g_n^{\rm tr}+\eps H(t,\Phi^{-1}(t,x)) \cdot\nb_s(s_0+\eps t,X(\Phi^{-1}(t,x))),\end{aligned} &&\text{in }Q_1,\\[2pt]
 &\begin{aligned}[b] &c^{\rm tr}_t+U^0\cdot\nabla c^{\rm tr}=g_c^{\rm tr}+\eps H(t,\Phi^{-1}(t,x)) \cdot\cb_s(s_0+\eps t,X(\Phi^{-1}(t,x))),\end{aligned} &&\text{in }Q_1.
\end{aligned}\right.
\end{equation*}
Since \(H=1\) and \(\Phi(t,x_0)=X(x_0)\) for \(t\) sufficiently close to \(1\), we have \(n^{\rm tr}=\bar n^\varepsilon\) and \(c^{\rm tr}=\bar c^\varepsilon\). We therefore extend the profiles to \(t\geq1\) by setting \(n^{\rm tr}=\bar n^\varepsilon\) and \(c^{\rm tr}=\bar c^\varepsilon\).

On the other hand, by \eqref{3.9}-\eqref{3.10}, for
$
 m_*:=\min\{\min_{\overline\OO}n_{\rm in},\min_{\overline\OO}c_{\rm in},m_{s_*}\}>0,
$
we have
\begin{equation*}
 n^{\rm tr}\ge m_*,\qquad c^{\rm tr}\ge m_*.
\end{equation*}
Moreover, the smoothness of the flow map, its inverse, and the reference trajectory implies that, for every integer \(m\geq0\),
\[
\sup_{0<\varepsilon\leq\varepsilon_0}
\left(
\|n^{\rm tr}\|_{C^m(\overline{Q_1})}
+\|c^{\rm tr}\|_{C^m(\overline{Q_1})}
\right)\leq C_m,
\]
where \(C_m\) is independent of \(\varepsilon\).

\subsection{Scalar Neumann corrections and defect analysis}
We next construct scalar boundary correctors to enforce the homogeneous Neumann conditions, which are not necessarily satisfied by the transported profiles.  
Let $\widetilde d(x)=\operatorname{dist}(x,\partial\OO)$
in a neighborhood of $\partial\OO$, and $\tilde{\pi}(x)$ the nearest boundary point. Choosing $\psi\in C_c^\infty([0,\infty))$ with $\psi'(0)=1$, for $z\in\{n,c\}$ we define the Neumann corrector for $z^{\rm tr}$ and the corrected profiles by 
\begin{equation*}
 \ell_z^\eps=\sqrt\eps\,\partial_\nu z^{\rm tr}(t,\tilde{\pi}(x))\psi(\tilde{d}(x)/\sqrt\eps),
 \qquad n^{\rm app}=n^{\rm tr}+\ell_n^\eps,\qquad c^{\rm app}=c^{\rm tr}+\ell_c^\eps,
\end{equation*}
and the correction is extended by zero outside the collar neighborhood. Since $\partial_\nu \tilde{d}=-1$,
\begin{equation*}
 \partial_\nu\ell_z^\eps|_{\tilde{d}=0}=-\partial_\nu z^{\rm tr}\,\psi'(0)=-\partial_\nu z^{\rm tr}.
\end{equation*}
Thus, one can see
\begin{equation*}
\left\{\begin{aligned}
 &\partial_\nu n^{\rm app}=\partial_\nu c^{\rm app}=0, &&\text{on }\Sigma^\varepsilon,\\[2pt]
 &(n^{\rm app},c^{\rm app})|_{t=0}=(n_{\rm in},c_{\rm in}), &&\text{in }\OO,\\[2pt]
 &(n^{\rm app},c^{\rm app})=(\bar n^\eps,\bar c^\eps), &&\text{in }(1,S_\eps)\times\OO.
\end{aligned}\right.
\end{equation*}
Differentiating and integrating directly in the normal variable, we have
\begin{align}
 \norm[2]{\ell_z^\eps}+\norm[2]{\partial_t\ell_z^\eps}\le C\eps^{3/4}, \quad
 \norm[2]{\nabla\ell_z^\eps}\le C\eps^{1/4},  \quad \norm[2]{\Delta\ell_z^\eps}\le C\eps^{-1/4},  \quad\norm[\infty]{\ell_z^\eps}\le C\sqrt\eps,  \quad
 \norm[\infty]{\nabla\ell_z^\eps}\le C,\label{3.11}
\end{align}
and
\begin{equation}\label{3.12}
 \norm[W^{1,\infty}]{n^{\rm app}}+\norm[W^{1,\infty}]{c^{\rm app}}\le C,\qquad
 n^{\rm app},c^{\rm app}\ge m_*/2.
\end{equation}

Define
\begin{equation}\label{3.13}
 \alpha_n^\eps=\frac{g_n^{\rm tr}+\eps\bar f_n^\eps}{n^{\rm app}},\qquad
 \alpha_c^\eps=\frac{g_c^{\rm tr}+\eps\bar f_c^\eps}{c^{\rm app}},\qquad
 G_n=\alpha_n^\eps N,\quad G_c=\alpha_c^\eps \mathcal{C}.
\end{equation}
Since \(g_n^{\rm tr}=g_c^{\rm tr}=0\) for \(t\geq1\) and \(g_n^{\rm tr},g_c^{\rm tr},\bar f_n^\varepsilon,\bar f_c^\varepsilon\) are supported in the exterior region, we obtain, for \(1\leq t\leq T_{\rm app}/\varepsilon\),
\begin{equation*}
 \alpha_n^\eps=\eps\an(s_0+\eps t),\qquad
 \alpha_c^\eps=\eps\ac(s_0+\eps t).
\end{equation*}
Thus, using \eqref{3.12} and \eqref{3.13}, we obtain
\begin{equation}\label{3.14}
 \norm[W^{1,\infty}]{\alpha_n^\eps(t)}+\norm[W^{1,\infty}]{\alpha_c^\eps(t)}
 \le C\big(\mathbf1_{[0,1]}(t)+\eps\big)
\end{equation}
and
\[
\varepsilon\int_0^{S_\varepsilon}
\left(
\|\nabla n^{\rm app}(t)\|_{L^\infty(\mathcal O)}
+\|\nabla c^{\rm app}(t)\|_{L^\infty(\mathcal O)}
\right)\,dt
\leq C(1+T_{\rm app}),
\]
where \(C\) is independent of \(\varepsilon\).

The corresponding scalar sources are
\begin{equation}\label{3.15}
 f_n(s,x)=\eps^{-1}\alpha_n^\eps((s-s_0)/\eps,x)n(s,x),\quad
 f_c(s,x)=\eps^{-1}\alpha_c^\eps((s-s_0)/\eps,x)c(s,x),
\end{equation}
which vanish in $\Omega$. For fixed $\varepsilon$, Lemma~\ref{lemma 2.2} preserves the signs of the actual scalar states.
 
To estimate the scalar residuals arising from the fluid transport, we introduce the following quantities for the fluid profile $V^\eps$  
\begin{align}
 R_n&=g_n^{\rm tr}+\eps\bar f_n^\eps-
       \big(n^{\rm app}_t+V^\eps\cdot\nabla n^{\rm app}-\eps D_n\Delta n^{\rm app}+\eps\chi\Div(n^{\rm app}\nabla c^{\rm app})\big),\label{3.16}\\
 R_c&=g_c^{\rm tr}+\eps\bar f_c^\eps-
       \big(c^{\rm app}_t+V^\eps\cdot\nabla c^{\rm app}-\eps D_c\Delta c^{\rm app}+\eps n^{\rm app}c^{\rm app}\big).\label{3.17}
\end{align}
Since $V^\eps$ is tangent and uniformly Lipschitz on $[0,1]$,
\begin{equation*}
 |V^\eps(t,x)\cdot\nabla \tilde{d}(x)|\le C \tilde{d}(x),\qquad
 \norm[2]{V^\eps\cdot\nabla\ell_z^\eps}\le C\eps^{3/4}.
\end{equation*}
Using the fact \(V^\varepsilon\cdot\nu=0\) on \(\partial\mathcal O\), we also have \(|V^\varepsilon\cdot\nu|\leq C\widetilde d\leq C\sqrt{\varepsilon}\) in the boundary neighborhood, which, together with the tangential derivative estimates, yields
\[
\|V^\varepsilon\cdot\nabla\ell_z^\varepsilon\|_{L^\infty}
\leq C\sqrt{\varepsilon}.
\]

On $[0,1]$, for $R_n$, one has
\begin{align*}
&\left\|
\varepsilon H(t)\,
\partial_s\bar n
\left(
s_0+\varepsilon t,\,
X\bigl(\Phi^{-1}(t,x)\bigr)
\right)
\right\|_{L^2L^2}+ \norm[L^2L^2]{(V^\eps-U^0)\cdot\nabla n^{\rm tr}}
\le C\varepsilon,\\
& \norm[L^2L^2]{(\partial_t+V^\eps\cdot\nabla)\ell_n^\eps}  + \eps D_n\norm[L^2L^2]{\Delta n^{\rm app}} \le C\eps^{3/4},\\
& \eps\chi\norm[L^2L^2]{\Div(n^{\rm app}\nabla c^{\rm app})}
 \le C\eps\big(\norm[\infty]{\nabla n^{\rm app}}\norm[L^2L^2]{\nabla c^{\rm app}}
             +\norm[\infty]{n^{\rm app}}\norm[L^2L^2]{\Delta c^{\rm app}}\big) \le C\eps^{3/4},\\ &\eps\norm[L^2L^2]{\bar f_n^\eps}\le C\eps,
\end{align*}
which implies 
$$
\|R_n\|_{L^2(0,1;L^2)}
\le C\varepsilon^{3/4}.
$$ 
Similarly, using $\varepsilon\|n^{\rm app}c^{\rm app}\|_{L^2(0,1;L^2)}
\le C\varepsilon,$  we have
\[
\|R_c\|_{L^2(0,1;L^2)}
\le C\varepsilon^{3/4}.
\]
While on $t\ge1$, we get
\begin{equation*}
\left\{\begin{aligned}
 &\begin{aligned}[b] &n^{\rm app}_t+\eps\bar u^\eps\cdot\nabla n^{\rm app}-\eps D_n\Delta n^{\rm app} +\eps\chi\Div(n^{\rm app}\nabla c^{\rm app})=\eps\bar f_n^\eps,\end{aligned} &&\text{in }(1,S_\eps)\times\OO,\\[2pt]
 &c^{\rm app}_t+\eps\bar u^\eps\cdot\nabla c^{\rm app}-\eps D_c\Delta c^{\rm app}+\eps n^{\rm app}c^{\rm app}=\eps\bar f_c^\eps, &&\text{in }(1,S_\eps)\times\OO.
\end{aligned}\right.
\end{equation*}
Since $R_n=R_c=0$ for $t\ge1$, we obtain
\begin{equation*}
 \norm[L^2(0,S_\eps;L^2)]{R_n}+\norm[L^2(0,S_\eps;L^2)]{R_c}\le C\eps^{3/4}.
\end{equation*}

For the fluid defect, one has 
\begin{equation}\label{3.18}
 F_u=\eps^{-1}\big[G_u+\eps^2n^{\rm app} e_d-
 (V_t^\eps+(V^\eps\cdot\nabla)V^\eps-\eps\mu\Delta V^\eps+\nabla P_V)\big].
\end{equation}
From \eqref{3.6}-\eqref{3.7}, we get
\begin{equation}\label{3.19}
 F_u=\eps\big(n^{\rm app}e_d+\mu\Delta z^\eps-(z^\eps\cdot\nabla)z^\eps\big)
 \quad(t<1),\qquad F_u=0\quad(t>1),
\end{equation}
and since $\NN U^0=\NN\bar u^\eps=0$,
\begin{equation}\label{3.20}
 h^\eps=-\eps^{-1}\NN V^\eps
 =-\NN z^\eps\quad(t<1),\qquad h^\eps=0\quad(t>1).
\end{equation}
Thus, it holds that
\begin{equation*}
 \norm[L^1L^2]{F_u}\le C\eps,\qquad
 \eps\norm[L^2L^2(\partial\OO)]{h^\eps}^2\le C\eps
\end{equation*}
and
\begin{equation*}
 \int_0^{S_\eps}\big(\norm[\infty]{\nabla V^\eps}+\norm[\infty]{\sigma^\varepsilon}\big)\dd t
 \le C+C\eps(S_\eps-1)\le C(1+T_{\rm app}).
\end{equation*}

\subsection{Estimates for the remainder}
\subsubsection{Preliminaries}
Let
\begin{equation*}
 N=\eps(n^{\rm app}+r),\qquad \mathcal{C}=\eps(c^{\rm app}+q),\qquad
 U=V^\eps+\eps y,\qquad P=P_V+\eps\pi.
\end{equation*}
The remainders $r,q,y$ satisfy the following system
\begin{equation}\label{3.21}
\left\{\begin{aligned}
 &\begin{aligned}[b] &r_t+U\cdot\nabla r-\eps D_n\Delta r =R_n+\alpha_n^\eps r-\eps y\cdot\nabla n^{\rm app} -\eps\chi\Div(n^{\rm app}\nabla q+r\nabla c^{\rm app}+r\nabla q),\end{aligned} &&\text{in }Q^\varepsilon,\\[2pt]
 &\begin{aligned}[b] &q_t+U\cdot\nabla q-\eps D_c\Delta q =R_c+\alpha_c^\eps q-\eps y\cdot\nabla c^{\rm app}-\eps(n^{\rm app}q+c^{\rm app}r+rq),\end{aligned} &&\text{in }Q^\varepsilon,\\[2pt]
 &\begin{aligned}[b] &y_t+(U\cdot\nabla)y+(y\cdot\nabla)V^\eps -\eps\mu\Delta y+\nabla\pi=F_u+\eps r e_d,\end{aligned}  ~~~~\Div y=0, &&\text{in }Q^\varepsilon,\\[2pt]
 &\partial_\nu r=\partial_\nu q=0,\quad y\cdot\nu=0,\quad\NN y=h^\eps, &&\text{on }\Sigma^\varepsilon,\\[2pt]
 &(r,q,y)|_{t=0}=0, &&\text{in }\OO.
\end{aligned}\right.
\end{equation}
 
We now derive uniform \(L^4\) estimates for the remainder system. In order to control the nonlinear coupling between the scalar and fluid equations, we use the equivalent energy quantities 
\begin{equation}\label{3.22}
 \begin{aligned}
 (\|r\|_4^4+\|q\|_{W^{1,4}}^4)^{1/2}+\|y\|_4^2
 \le \|r\|_4^2+\|q\|_{W^{1,4}}^2+\|y\|_4^2 \le\sqrt2\big[(\|r\|_4^4+\|q\|_{W^{1,4}}^4)^{1/2}+\|y\|_4^2\big].
 \end{aligned}
\end{equation}
In what follows, the estimates require the following uniform bounds on the prescribed profiles and control coefficients
\begin{equation}\label{3.23}
 \begin{aligned}
 \|n^{\rm app}\|_{W^{1,\infty}}+\|c^{\rm app}\|_{W^{1,\infty}}\le C_{\rm prof}, 
 ~~~~\int_0^{S_\varepsilon}\big(\|\nabla V^\varepsilon\|_\infty
   +\|\Div V^\varepsilon\|_\infty+\|\alpha_n^\varepsilon\|_\infty
   +\|\alpha_c^\varepsilon\|_{W^{1,\infty}}\big)\dd t\le C_{\rm prof}.
 \end{aligned}
\end{equation}
The constant $C_{\rm prof}$ is independent of $\varepsilon$ and
$S_\varepsilon=T_{\rm app}/\varepsilon$.

Thus, by the Gagliardo-Nirenberg and trace inequalities,
for $d=2,3$, we have
\begin{align}
&\|f\|_4^2
\le C\|f\|_2^{1/2}\|\nabla f\|_2^{3/2}
   +C\|f\|_2^2, \label{3.24}\\
&\|f\|_{L^2(\partial\mathcal O)}^2
\le \eta\|\nabla f\|_2^2
   +C_\eta\|f\|_2^2. \label{3.25}
\end{align}
In particular, applying \eqref{3.24} to $f=r^2$
and using Young's inequality, we obtain
\begin{align*}\begin{aligned}
\int_{\mathcal O}r^4|\nabla q|^2\,dx
&\le \|\nabla q\|_4^2\|r^2\|_4^2 \le C\|\nabla q\|_4^2\|r\|_4
   \left(\int_{\mathcal O}
       r^2|\nabla r|^2\,dx\right)^{3/4}
   +C\|\nabla q\|_4^2\|r\|_4^4\\
&\le \eta\int_{\mathcal O}
       r^2|\nabla r|^2\,dx
   +C_\eta(1+\|\nabla q\|_4^8)\|r\|_4^4.
\end{aligned}\end{align*}

\subsubsection{Estimate for the scalar remainders}
Multiplying \eqrow{3.21}{1} by $r^3$, one has
\begin{align}
 \frac14\frac{\dd}{\dd t}\|r\|_4^4+3\varepsilon D_n\left(\int r^2|\nabla r|^2\right)
 ={}&\frac14\int\sigma^\varepsilon r^4+\int\alpha_n^\varepsilon r^4+\int R_nr^3
       -\varepsilon\int(y\cdot\nabla n^{\rm app})r^3\nonumber\\
 &+3\varepsilon\chi\int n^{\rm app} r^2\nabla q\cdot\nabla r
  +3\varepsilon\chi\int r^3\nabla c^{\rm app}\cdot\nabla r
  +3\varepsilon\chi\int r^3\nabla q\cdot\nabla r.\label{3.26}
\end{align}
From
$\int U\cdot\nabla(r^4)=-\int\sigma^\varepsilon r^4$, the first four terms satisfy
\begin{align*}
 \frac14\left|\int\sigma^\varepsilon r^4\right|+\left|\int\alpha_n^\varepsilon r^4\right|
 &\le C(\|\Div V^\varepsilon\|_\infty+\|\alpha_n^\varepsilon\|_\infty)\|r\|_4^4,~~~~
 \left|\int R_nr^3\right|\le\|R_n\|_4\|r\|_4^{3},\\
 \varepsilon\left|\int(y\cdot\nabla n^{\rm app})r^3\right|&\le C\varepsilon \|y\|_4 \|r\|_4^{3}.
\end{align*}
By Young's inequality, we also have
\begin{align*}
 3\varepsilon\chi\left|\int n^{\rm app} r^2\nabla q\cdot\nabla r\right|
 &\le\eta\varepsilon \left(\int r^2|\nabla r|^2\right)+C_\eta\varepsilon\int r^2|\nabla q|^2
\le\eta\varepsilon \int r^2|\nabla r|^2+C_\eta\varepsilon(\|r\|_4^4+\|\nabla q\|_4^4),\\
 3\varepsilon\chi\left|\int r^3\nabla c^{\rm app}\cdot\nabla r\right|
 &\le\eta\varepsilon \left(\int r^2|\nabla r|^2\right)+C_\eta\varepsilon \|r\|_4^{4},\\
 3\varepsilon\chi\left|\int r^3\nabla q\cdot\nabla r\right|
 &\le\eta\varepsilon \left(\int r^2|\nabla r|^2\right)+C_\eta\varepsilon\int r^4|\nabla q|^2 
 \le2\eta\varepsilon \left(\int r^2|\nabla r|^2\right)+C_\eta\varepsilon(1+\|\nabla q\|_4^{8})\|r\|_4^{4}.
\end{align*}
 
Similarly, testing \eqref{3.21}$_2$ by $q^3$,
we obtain
\begin{equation}\label{3.27}
 \frac14\frac{\dd}{\dd t}\|q\|_4^{4}+3\varepsilon D_c\left(\int q^2|\nabla q|^2\right)
 =\frac14\int\sigma^\varepsilon q^4+\int\alpha_c^\varepsilon q^4+\int b_cq^3,
\end{equation}
where 
$ b_c=R_c-\varepsilon y\cdot\nabla c^{\rm app}-\varepsilon n^{\rm app}q-\varepsilon c^{\rm app}r-\varepsilon rq.$
All the terms in $b_c$  satisfy
\begin{align*}
 &\|\varepsilon y\cdot\nabla c^{\rm app}\|_4\le C\varepsilon \|y\|_4,~~~~
 \|\varepsilon n^{\rm app}q\|_4\le C\varepsilon \|q\|_4,\\
 &\|\varepsilon c^{\rm app}r\|_4\le C\varepsilon \|r\|_4,~~~~
 \|\varepsilon rq\|_4\le C\varepsilon \|r\|_4\|q\|_{W^{1,4}}.
\end{align*}

For the gradient, let $g=\nabla q$ and then multiply the same equation by
$-\operatorname{div}(|g|^2g)$. It holds that 
\begin{equation}\label{3.28}
\begin{aligned}
 &\frac14\frac{\dd}{\dd t}\|\nabla q\|_4^4
 +\varepsilon D_c\int|\nabla q|^2|D^2q|^2
 +\frac{\varepsilon D_c}{2}\int|\nabla|\nabla q|^2|^2\\
 &\quad=\frac14\int\sigma^\varepsilon|g|^4
 -\sum_{i,j=1}^d\int\partial_iU_j\,g_jg_i|g|^2
 +\int\alpha_c^\varepsilon|g|^4 +\int q\nabla\alpha_c^\varepsilon\cdot g|g|^2
 -\int b_c\operatorname{div}(|g|^2g)
 +\frac{\varepsilon D_c}{2}\int_{\partial\OO}|g|^2\partial_\nu|g|^2.
\end{aligned}
\end{equation}
The diffusion term is given by
\[
 \int\nabla\Delta q\cdot g|g|^2
 =-\int|g|^2|D^2q|^2-\frac12\int|\nabla|g|^2|^2
       +\frac12\int_{\partial\mathcal O}|g|^2\partial_\nu|g|^2.
\]
Since
\[
 \int |q||g|^3\le \|q\|_4\|\nabla q\|_4^{3}\le\tfrac14\|q\|_4^{4}+\tfrac34\|\nabla q\|_4^{4},
\]
the terms involving $\sigma^\varepsilon$, $\nabla V^\varepsilon$ and $\alpha_c^\varepsilon$ can be bounded by $C\big(\|\nabla V^\varepsilon\|_\infty+\|\Div V^\varepsilon\|_\infty+\|\alpha_n^\varepsilon\|_\infty+\|\alpha_c^\varepsilon\|_{W^{1,\infty}}\big)(\|q\|_4^{4}+\|\nabla q\|_4^{4})$.
By integration by parts, one has
\begin{align*}
  \left|\sum_{i,j}\int y_j\partial_i(g_jg_i|g|^2)\right|
 &\le C\int|y||g|^3|D^2q| \le\eta \left(\int |\nabla q|^2|D^2q|^2\right)+C_\eta\int|y|^2|g|^4\\
 &\le2\eta \left(\int |\nabla q|^2|D^2q|^2\right)+C_\eta(1+\|y\|_4^{8})\|\nabla q\|_4^{4}.
\end{align*}
 
From $\partial_\nu q=0$, we have
\[
 \partial_\nu|g|^2=-2(M_wg)\cdot g\quad\hbox{on }\partial\mathcal O.
\]
Applying \eqref{3.25} with $|g|^2$, we have
\begin{equation*}
 \left|\frac{\varepsilon D_c}{2}\int_{\partial\mathcal O}|g|^2\partial_\nu|g|^2\right|
 \le\eta\varepsilon \left(\int |\nabla q|^2|D^2q|^2\right)+C_\eta\varepsilon \|\nabla q\|_4^{4}.
\end{equation*}
Finally, $|\operatorname{div}(|g|^2g)|\le C|g|^2|D^2q|$. For any one of the five sources $f$ in $b_c$, we can get
\begin{equation*}
\begin{aligned}
 \left|\int f\operatorname{div}(|g|^2g)\right|
 &\le C\|f\|_4\|\nabla q\|_4
           \left(\int|\nabla q|^2|D^2q|^2\right)^{1/2} \le\eta\varepsilon\int|\nabla q|^2|D^2q|^2
       +C_\eta\varepsilon^{-1}\|f\|_4^2\|\nabla q\|_4^2.
\end{aligned}
\end{equation*}
The five resulting terms are bounded by, respectively,
\begin{equation*}
\begin{gathered}
 C_\eta\varepsilon^{-1}\|R_c\|_4^2\|\nabla q\|_4^2,\qquad
 C_\eta\varepsilon\|y\|_4^2\|\nabla q\|_4^2,\qquad
 C_\eta\varepsilon\|q\|_4^2\|\nabla q\|_4^2,\\
 C_\eta\varepsilon\|r\|_4^2\|\nabla q\|_4^2,\qquad
 C_\eta\varepsilon\|r\|_4^2\|q\|_{W^{1,4}}^2\|\nabla q\|_4^2,
\end{gathered}
\end{equation*}
where the bounds depend on $\|R_c\|_4$ and $\|c^{\rm app}\|_{W^{1,\infty}}$.

Thus, adding \eqref{3.26}, \eqref{3.27}
and \eqref{3.28} with their bounds, and choosing $\eta$ sufficiently small,
we obtain
\begin{align*}
\frac14\frac{d}{dt}
\left(\|r\|_4^4+\|q\|_{W^{1,4}}^4\right)
&+\gamma_0\varepsilon\int_{\mathcal O}
      ( r^2|\nabla r|^2\, +
       q^2|\nabla q|^2\, +
       |\nabla q|^2|D^2q|^2\,)dx\\
&\le C\Big(
       \|\nabla V^\varepsilon\|_\infty
      +\|\operatorname{div}V^\varepsilon\|_\infty
      +\|\alpha_n^\varepsilon\|_\infty
      +\|\alpha_c^\varepsilon\|_{W^{1,\infty}}
      \Big)
      \left(\|r\|_4^4+\|q\|_{W^{1,4}}^4\right)\\
&\quad+C\varepsilon(1+\|y\|_4^8+\|\nabla q\|_4^8)
      \left(\|r\|_4^4+\|q\|_{W^{1,4}}^4\right)+C\varepsilon\|y\|_4
      \left(\|r\|_4^4+\|q\|_{W^{1,4}}^4\right)^{3/4}\\
&\quad+C\varepsilon\|r\|_4^2
      \|q\|_{W^{1,4}}^2
      \|\nabla q\|_4^2 +C\|R_n\|_4\|r\|_4^3 + C\varepsilon^{-1}\|R_c\|_4^2
      \left(\|r\|_4^4+\|q\|_{W^{1,4}}^4\right)^{1/2}.
\end{align*}

Define
\[
E(t)=\|r(t)\|_4^4+\|q(t)\|_{W^{1,4}}^4.
\]
To justify the differentiation of the square-root
energies when they may vanish, we first replace
$E^{1/2}$ by $(E+\kappa)^{1/2}$, with $\kappa>0$.
We derive the corresponding integrated inequality
and then pass to the limit as $\kappa\to0^+$.
The same regularization applies to
$\|y\|_4^2$ in the velocity estimate.
 If $0<E^{1/2}\le1$, it holds that
\[
\varepsilon
\frac{\|r\|_4^2\|q\|_{W^{1,4}}^2
      \|\nabla q\|_4^2}{E^{1/2}}
\le C\varepsilon E
\le C\varepsilon E^{1/2}.
\]
Thus, we obtain
\begin{equation}\label{3.29}
\begin{aligned}
\frac{d}{dt}E^{1/2}
&\le C\Big(
       \|\nabla V^\varepsilon\|_\infty
      +\|\operatorname{div}V^\varepsilon\|_\infty
      +\|\alpha_n^\varepsilon\|_\infty
      +\|\alpha_c^\varepsilon\|_{W^{1,\infty}}
      \Big)E^{1/2}\\
&\quad+C\varepsilon(1+\|y\|_4^8)E^{1/2}
      +C\varepsilon\|y\|_4^2 +C\|R_n\|_4E^{1/4}
      +C\varepsilon^{-1}\|R_c\|_4^2.
\end{aligned}
\end{equation}

\subsubsection{Estimate for the velocity remainders}

In order to establish the \(L^4\) estimate of the velocity remainder, we first derive bounds for the convective and viscous pressure contributions under the Navier boundary condition.

For a given vector field $f$, let us consider zero-mean Neumann pressures.
 \begin{equation*}
\left\{\begin{aligned}
 &\Delta p=\Div f, &&\text{in }\OO,\\[2pt]
 &\partial_\nu p=f\cdot\nu, &&\text{on }\partial\OO,\\[2pt]
 &\displaystyle\int_\OO p\,\dd x=0. &&
\end{aligned}\right.
\end{equation*}
From standard elliptic regularity estimates, one has $\|\nabla p\|_4\le C\|f\|_4$. The vector $f-\nabla p$ is divergence free and tangent.
Let 
\[
 b_0=(V^\varepsilon\cdot\nabla)y-(y\cdot\nabla)V^\varepsilon+\operatorname{div}V^\varepsilon y.
\]
It holds that
\begin{align*}
 \operatorname{div}\big((V^\varepsilon\cdot\nabla)y-(y\cdot\nabla)V^\varepsilon\big)
 =-y\cdot\nabla\operatorname{div}V^\varepsilon,~~~~
 \operatorname{div}(\operatorname{div}V^\varepsilon y)=y\cdot\nabla\operatorname{div}V^\varepsilon.
\end{align*}
On the boundary, differentiating tangency along a tangent vector gives
$(V^\varepsilon\cdot\nabla)y\cdot\nu=-y\cdot M_w V^\varepsilon$ and
$(y\cdot\nabla)V^\varepsilon\cdot\nu=-V^\varepsilon\cdot M_w y$, 
which are equal because $M_w$ is symmetric. Hence $b_0$ is divergence free and tangent. Define the zero-mean pressure $\pi_L$ by
\[
 \nabla\pi_L=(I-\mathbb P_4)\big((V^\varepsilon\cdot\nabla)y+(y\cdot\nabla)V^\varepsilon\big)
 =(I-\mathbb P_4)\big(2(y\cdot\nabla)V^\varepsilon-(\Div V^\varepsilon)y\big).
\]
It follows that
\begin{equation}\label{3.30}
\begin{aligned}
 & \|\nabla\pi_L\|_4
\le C(\|\nabla V^\varepsilon\|_\infty+\|\Div V^\varepsilon\|_\infty)\|y\|_4,\\
 &\left|\int\nabla\pi_L\cdot|y|^2y\right|
\le C(\|\nabla V^\varepsilon\|_\infty+\|\Div V^\varepsilon\|_\infty)\|y\|_4^4.
\end{aligned}
\end{equation}

For the viscous pressure $\pi_S$, the equations are
\begin{equation}\label{3.31}
\left\{\begin{aligned}
 &\Delta\pi_S=0, &&\text{in }\OO,\\[2pt]
 &\partial_\nu\pi_S=-\Div_\Gamma((2\mu M_w-\beta_\tau)y+h^\varepsilon), &&\text{on }\partial\OO,\\[2pt]
 &\displaystyle\int_\OO\pi_S\,\dd x=0. &&
\end{aligned}\right.
\end{equation}
In 3D the Neumann boundary condition follows from $\Delta y=-\operatorname{curl}\operatorname{curl}y$ and
\[
 (\operatorname{curl}y)\times\nu=(2M_w-\beta_\tau/\mu)y+h^\varepsilon/\mu.
\]
In 2D, it follows from the corresponding scalar vorticity formula. 
To prove its required bound, take a mean-zero $f\in L^{4/3}$ and solve
\begin{equation*}
\left\{\begin{aligned}
 &-\Delta\phi=f, &&\text{in }\OO,\\[2pt]
 &\partial_\nu\phi=0, &&\text{on }\partial\OO,\\[2pt]
 &\displaystyle\int_\OO\phi\,\dd x=0. &&
\end{aligned}\right.
\end{equation*}
Then $\|\phi\|_{W^{2,4/3}}\le C\|f\|_{4/3}$.
By Green's formula and integration by parts, one gets
\[
 \left|\int\pi_S f\right|
 =\left|\int_{\partial\mathcal O}\big((2\mu M_w-\beta_\tau)y+h^\varepsilon\big)\cdot\nabla_\Gamma\phi\right|
 \le C(\|y\|_{4,\partial\mathcal O}+\|h^\varepsilon\|_{4,\partial\mathcal O})\|f\|_{4/3}.
\]
Here the trace of $\nabla\phi$ is bounded in $L^{4/3}(\partial\mathcal O)$ by its $W^{1,4/3}$ norm. Using the duality between $L^4$ and $L^{4/3}$, we have 
\begin{equation*}
 \|\pi_S\|_4\le C(\|y\|_{4,\partial\mathcal O}+\|h^\varepsilon\|_{4,\partial\mathcal O}).
\end{equation*}
Applying the multiplicative trace inequality used in \eqref{3.25}, for $|y|^2$, we have 
\begin{equation*}
 \|y\|_{4,\partial\mathcal O}\le C \|y\|_4^{1/2}\left(\int |y|^2|\nabla y|^2\right)^{1/8}+C\|y\|_4.
\end{equation*}
It follows from $\operatorname{div}y=0$ that
$|\operatorname{div}(|y|^2y)|\le2|y|^2|\nabla y|$. Hence, one has
\begin{align}
 &\left|\int\pi_S\operatorname{div}(|y|^2y)\right|\\
\le& C\|\pi_S\|_4\|y\|_4 \left(\int |y|^2|\nabla y|^2\right)^{1/2}\nonumber\\
\le& C \|y\|_4^{3/2}\left(\int |y|^2|\nabla y|^2\right)^{5/8}+C\|y\|_4^{2}\left(\int |y|^2|\nabla y|^2\right)^{1/2} +C\|h^\varepsilon\|_{4,\partial\mathcal O}\|y\|_4 \left(\int |y|^2|\nabla y|^2\right)^{1/2}\nonumber\\
\le&\eta \left(\int |y|^2|\nabla y|^2\right)+C_\eta \|y\|_4^{4}+C_\eta \|h^\varepsilon\|_{4,\partial\mathcal O}^{2}\|y\|_4^{2}.\label{3.33}
\end{align}

The boundary term from $-\mu\Delta y$ must also be retained. From the boundary condition, one has
\[
 [\partial_\nu y]_{\rm tan}=(M_w-\beta_\tau/\mu)y+h^\varepsilon/\mu,
\]
which implies
\begin{align*}
 &\int(-\mu\Delta y)\cdot|y|^2y
 =\mu \left(\int |y|^2|\nabla y|^2\right)+2\mu\sum_i\!\int(y\cdot\partial_i y)^2
   -\int_{\partial\mathcal O}|y|^2y\cdot((\mu M_w-\beta_\tau)y+h^\varepsilon),\\
& \left|\int_{\partial\mathcal O}|y|^2y\cdot((\mu M_w-\beta_\tau)y+h^\varepsilon)\right| \le C\|y\|_{4,\partial\mathcal O}^4+\|h^\varepsilon\|_{4,\partial\mathcal O}\|y\|_{4,\partial\mathcal O}^3\\
&\qquad\qquad\qquad\qquad\qquad\qquad \qquad ~  \le\eta \left(\int |y|^2|\nabla y|^2\right)+C_\eta \|y\|_4^{4}+C_\eta \|h^\varepsilon\|_{4,\partial\mathcal O}^{2}\|y\|_4^{2}.
\end{align*}
By Young's inequality, we obtain
\begin{align*}
\|h^\varepsilon\|_{4,\partial\mathcal O}
 \|y\|_4^{3/2}
 \left(\int_{\mathcal O}|y|^2|\nabla y|^2\right)^{3/8} &\le
\eta\int_{\mathcal O}|y|^2|\nabla y|^2
+C_\eta
 \|h^\varepsilon\|_{4,\partial\mathcal O}^{8/5}
 \|y\|_4^{12/5}\\
&\le
\eta\int_{\mathcal O}|y|^2|\nabla y|^2
+C_\eta
 \|h^\varepsilon\|_{4,\partial\mathcal O}^{2}
 \|y\|_4^{2}
+C_\eta\|y\|_4^{4}.
\end{align*}
Combining with \eqref{3.33}, one can find
\begin{equation}\label{3.34}
 \int(-\mu\Delta y+\nabla\pi_S)\cdot|y|^2y
 \ge \gamma_{\rm E}\left(\int |y|^2|\nabla y|^2\right)-C\|y\|_4^{4}-C\|h^\varepsilon\|_{4,\partial\mathcal O}^{2}\|y\|_4^{2}.
\end{equation}

For the pressure $\pi_N$ arising from the nonlinear
convection term, we have
\[
-\int_{\mathcal O}\pi_N\Delta\phi\,dx
=\varepsilon\sum_{i,j=1}^d
\int_{\mathcal O}y_i y_j\partial_i\partial_j\phi\,dx,
\qquad \partial_\nu\phi=0.
\]
Applying the Neumann elliptic estimate and the
duality between $L^4$ and $L^{4/3}$, we obtain
\[
\|\pi_N\|_4
\le C\varepsilon\|y\|_8^2.
\]
By the Gagliardo-Nirenberg inequality, we further have
\[
\|y\|_8^2
\le C\|y\|_4^{1/2}
\left(\int_{\mathcal O}|y|^2|\nabla y|^2\,dx
\right)^{3/8}
+C\|y\|_4^2.
\]
Consequently,
\begin{align*}
\left|\int_{\mathcal O}
\pi_N\operatorname{div}(|y|^2y)\,dx\right|
&\le 2\|\pi_N\|_4\|y\|_4
\left(\int_{\mathcal O}
|y|^2|\nabla y|^2\,dx\right)^{1/2}\\
&\le C\varepsilon\|y\|_4^{3/2}
\left(\int_{\mathcal O}
|y|^2|\nabla y|^2\,dx\right)^{7/8} +C\varepsilon\|y\|_4^3
\left(\int_{\mathcal O}
|y|^2|\nabla y|^2\,dx\right)^{1/2}.
\end{align*}
By Young's inequality, we deduce
\begin{equation}\label{3.35}
\begin{aligned}
\left|\int_{\mathcal O}
\pi_N\operatorname{div}(|y|^2y)\,dx\right|
&\le \eta\varepsilon
\int_{\mathcal O}|y|^2|\nabla y|^2\,dx +C_\eta\varepsilon
\left(\|y\|_4^{12}+\|y\|_4^6\right).
\end{aligned}
\end{equation}

Finally, we estimate the pressure contribution
associated with the forcing term
$F_u+\varepsilon r e_d$.
Let $\pi_F$ be the zero-mean solution of
\[
\begin{cases}
\Delta\pi_F
=\operatorname{div}(F_u+\varepsilon r e_d),
&\text{in }\mathcal O,\\
\partial_\nu\pi_F
=(F_u+\varepsilon r e_d)\cdot\nu,
&\text{on }\partial\mathcal O,\\
\displaystyle\int_{\mathcal O}\pi_F\,dx=0.
\end{cases}
\]
By the Neumann elliptic estimate and
Hölder's inequality, one has
\begin{align*}
\left|\int_{\mathcal O}
\nabla\pi_F\cdot|y|^2y\,dx\right|
&\le \|\nabla\pi_F\|_4\|y\|_4^3 \le C\bigl(\|F_u\|_4+\varepsilon\|r\|_4\bigr)
\|y\|_4^3.
\end{align*}

Combining the preceding estimates, for the total pressure $\pi=-\pi_L+\varepsilon\pi_S+\pi_N+\pi_F$, we obtain
\begin{align*}
\left|\int_{\mathcal O}
\nabla\pi\cdot |y|^2y\,dx\right|
&\le \eta\varepsilon
\int_{\mathcal O}|y|^2|\nabla y|^2\,dx +C\big(
\|\nabla V^\varepsilon\|_\infty
+\|\operatorname{div}V^\varepsilon\|_\infty
+\varepsilon\big)\|y\|_4^4 +C_\eta\varepsilon
\big(\|y\|_4^{12}+\|y\|_4^6\big)\\
&\quad+C_\eta\varepsilon
\|h^\varepsilon\|_{4,\partial\mathcal O}^2
\|y\|_4^2+C\big(\|F_u\|_4+\varepsilon\|r\|_4\big)
\|y\|_4^3.
\end{align*}

We are now ready to study the $L^4$ estimate for the velocity.  
Testing \eqrow{3.21}{3} by $|y|^2y$, and using \eqref{3.30}, \eqref{3.34}, \eqref{3.35}, one has
\begin{equation*}
\begin{aligned}
 \frac{\dd}{\dd t}\|y\|_4^4+\gamma_0\varepsilon\int|y|^2|\nabla y|^2 &\le C(\|\nabla V^\varepsilon\|_\infty+\|\Div V^\varepsilon\|_\infty+\varepsilon)\|y\|_4^4
       +C\varepsilon(\|y\|_4^{12}+\|y\|_4^6)\\
 &\quad+C(\|F_u\|_4+\varepsilon\|r\|_4)\|y\|_4^3
       +C\varepsilon\|h^\varepsilon\|_{4,\partial\OO}^2\|y\|_4^2.
\end{aligned}
\end{equation*}
By using the fact $\|y\|_4^{4}\le \|y\|_4^{2}+\|y\|_4^{10}$, we obtain
\begin{equation}\label{3.36}
\begin{aligned}
 \frac{\dd}{\dd t}\|y\|_4^2
 &\le C(\|\nabla V^\varepsilon\|_\infty+\|\Div V^\varepsilon\|_\infty+\varepsilon)\|y\|_4^2
       +C\varepsilon\|y\|_4^{10}\\
 &\quad+C\varepsilon(\|r\|_4^4+\|q\|_{W^{1,4}}^4)^{1/2}
       +C\|F_u\|_4\|y\|_4+C\varepsilon\|h^\varepsilon\|_{4,\partial\OO}^2.
\end{aligned}
\end{equation}

\begin{proposition} \label{prop 3.2}
Assume that the prescribed profiles, control coefficients and source terms are piecewise smooth, satisfy the boundary compatibility conditions in Lemma~\ref{lemma 3.3}, and  \eqref{3.23} holds.
Define the defect size
\begin{equation}\label{3.37}
 d_\varepsilon=
 \left(\int_0^{S_\varepsilon}(\|R_n\|_4+\|F_u\|_4)\dd t\right)^2
 +\varepsilon^{-1}\int_0^{S_\varepsilon}\|R_c\|_4^2\dd t
 +\varepsilon\int_0^{S_\varepsilon}\|h^\varepsilon\|_{4,\partial\mathcal O}^2\dd t.
\end{equation}
If $d_\varepsilon$ is sufficiently small, then the solution continues to $S_\varepsilon$ and
\begin{equation}\label{3.38}
 \sup_{0\le t\le S_\varepsilon}
 (\|r(t)\|_4^2+\|q(t)\|_{W^{1,4}}^2+\|y(t)\|_4^2)
 \le C d_\varepsilon.
\end{equation}
The constant depends on the fixed domain, coefficients, $C_{\rm prof}$ and $T_{\rm app}$, but not on $\varepsilon$.
\end{proposition}
\begin{proof}
We first derive the following estimates under the
bootstrap assumption
\[
\left(\|r\|_4^4+\|q\|_{W^{1,4}}^4\right)^{1/2}
+\|y\|_4^2\leq 1.
\]
Under this assumption, we have
\[
\begin{aligned}
&\|y\|_4^8
\left(\|r\|_4^4+\|q\|_{W^{1,4}}^4\right)^{1/2}
+\|y\|_4^{10} \leq
\left(\|r\|_4^4+\|q\|_{W^{1,4}}^4\right)^{1/2}
+\|y\|_4^2.
\end{aligned}
\]
Adding \eqref{3.29} and \eqref{3.36}, we have
\begin{align*}
 &\frac{\dd}{\dd t}\Big[(\|r\|_4^4+\|q\|_{W^{1,4}}^4)^{1/2}+\|y\|_4^2\Big] \\
 \le& C\big(\|\nabla V^\varepsilon\|_\infty+\|\Div V^\varepsilon\|_\infty
       +\|\alpha_n^\varepsilon\|_\infty+\|\alpha_c^\varepsilon\|_{W^{1,\infty}}+\varepsilon\big) \cdot\Big[(\|r\|_4^4+\|q\|_{W^{1,4}}^4)^{1/2}+\|y\|_4^2\Big]\\
 &+C(\|R_n\|_4+\|F_u\|_4)
      \Big[(\|r\|_4^4+\|q\|_{W^{1,4}}^4)^{1/2}+\|y\|_4^2\Big]^{1/2} +C\varepsilon^{-1}\|R_c\|_4^2
      +C\varepsilon\|h^\varepsilon\|_{4,\partial\OO}^2.
\end{align*}
By \eqref{3.23} and Young's inequality, we can get
\begin{align*}
 &\sup_{0\le s\le t}\Big[(\|r(s)\|_4^4+\|q(s)\|_{W^{1,4}}^4)^{1/2}+\|y(s)\|_4^2\Big]\\
 &\quad\le C\Big\{\sup_{0\le s\le t}
      \big[(\|r(s)\|_4^4+\|q(s)\|_{W^{1,4}}^4)^{1/2}+\|y(s)\|_4^2\big]\Big\}^{1/2}
       \int_0^t(\|R_n\|_4+\|F_u\|_4)\dd s\\
 &\qquad+C\varepsilon^{-1}\int_0^t\|R_c\|_4^2\dd s
      +C\varepsilon\int_0^t\|h^\varepsilon\|_{4,\partial\OO}^2\dd s\\
 &\quad\le\tfrac12\sup_{0\le s\le t}
      \Big[(\|r(s)\|_4^4+\|q(s)\|_{W^{1,4}}^4)^{1/2}+\|y(s)\|_4^2\Big]+Cd_\varepsilon.
\end{align*}
Therefore,
\[
 \sup_{0\le s\le t}\Big[(\|r(s)\|_4^4+\|q(s)\|_{W^{1,4}}^4)^{1/2}+\|y(s)\|_4^2\Big]
 \le2Cd_\varepsilon< \frac{1}{2}.
\]
Choosing \(d_\varepsilon>0\) sufficiently small, we close the bootstrap assumption
\[
(\|r\|_4^4+\|q\|_{W^{1,4}}^4)^{1/2}
+\|y\|_4^2\leq1.
\]
The resulting uniform \(L^4\) bound, together with the continuation criterion in Lemma~\ref{lemma 3.3} below, ensures that the solution extends to \(S_\varepsilon\). Finally, \eqref{3.22} yields \eqref{3.38}, which completes the proof.
\end{proof}

\subsection{Local strong solvability and continuation criterion}\label{subsection 3.5}
In this subsection, for each fixed \(\varepsilon>0\) and $L<\infty$, we establish the local strong solvability and a continuation criterion for the remainder system.

\begin{lemma} \label{lemma 3.3}
Assume that the prescribed coefficients and source terms in \(\eqref{3.21}\) are sufficiently smooth on \([0,L]\) and satisfy the corresponding boundary compatibility conditions. Suppose that \(\beta\) is smooth, symmetric and independent of time.
Compatibility of the initial data is understood in the trace spaces defining \(\mathcal Y\). Then, for every compatible initial datum \((r_0,q_0,y_0)\in\mathcal Y\), the remainder system admits a unique local strong solution. Moreover, if \(T_*\le L\) is the maximal existence time and
\[
\sup_{0\le t<T_*}
\left(
\|r(t)\|_{L^4}
+\|q(t)\|_{W^{1,4}}
+\|y(t)\|_{L^4}
\right)<\infty,
\]
then \(T_*=L\), and the solution extends continuously to \(t=L\) in \(\mathcal Y\).
\end{lemma}

\begin{proof}
By \eqref{2.4}, we can choose $\lambda_S>0$, such that for $v\in\Vv$
\begin{equation}\label{3.39}
 \gamma_{\rm E}\|v\|_{H^1}^2\le a_\beta(v,v)+\lambda_S\|v\|_2^2
 \le C\|v\|_{H^1}^2.
\end{equation}
For each $t$, let $(u_h,p_h)$ satisfy
\begin{equation}\label{3.40}
\left\{\begin{aligned}
 &\lambda_S u_h-\mu\Delta u_h+\nabla p_h=0,\quad\Div u_h=0,
       &&\text{in }\OO,\\
 &u_h\cdot\nu=0,\quad\NN u_h=h^\varepsilon,
       &&\text{on }\partial\OO,\\
 &\displaystyle\int_\OO p_h\,\dd x=0.
\end{aligned}\right.
\end{equation}
It follows from Lax-Milgram theorem and the tangential-stress regularity theorem
\cite{Saito2004} that there is a unique smooth lifting, with
\begin{equation}\label{3.41}
 \|u_h\|_{H^{j+2}}+\|p_h\|_{H^{j+1}}
 \le C_j\big(\|h^\varepsilon\|_{H^{j+1/2}(\partial\OO)}+\|u_h\|_2\big),
 \qquad j\in\mathbb N_0.
\end{equation}
We rewrite the Navier boundary condition as
\[
[2\mu D(u_h)\nu]_\tau
=h^\varepsilon-(\beta u_h)_\tau.
\]
Since \(\beta u_h\) contains no derivatives of \(u_h\), one can control this term by the trace inequality and absorb it into the elliptic regularity estimate.
Differentiating \eqref{3.40}, we obtain the same estimate as \eqref{3.41} for $\partial_t u_h$, with
$h_t^\varepsilon$ in place of $h^\varepsilon$.

Let $v=y-u_h$ and $b^\varepsilon=V^\varepsilon+\varepsilon u_h$. Then, \eqref{3.21} can be rewritten as
\begin{equation*}
\left\{\begin{aligned}
 &r_t-\varepsilon D_n\Delta r
=-b^\varepsilon\cdot\nabla r-\varepsilon v\cdot\nabla r+\alpha_n^\varepsilon r
   -\varepsilon(v+u_h)\cdot\nabla n^{\rm app}+R_n -\varepsilon\chi\Div(n^{\rm app}\nabla q+r\nabla c^{\rm app}+r\nabla q),\\
 &q_t-\varepsilon D_c\Delta q
=-b^\varepsilon\cdot\nabla q-\varepsilon v\cdot\nabla q+\alpha_c^\varepsilon q
   -\varepsilon(v+u_h)\cdot\nabla c^{\rm app}+R_c -\varepsilon(n^{\rm app}q+c^{\rm app}r+rq),\\
 &v_t-\varepsilon\mu\Delta v+\nabla\pi
=-b^\varepsilon\cdot\nabla v-\varepsilon(v\cdot\nabla)v
   -(v\cdot\nabla)V^\varepsilon-\varepsilon(v\cdot\nabla)u_h
   +\varepsilon r e_d+F_h, ~~~~ \Div v=0,
\end{aligned}\right.
\end{equation*}
with $\partial_\nu r=\partial_\nu q=0$, $v\cdot\nu=0$ and $\NN v=0$, where
\begin{equation*}
 F_h=F_u-\partial_t u_h+\varepsilon\mu\Delta u_h
       -(b^\varepsilon\cdot\nabla)u_h-(u_h\cdot\nabla)V^\varepsilon.
\end{equation*}
The scalar sources belong to $L^2L^2$ and $L^2H^1$, respectively,
and $F_h\in L^2L^2$.

By the Stokes maximal regularity estimate in \cite {G2006} and \cite {K2025}, for
$z_t-\varepsilon\mu\Delta z+\nabla p=f$, $\Div z=0$, with homogeneous
Navier conditions, one has
\begin{equation}\label{3.42}
 \|z\|_{\X_u(0,\tau)}+\|\nabla p\|_{L^2(0,\tau;L^2)}
 \le C_{\varepsilon,L}\big(\|z(0)\|_{H^1}+\|f\|_{L^2(0,\tau;L^2)}\big).
\end{equation}
The standard Neumann estimates give us
\begin{equation}\label{3.43}
 \begin{aligned}
 \|r\|_{\X_n(0,\tau)}
 &\le C_{\varepsilon,L}\big(\|r(0)\|_{H^1}+\|f_n\|_{L^2L^2}\big),\\
 \|q\|_{\X_c(0,\tau)}
 &\le C_{\varepsilon,L}\big(\|q(0)\|_{H^2}+\|f_c\|_{L^2H^1}\big).
 \end{aligned}
\end{equation}
The constants in
\eqref{3.42}-\eqref{3.43} are uniform for
$0<\tau\le\min\{1,L\}$.

By spatial interpolation and H\"older's inequality in time, we have 
\[
\|D^2q\|_{L^2(0,\tau;L^3)}
\le C\tau^{1/4}
\|q\|_{L^\infty(0,\tau;H^2)}^{1/2}
\|q\|_{L^2(0,\tau;H^3)}^{1/2}.
\]
 Hence, with all norms below taken
on $(0,\tau)$,
\begin{align}
& \|v\cdot\nabla r\|_{L^2L^2}
\le C\tau^{1/4}\|v\|_{\X_u}\|r\|_{\X_n}, 
~~~~ \|\nabla r\cdot\nabla q\|_{L^2L^2}+\|r\Delta q\|_{L^2L^2}
 \le C\tau^{1/4}\|r\|_{\X_n}\|q\|_{\X_c},\nonumber\\
& \|(v\cdot\nabla)v\|_{L^2L^2}
\le C\tau^{1/4}\|v\|_{\X_u}^2,
~~~~ \|v\cdot\nabla q\|_{L^2H^1}
 \le C\tau^{1/4}\|v\|_{\X_u}\|q\|_{\X_c}, 
 ~~~~\|rq\|_{L^2H^1}
 \le C\tau^{1/2}\|r\|_{\X_n}\|q\|_{\X_c}.\label{3.44}
\end{align}
The linear terms can be estimated by
\[
\|\operatorname{div}(n^{\rm app}\nabla q)\|_{L^2(0,\tau;L^2)}
\le C\tau^{1/2}\|q\|_{L^\infty(0,\tau;H^2)},
\]
and the remaining linear terms satisfy similar estimates. We choose \(M_*>0\) sufficiently large and then take \(\tau>0\) sufficiently small such that
\[
C_{\varepsilon,L}\bigl(\tau^{1/2}+M_*\tau^{1/4}\bigr)<\frac12.
\]
Then, the contraction mapping principle yields a unique local strong solution. The existence time depends on the initial data, the coefficient bounds and the $L^2$ norms of the lifted sources, including $\partial_t u_h$ when $h^\varepsilon\ne0$. When $h^\varepsilon=0$, it does not require time derivatives of the prescribed bulk sources.

Now, let us turn to the continuation.
Since $u_h$ is smooth, there is $B_*<\infty$ such that
\[
 \|r(t)\|_4+\|q(t)\|_{W^{1,4}}+\|v(t)\|_4\le B_*,
 \qquad 0\le t<T_*.
\]
Due to the Gagliardo-Nirenberg inequality, we have
\begin{equation}\label{3.45}
 \begin{aligned}
 \|\nabla r\|_4&\le C_{B_*}(1+\|r\|_{H^2}^{4/5}),&
 \|\nabla v\|_4&\le C_{B_*}(1+\|v\|_{H^2}^{4/5}),\\
 \|D^2q\|_4&\le C_{B_*}(1+\|q\|_{H^3}^{4/5}),&
 \|q\|_\infty&\le CB_*.
 \end{aligned}
\end{equation}
Next, we estimate the nonlinear terms arising in the higher order energy estimates.
For every $\eta>0$, by H\"older's and Young's inequalities, we have 
\begin{align}
 &|(v\cdot\nabla r,\Delta r)|+| (\nabla r\cdot\nabla q,\Delta r)|
 \le CB_*\|\nabla r\|_4\|r\|_{H^2}
 \le\eta\|r\|_{H^2}^2+C_{\eta,B_*},\nonumber\\
 & |(r\Delta q,\Delta r)|
\le C_{B_*}(1+\|q\|_{H^3}^{4/5})\|r\|_{H^2}
 \le\eta(\|r\|_{H^2}^2+\|q\|_{H^3}^2)+C_{\eta,B_*},
 \label{3.46}\\
 &\|\nabla(v\cdot\nabla q)\|_2\|\nabla\Delta q\|_2
\le C_{B_*}(1+\|v\|_{H^2}^{4/5}+\|q\|_{H^3}^{4/5})\|q\|_{H^3} \le\eta(\|v\|_{H^2}^2+\|q\|_{H^3}^2)+C_{\eta,B_*},
 \label{3.47}\\
& |(\nabla(rq),\nabla\Delta q)|
\le C_{B_*}(1+\|\nabla r\|_2)\|q\|_{H^3} \le\eta\|q\|_{H^3}^2+C_{\eta,B_*}(1+\|r\|_{H^1}^2).
 \label{3.48}
\end{align}
It is clear that in each case the high-order powers are strictly less than two and are
therefore absorbed by Young's inequality, without a smallness assumption
on $B_*$.

For the velocity, set
$v^\sharp=-\mu\Delta v+\nabla p_v+\lambda_Sv$, where $p_v$ is the
zero-mean Neumann pressure making $v^\sharp$ divergence free and tangent.
From the stationary Stokes estimate and \eqref{3.39}, we have
\[
 \|v^\sharp\|_2\le C\|v\|_{H^2},\qquad
 \|v\|_{H^2}\le C\|v^\sharp\|_2,
\]
which implies
\[
 |((v\cdot\nabla)v,v^\sharp)|
 \le C B_*\|\nabla v\|_4\|v\|_{H^2}
 \le\eta\|v\|_{H^2}^2+C_{\eta,B_*}.
\]
It follows from the symmetry and time independence of $\beta$ that
\[
 (v_t,v^\sharp)=\frac12\frac{\dd}{\dd t}
       \big[a_\beta(v,v)+\lambda_S\|v\|_2^2\big].
\]

We test \eqref{3.21}$_1$ by $r-\Delta r$, \eqref{3.21}$_2$ by
$q$ and $-\Delta q$ and its gradient equation by $-\nabla\Delta q$,
and \eqref{3.21}$_3$ by $v^\sharp$.
For the highest-order chemical estimate, from integration
by parts and $\partial_\nu q_t=0$, one has
\[
\begin{aligned}
-\int_{\mathcal O}\nabla q_t\cdot\nabla\Delta q\,dx
&=\int_{\mathcal O}\Delta q_t\,\Delta q\,dx
 -\int_{\partial\mathcal O}
       \partial_\nu q_t\,\Delta q\,dS =\frac12\frac{d}{dt}\|\Delta q\|_2^2.
\end{aligned}
\]
For the terms involving the prescribed coefficients,
we have
\begin{align*}
&\|\operatorname{div}(n^{\rm app}\nabla q)\|_2
\le C_{\varepsilon,L}\|q\|_{H^2}, 
~~~~\|b^\varepsilon\cdot\nabla q\|_{H^1}
\le C_{\varepsilon,L}\|q\|_{H^2}, 
~~~~\|v\cdot\nabla c^{\rm app}\|_{H^1}
\le C_{\varepsilon,L}\|v\|_{H^1}, \\
&\|\alpha_c^\varepsilon q-n^{\rm app}q-c^{\rm app}r\|_{H^1}
\le C_{\varepsilon,L}
  (\|q\|_{H^1}+\|r\|_{H^1}), ~~~~\|R_c-\varepsilon u_h\cdot\nabla c^{\rm app}\|_{H^1}
\le C_{\varepsilon,L}.
\end{align*}
Thus,
\begin{align*}
&\big\|-b^\varepsilon\cdot\nabla q+\alpha_c^\varepsilon q
-\varepsilon v\cdot\nabla c^{\rm app}
+R_c-\varepsilon u_h\cdot\nabla c^{\rm app}
-\varepsilon n^{\rm app}q-\varepsilon c^{\rm app}r\big\|_{H^1} 
\\ &\qquad\le C_{\varepsilon,L}
(1+\|r\|_{H^1}+\|q\|_{H^2}+\|v\|_{H^1}).
\end{align*}
For the nonlinear terms, by the interpolation estimates, we get
\begin{align*}
\|\nabla(v\cdot\nabla q)\|_2
&\le \|\nabla v\|_4\|\nabla q\|_4
   +\|v\|_4\|D^2q\|_4 \le C_{B_*}
(1+\|v\|_{H^2}^{4/5}+\|q\|_{H^3}^{4/5}),\\
\|\nabla(rq)\|_2
&\le \|q\|_\infty\|\nabla r\|_2
   +\|r\|_4\|\nabla q\|_4 \le C_{B_*}(1+\|r\|_{H^1}).
\end{align*}
By Young's inequality, for every $\eta>0$, we have
\begin{align*}
&\varepsilon\left|
\int_{\mathcal O}
\nabla(v\cdot\nabla q+rq)
\cdot\nabla\Delta q\,dx\right| \le
\eta(\|r\|_{H^2}^2+\|q\|_{H^3}^2
     +\|v\|_{H^2}^2)
+C_{\eta,\varepsilon,B_*}
(1+\|r\|_{H^1}^2).
\end{align*}
Therefore, we have
\begin{align*}
&\frac12\frac{d}{dt}
(\|q\|_2^2+\|\nabla q\|_2^2+\|\Delta q\|_2^2)
+\gamma_{\rm E}\varepsilon\|q\|_{H^3}^2\\
 &\quad\le
\eta(\|r\|_{H^2}^2+\|q\|_{H^3}^2
     +\|v\|_{H^2}^2) +C_{\eta,\varepsilon,B_*,L}
(1+\|r\|_{H^1}^2+\|q\|_{H^2}^2
     +\|v\|_{H^1}^2).
\end{align*}
Here,  
$\|\Div(n^{\rm app}\nabla q)\|_2\le C\|q\|_{H^2}$ and
$\|b^\varepsilon\cdot\nabla q\|_{H^1}+\|v\cdot\nabla c^{\rm app}\|_{H^1}
\le C(\|q\|_{H^2}+\|v\|_{H^1})$.
The lower-order nonlinear tests are bounded using
$\|v\cdot\nabla q\|_2+\|rq\|_2+\|r\nabla q\|_2\le C_{B_*}$.
Using \eqref{3.46}-\eqref{3.48},
the velocity bound above, and the Neumann elliptic estimates, we choose
$\eta$ after the fixed $\varepsilon$ and obtain
\begin{equation*}
 \begin{aligned}
 &\frac{\dd}{\dd t}\Big[\|r\|_{H^1}^2+\|q\|_2^2
       +\|\nabla q\|_2^2+\|\Delta q\|_2^2
       +a_\beta(v,v)+\lambda_S\|v\|_2^2\Big] +\gamma_\varepsilon
       \big(\|r\|_{H^2}^2+\|q\|_{H^3}^2+\|v\|_{H^2}^2\big)\\
 &\qquad\le C_{\varepsilon,B_*,L}
       \big(1+\|r\|_{H^1}^2+\|q\|_{H^2}^2+\|v\|_{H^1}^2\big),
 \end{aligned}
\end{equation*}
where the differentiated quantity is equivalent to
$\|r\|_{H^1}^2+\|q\|_{H^2}^2+\|v\|_{H^1}^2$ by
\eqref{3.39} and Neumann regularity. Using Gronwall's inequality, we obtain
\[
 \sup_{t<T_*}\big(\|r(t)\|_{H^1}^2+\|q(t)\|_{H^2}^2+\|v(t)\|_{H^1}^2\big)
 +\int_0^{T_*}\big(\|r\|_{H^2}^2+\|q\|_{H^3}^2+\|v\|_{H^2}^2\big)\dd t
 <\infty.
\]

Finally, \eqref{3.45} bounds the squared source norms by
\[
 C_{\varepsilon,B_*,L}\big(
 1+\|r\|_{H^1}^2+\|q\|_{H^2}^2+\|v\|_{H^1}^2
 +\|r\|_{H^2}^{8/5}+\|q\|_{H^3}^{8/5}+\|v\|_{H^2}^{8/5}\big).
\]
They are integrable since $T_*\le L<\infty$ and $8/5<2$.
The time continuity in \eqref{1.5} yields
$(r,q,v)(T_*)\in\Y$. If \(T_*<L\), the local existence result applied at \(t=T_*\) yields a time \(\delta>0\) such that the solution can be extended to \([0,T_*+\delta]\), where \(T_*+\delta\le L\). The same argument applies at each prescribed junction, using the common value of \(y=v+u_h\) as the initial velocity for the next interval. This contradicts the maximality of \(T_*\) and completes the proof.
\end{proof}

\begin{remark} Taking \(L=S_\varepsilon\), Lemma~\ref{lemma 3.3} and
\eqref{3.38} ensure the existence of the strong solution on \([0,S_\varepsilon]\). The continuation argument applies for fixed \(\varepsilon>0\), whereas the joint \(L^4\) estimate is uniform in \(\varepsilon\).
\end{remark} 

\subsection{Conclusion in the case without fluid layer}\label{subsection 3.6}
We now verify that the approximate profiles constructed
above satisfy the small defect condition required in
Proposition~\ref{prop 3.2}. Then, we can obtain the
terminal approximation from the joint remainder estimate.

We first estimate the scalar Neumann corrections.
Indeed, by the definition of $\ell_z^\varepsilon$ for $z\in\{n,c\}$
and the smoothness of the boundary coordinates,
for $j=0,1,2$ we have
\[
\|D^j\ell_z^\varepsilon(t)\|_{L^\infty(\mathcal O)}
\le C\varepsilon^{(1-j)/2}.
\]
Since $\ell_z^\varepsilon$ is supported in a boundary
neighborhood of thickness $O(\sqrt{\varepsilon})$,
we obtain
\[
\begin{aligned}
\|D^j\ell_z^\varepsilon(t)\|_4^4
\le C\varepsilon^{2(1-j)}
\int_{\partial\mathcal O}
\int_0^{C\sqrt{\varepsilon}}
ds\,dS \le C\varepsilon^{5/2-2j},
\end{aligned}
\]
which implies
\[
\|D^j\ell_z^\varepsilon(t)\|_4
\le C\varepsilon^{5/8-j/2},
\qquad j=0,1,2.
\]
The same argument applied to
$\partial_t\ell_z^\varepsilon$
yields
\[
\|\partial_t\ell_z^\varepsilon(t)\|_4
\le C\varepsilon^{5/8}.
\]
Since $V^\varepsilon$ is tangent, one has 
$
|V^\varepsilon\cdot\nabla \tilde{d}|
\le C\tilde{d}
 $ in the boundary neighborhood. 
Thus,
\[
\|V^\varepsilon\cdot\nabla\ell_z^\varepsilon\|_4
\le C\varepsilon^{5/8}.
\]
Using \eqref{3.11}, \eqref{3.16} and \eqref{3.17},
we obtain
\[
\|R_n(t)\|_4+\|R_c(t)\|_4
\le C\varepsilon^{5/8}\mathbf1_{(0,1)}(t).
\]
Since the scalar defects vanish for $t\ge 1$,
the preceding estimate gives
\begin{align*}
\|R_n\|_{L^1(0,S_\varepsilon;L^4)}
=\int_0^1\|R_n(t)\|_4\,dt
\le C\varepsilon^{5/8}, 
~~~~\|R_c\|_{L^2(0,S_\varepsilon;L^4)}^2
=\int_0^1\|R_c(t)\|_4^2\,dt
\le C\varepsilon^{5/4}.
\end{align*}
Then, by \eqref{3.19}, we have
\begin{align*}
&\|F_u\|_{L^1(0,S_\varepsilon;L^4)}
=\varepsilon\int_0^1
\|n^{\rm app}e_d+\mu\Delta z^\varepsilon
-(z^\varepsilon\cdot\nabla)z^\varepsilon\|_4\,dt\\ &\qquad\le C\varepsilon\int_0^1
\big(
\|n^{\rm app}\|_4+\|z^\varepsilon\|_{W^{2,4}}
+\|z^\varepsilon\|_\infty
 \|\nabla z^\varepsilon\|_4
\big)\,dt \le C\varepsilon
\end{align*}
and similarly, by \eqref{3.20}, we can get
\begin{align*}
&\varepsilon
\|h^\varepsilon\|_{
L^2(0,S_\varepsilon;L^4(\partial\mathcal O))}^2 =\varepsilon\int_0^1
\|\mathcal N_\beta z^\varepsilon\|_{
L^4(\partial\mathcal O)}^2\,dt \le C\varepsilon\int_0^1
\|z^\varepsilon\|_{W^{2,4}(\mathcal O)}^2\,dt \le C\varepsilon.
\end{align*}

We next verify the smallness condition in
Proposition~\ref{prop 3.2}.
Substituting the preceding estimates into
\eqref{3.37}, we obtain
\begin{align*}
d_\varepsilon
&=
\left(\int_0^{S_\varepsilon}
(\|R_n\|_4+\|F_u\|_4)\,dt\right)^2 +\varepsilon^{-1}
\int_0^{S_\varepsilon}\|R_c\|_4^2\,dt
+\varepsilon\int_0^{S_\varepsilon}
\|h^\varepsilon\|_{4,\partial\mathcal O}^2\,dt \le C\left(
\varepsilon^{5/4}
+\varepsilon^{1/4}
+\varepsilon
\right)
\le C\varepsilon^{1/4}.
\end{align*}
Thus, one can see $d_\varepsilon\to0$ as $\varepsilon\to0$,
and Proposition~\ref{prop 3.2} applies.
It follows that
\begin{equation}\label{3.49}
\sup_{0\le t\le S_\varepsilon}
\left(
\|r(t)\|_4
+\|q(t)\|_{W^{1,4}}
+\|y(t)\|_4
\right)
\le C\varepsilon^{1/8}.
\end{equation}

Finally, at $t=S_\varepsilon$, the approximate
profiles satisfy
$
V^\varepsilon=\varepsilon\bar u^\varepsilon, n^{\rm app}=\bar n^\varepsilon, c^{\rm app}=\bar c^\varepsilon.
$
Hence, from \eqref{3.49}, we find the terminal
approximation of the prescribed reference trajectory.
The nonnegativity of the scalar states follows from
Lemma~\ref{lemma 2.2} and the multiplicative form
of the scalar controls.

 \section{Global approximation with a weak fluid boundary layer}\label{section 4}
For a general symmetric friction matrix \(\beta\), the leading return flow \(U^0\) does not necessarily satisfy the Navier boundary condition. 
Therefore, in this section we introduce a boundary layer corrector to compensate for this mismatch and construct an approximate velocity that compensates for the leading
Navier boundary mismatch.
The scalar profiles $n^{\rm app},c^{\rm app}$ and the gains $\alpha_n^\varepsilon,\alpha_c^\varepsilon$ are those of Section~\ref{section 3}. Let
\[
 \delta=\mu\varepsilon,\qquad S_\varepsilon=T_{\rm app}/\varepsilon,\qquad
 0<\varepsilon<\min\{1,\mu^{-1},T_{\rm app}/2\}.
\]
All constants in this section are independent of $\varepsilon$.

\subsection{Boundary layer profile}\label{subsection 4.1} 
Use the distance $\tilde{d}(x)$ to $\partial\OO$ in a neighborhood
$V_\eta=\{x\in\overline\OO:\tilde{d}(x)<\eta\}$. In this collar, extend the friction matrix by $\beta(x)=\beta(\tilde\pi(x))$, where $\tilde\pi$ is the nearest boundary projection. Notice that $\nu=-\nabla \tilde{d}$, $|\nu|=1$ and $\partial_\nu\nu=0$. As in \cite{TYY2025}, let
\[
 U_b^0=-\frac{U^0\cdot\nu}{\tilde{d}}
 =\int_0^1\partial_\nu U^0(t,y-\lambda \tilde{d}\nu(y))\cdot\nu(y)\dd\lambda,
 \qquad x=y-\tilde{d}\nu(y).
\]
Thus $U_b^0$ extends smoothly to $\tilde{d}=0$. Denote $\Phi(t,t',x)=\Phi(t,\Phi(t',\cdot)^{-1}(x))$ for the two time flow. There is $K>0$ such that
\begin{equation*}
 |U^0\cdot\nabla \tilde{d}|\le K\tilde{d},\qquad
 \tilde{d}(\Phi(t,t',x))\le e^{K|t-t'|}\tilde{d}(x).
\end{equation*}
Choose $\rho>0$ and a smooth collar cutoff $\zeta_{\rm bl}$ satisfying
\begin{equation}\label{4.1}
 e^{3K}\rho<\eta/2,\qquad
 \zeta_{\rm bl}=1\quad(\tilde{d}\le\rho/3),\qquad
 0\le\zeta_{\rm bl}\le1,\qquad
 \supp\zeta_{\rm bl}\subset\{\tilde{d}<\rho\},
\end{equation}
and define
\begin{equation}\label{4.2}
 g^0=\zeta_{\rm bl}[2D(U^0)\nu+(\beta/\mu)U^0]_{\tr},
 \qquad \mu g^0|_{\partial\OO}=\NN U^0.
\end{equation}
Let $Q_+=(0,\infty)\times\OO\times\mathbb R_+$ and $\Sigma_+=(0,\infty)\times\OO\times\{0\}$. The boundary profile is the solution of
\begin{equation}\label{4.3}
\left\{\begin{aligned}
 &\begin{aligned}[b] &v_t+[(U^0\cdot\nabla_x)v+(v\cdot\nabla_x)U^0]_{\tr} + U_b^0\zeta\partial_\zeta v-\partial_{\zeta\zeta}v=\xi^v,\end{aligned} &&\text{in }Q_+,\\[2pt]
 &\partial_\zeta v(t,x,0)=g^0(t,x), &&\text{on }\Sigma_+,\\[2pt]
 &v(0,x,\zeta)=0, &&\text{in }\OO\times\mathbb R_+.
\end{aligned}\right.
\end{equation}
For a scalar or vector profile $f(t,x,\zeta)$, let $f^\delta(t,x)=f(t,x,\tilde{d}(x)/\sqrt\delta)$. 
The control source $\xi^v$ will be chosen tangent and supported in $(0,1)\times(\overline\OO\setminus\overline\Omega)\times\mathbb R_+$. Moreover, $v\cdot\nu=0$, and
\[
 (v^\delta\cdot\nabla)v^\delta
 =[(v\cdot\nabla_x)v]^\delta
   -\delta^{-1/2}(v^\delta\cdot\nu)(v_\zeta)^\delta
 =[(v\cdot\nabla_x)v]^\delta.
\]

For $t\ge1$, system \eqref{4.3} reduces to the heat equation on the half-line
\begin{equation}\label{4.4}
\left\{\begin{aligned}
 &v_t-v_{\zeta\zeta}=0, &&\text{in }(1,\infty)\times\OO\times\mathbb R_+,\\[2pt]
 &v_\zeta(t,x,0)=0, &&\text{on }(1,\infty)\times\OO\times\{0\}.
\end{aligned}\right.
\end{equation}
For nonnegative integers $p,s,m_0$, we introduce the weighted Sobolev norms
\begin{equation*}
 \|f\|_{H_x^p(H_\zeta^{s,m_0})}^2
 =\sum_{|\alpha|\le p}\sum_{j=0}^s\int_\OO\int_0^\infty
 (1+\zeta^2)^{m_0}|\partial_x^\alpha\partial_\zeta^j f|^2\dd\zeta\dd x.
\end{equation*}
The return flow is smooth, tangent and flushing, and the boundary datum in
\eqref{4.3} is smooth and vanishes near $t=0,1$.
Thus Lemma 7 in \cite{CMS2020}, with $n=10$ prepared moments and $(p,s)=(7,6)$ in its notation, gives the decay
below. The spatial localization of the
boundary layer follows from Section 3.4 in \cite{CMS2020}
and the choice of $\rho$ in \eqref{4.1}.
\begin{lemma} \label{lemma 4.1}
For the return flow $U^0$, there is a smooth exterior source $\xi^v$, supported in $0<t<1$, such that the solution of \eqref{4.3} satisfies
\begin{equation}\label{4.5}
 v\cdot\nu=0,\qquad
 \supp_xv(t)\subset\{\tilde{d}\le e^{3K}\rho\}\subset\subset V_\eta\qquad(t\ge0),
\end{equation}
and
\begin{equation*}
 \|v(t)\|_{H_x^7(H_\zeta^{6,m_0})}
 \le C_{m_0}\left(\frac{\ln(2+t)}{2+t}\right)^{1/4+5-m_0/2}
 \qquad(0\le m_0\le10).
\end{equation*}
The constants are independent of $t$ and $\varepsilon$.
\end{lemma}
 
Taking $m_0=6$, one has
\begin{equation*}
 \|v(t)\|_{H_x^7(H_\zeta^{6,6})}\le Cm(t),\qquad
 m(t)=\begin{cases}1,&0\le t\le1,\\(1+t)^{-2},&t>1.\end{cases}
\end{equation*}
By Sobolev embeddings in $x$ and $\zeta$, for $|\alpha|,j\le5$, we can find
\begin{equation}\label{4.6}
 \sup_x\|(1+\zeta)^6\partial_x^\alpha\partial_\zeta^jv\|_{L^2_\zeta}
 +\sup_{x,\zeta}(1+\zeta)^6|\partial_x^\alpha\partial_\zeta^jv|
 \le Cm(t).
\end{equation}
and
\[
\|\zeta^j\partial_x^\alpha\partial_\zeta^k v(t)\|
_{L^p_\zeta L^\infty_x}
\le C_{\alpha,k,j,p}m(t),
\qquad j=0,1,2,\quad p=2,4.
\]

\subsection{Definitions and estimates of the auxiliary profiles}\label{subsection 4.2}
Direct computations show that
\begin{align}
 &\partial_i f^\delta
=(\partial_{x_i}f)^\delta-\delta^{-1/2}\nu_i(\partial_\zeta f)^\delta,
       \qquad 1\le i\le d,\label{4.7}\\
 &\Div f^\delta
=(\Div_xf)^\delta-\delta^{-1/2}(\partial_\zeta f\cdot\nu)^\delta,
       \label{4.8}\\
 & \Delta f^\delta
=(\Delta_xf)^\delta
      -2\delta^{-1/2}(\nu\cdot\nabla_x\partial_\zeta f)^\delta
      +\delta^{-1/2}(\Delta \tilde{d})(\partial_\zeta f)^\delta
      +\delta^{-1}(\partial_{\zeta\zeta}f)^\delta.
       \label{4.9}
\end{align}
For vector $f$, on the boundary
\begin{equation}\label{4.10}
 \NN(f^\delta)
 =[2\mu D_x(f)(t,x,0)\nu+\beta f(t,x,0)]_{\tr}
       -\mu\delta^{-1/2}[\partial_\zeta f(t,x,0)]_{\tr}.
\end{equation}
By the support property \(\eqref{4.5}\), the boundary layer profiles extend smoothly by zero outside the boundary neighborhood.

\paragraph{Boundary pressure.}
Define
\begin{equation}\label{4.11}
 Q(t,x,\zeta)
 =-\int_\zeta^\infty
   \big((v\cdot\nabla_x)U^0\cdot\nu
           -v\cdot(U^0\cdot\nabla_x)\nu\big)(t,x,s)\dd s.
\end{equation}
Then
\begin{align*}
 &Q_\zeta
=(v\cdot\nabla_x)U^0\cdot\nu-v\cdot(U^0\cdot\nabla_x)\nu, 
 (U^0\cdot\nabla_x)v\cdot\nu
 =-v\cdot(U^0\cdot\nabla_x)\nu,\\
 & v_t+(U^0\cdot\nabla_x)v+(v\cdot\nabla_x)U^0
       +U_b^0\zeta v_\zeta-v_{\zeta\zeta}-\xi^v
=Q_\zeta\nu,\\
 &\nabla(\delta Q^\delta)
=\delta(\nabla_xQ)^\delta-\sqrt\delta(Q_\zeta)^\delta\nu.
\end{align*}
 Therefore, the normal component of the
boundary layer equation is canceled by the
pressure correction. Moreover, $Q=0$ for $t\ge1$, since $U^0=0$
on this interval.

\paragraph{Second boundary corrector.}
We introduce the second boundary corrector
\begin{equation}\label{4.12}
 w(t,x,\zeta)
 =-e^{-\zeta}[2D_x(v)(t,x,0)\nu+(\beta/\mu)v(t,x,0)]_{\tr}
       -\nu(x)\int_\zeta^\infty\Div_xv(t,x,s)\dd s,
\end{equation}
which implies
\begin{equation}\label{4.13}
\left\{\begin{aligned}
 &\partial_\zeta w\cdot\nu=\Div_xv, &&\text{in }Q_+,\\[2pt]
 &[\partial_\zeta w(t,x,0)]_{\tr}=[2D_x(v)(t,x,0)\nu+(\beta/\mu)v(t,x,0)]_{\tr}, &&\text{on }\Sigma_+.
\end{aligned}\right.
\end{equation}
For $j=0$ and $|\alpha|\le4$, it follows from \eqref{4.6} that
\[
\begin{aligned}
 \left|
 \int_\zeta^\infty
 \partial_x^\alpha\operatorname{div}_x v(t,x,s)\,ds
 \right|
 &\le
 \left(\int_\zeta^\infty(1+s)^{-12}\,ds\right)^{1/2}
 \left\|
 (1+s)^6\partial_x^\alpha\operatorname{div}_x v(t,x,s)
 \right\|_{L^2(0,\infty)}
 \\ &\le C(1+\zeta)^{-11/2}m(t).
\end{aligned}
\]
For $1\le j\le4$, we use
\[
 \partial_\zeta^j
 \int_\zeta^\infty
 \partial_x^\alpha\operatorname{div}_x v(t,x,s)\,ds
 =
 -\partial_x^\alpha\partial_\zeta^{j-1}
 \operatorname{div}_x v(t,x,\zeta).
\]
By \eqref{4.6}, \eqref{4.11} and \eqref{4.12}, we obtain, for $|\alpha|\le4$ and $0\le j\le4$,
\begin{equation*}
 |\partial_x^\alpha\partial_\zeta^j w|
   +|\partial_x^\alpha\partial_\zeta^j Q|
 \le C(1+\zeta)^{-11/2}m(t).
\end{equation*}
For $t\ge1$, using \eqrow{4.4}{1}, we have
\[
 w_t=-e^{-\zeta}
 [2D_x(v_{\zeta\zeta})(t,x,0)\nu
                +(\beta/\mu)v_{\zeta\zeta}(t,x,0)]_{\tr}
 -\nu\int_\zeta^\infty\Div_xv_{\zeta\zeta}(t,x,s)\dd s.
\]
The same bound on $0\leq t\leq1$ follows directly
from the smooth prepared construction in Lemma \ref{lemma 4.1}. Therefore
\[
 |w_t(t,x,\zeta)|\le C(1+\zeta)^{-11/2}m(t).
\]
In addition, $w,Q,w_t$ have the same slow support as in \eqref{4.5}.

\paragraph{Evaluation at the fast variable.}
In the normal coordinates $x=y-r\nu(y)$, the Jacobian is bounded above and below. For $p=2,4$,
\begin{align*}
 \|f^\delta\|_p^p
 &\le C\int_{\partial\OO}\int_0^\eta
       |f(t,y-r\nu(y),r/\sqrt\delta)|^p\dd r\dd S_y \le C\sqrt\delta\int_0^\infty
                 \sup_{x\in V_\eta}|f(t,x,\zeta)|^p\dd\zeta.
\end{align*}
For the derivatives of $v,w,Q$ above, for $w_t$, and for the same terms multiplied by $\zeta$ or $\zeta^2$, the integrals are bounded by $Cm(t)^p$. Hence
\begin{equation}\label{4.14}
 \|f^\delta\|_4\le C\delta^{1/8}m(t),\qquad
 \|f^\delta\|_2\le C\delta^{1/4}m(t),\qquad
 \|f^\delta\|_\infty\le Cm(t).
\end{equation}
If $f,g$ are two such profiles, one has
\[
 \|(fg)^\delta\|_4
 \le\|f^\delta\|_\infty\|g^\delta\|_4
 \le C\delta^{1/8}m(t)^2.
\]

\paragraph{Interior potential corrector.}
Define the mean-zero solution
\begin{equation}\label{4.15}
\left\{\begin{aligned}
 &\Delta\theta^\delta=-(\Div_xw)^\delta, &&\text{in }\OO,\\[2pt]
 &\partial_\nu\theta^\delta=-w(t,x,0)\cdot\nu, &&\text{on }\partial\OO,\\[2pt]
 &\displaystyle\int_\OO\theta^\delta(t,x)\,\dd x=0, &&t\ge0.
\end{aligned}\right.
\end{equation}
It follows from \eqref{4.8} and \eqref{4.13} that
\begin{align*}
 \Div(\sqrt\delta v^\delta+\delta w^\delta)
 &=\sqrt\delta(\Div_xv)^\delta+\delta(\Div_xw)^\delta
       -\sqrt\delta(\partial_\zeta w\cdot\nu)^\delta
 =\delta(\Div_xw)^\delta,\\
 (\sqrt\delta v^\delta+\delta w^\delta)\cdot\nu|_{\partial\OO}
 &=\delta w(t,x,0)\cdot\nu,\\
 \int_\OO-(\Div_xw)^\delta\dd x
 &=-\int_{\partial\OO}w(t,x,0)\cdot\nu\dd S.
\end{align*}
The compatibility condition is thus satisfied, and by the standard Neumann theory, there exists a unique solution \(\theta^\delta\) with zero mean.
Moreover,
\begin{align}
 \|\theta^\delta\|_{W^{2,4}}
 &\le C\left(\|(\Div_xw)^\delta\|_4
       +\|w(t,\cdot,0)\cdot\nu\|_{W^{3/4,4}(\partial\OO)}\right)
 \le Cm(t),\label{4.16}\\
 \|\theta^\delta\|_{H^{j+2}}
 &\le C_j\left(\|(\Div_xw)^\delta\|_{H^j}
       +\|w(t,\cdot,0)\cdot\nu\|_{H^{j+1/2}(\partial\OO)}\right)
       \quad(j=0,1,2).\nonumber
\end{align}
For $j=0,1,2$, it holds that
\[
 \|(\Div_xw)^\delta\|_{H^j}
 \le C\sum_{\ell=0}^j\delta^{-\ell/2+1/4}m(t),
 \qquad
 \|w(t,\cdot,0)\cdot\nu\|_{H^{j+1/2}(\partial\OO)}\le Cm(t),
\]
which gives
\begin{equation*}
 \|\theta^\delta\|_{H^2}\le Cm(t),\qquad
 \|\theta^\delta\|_{H^3}\le C\delta^{-1/4}m(t),\qquad
 \|\theta^\delta\|_{H^4}\le C\delta^{-3/4}m(t).
\end{equation*}
Thanks to the embedding theorem, we have
\begin{equation}\label{4.17}
 \|\nabla\theta^\delta\|_\infty+\|D^2\theta^\delta\|_4\le Cm(t),
 \qquad
 \delta\|D^2\theta^\delta\|_\infty\le C\delta^{1/4}m(t).
\end{equation}

Now, if we define
\begin{equation}\label{4.18}
 L^\varepsilon=\sqrt\delta\,v^\delta+\delta w^\delta
                         +\delta\nabla\theta^\delta,
\end{equation}
one has
\begin{equation*}
\left\{\begin{aligned}
 &\Div L^\varepsilon=\delta(\Div_xw)^\delta+\delta\Delta\theta^\delta=0, &&\text{in }Q^\varepsilon,\\[2pt]
 &L^\varepsilon\cdot\nu=\delta w(t,x,0)\cdot\nu+\delta\partial_\nu\theta^\delta=0, &&\text{on }\Sigma^\varepsilon.
\end{aligned}\right.
\end{equation*}
Moreover,
\begin{align}\label{4.19}\begin{aligned}
 &\|L^\varepsilon\|_4
\le\sqrt\delta\|v^\delta\|_4
       +\delta\|w^\delta\|_4+\delta\|\nabla\theta^\delta\|_4
 \le C(\delta^{5/8}+\delta^{9/8}+\delta)m(t) \le C\varepsilon^{5/8}m(t),\\
& \|L^\varepsilon\|_\infty
\le\sqrt\delta\|v^\delta\|_\infty+\delta\|w^\delta\|_\infty
       +\delta\|\nabla\theta^\delta\|_\infty
 \le C\sqrt\delta\,m(t),\\
& \|\nabla L^\varepsilon\|_\infty
\le\sqrt\delta\|(\nabla_xv)^\delta\|_\infty
       +\|(\partial_\zeta v)^\delta\|_\infty
       +\delta\|(\nabla_xw)^\delta\|_\infty +\sqrt\delta\|(\partial_\zeta w)^\delta\|_\infty
       +\delta\|D^2\theta^\delta\|_\infty
 \le Cm(t).
\end{aligned}
\end{align}

\paragraph{Complete expansion.}
With the profiles defined above, the prescribed velocity is
\begin{equation}\label{4.20}
 V^\varepsilon=
 \begin{cases}
 U^0+\varepsilon z^\varepsilon+L^\varepsilon,
       &0\le t<1,\\
 \varepsilon\bar u^\varepsilon+L^\varepsilon,
       &1\le t\le S_\varepsilon.
 \end{cases}
\end{equation}
We introduce the following expansion
\begin{equation}\label{4.21}
 \begin{aligned}
 N=\varepsilon(n^{\rm app}+r),~~~~ \mathcal{C}=\varepsilon(c^{\rm app}+q),~~~~
 U=V^\varepsilon+\varepsilon y,~~~~
 P=P_V+\varepsilon\pi,~~~~
 G_n=\alpha_n^\varepsilon N,~~~~
 G_c=\alpha_c^\varepsilon \mathcal{C}.
 \end{aligned}
\end{equation}
Here $V^\varepsilon$ is given by \eqref{4.20}, and
\begin{equation}\label{4.22}
 P_V=
 \begin{cases}
 p^0+\varepsilon p^{1,\varepsilon}+\delta\sigma^0
       +\delta Q^\delta
       -\delta(\theta_t^\delta+U^0\cdot\nabla\theta^\delta),
       &0\le t<1,\\
 \varepsilon^2\bar p^\varepsilon-\delta\theta_t^\delta,
       &1<t\le S_\varepsilon.
 \end{cases}
\end{equation}
The pressure, at $t=1$, may be assigned arbitrarily,
since the equations are required to hold only
almost everywhere in time. Moreover, by \eqref{3.4}
and $U^0(1)=0$, we have
\[
 V^\varepsilon(1^-)
 =\varepsilon z^\varepsilon(1)+L^\varepsilon(1)
 =\varepsilon\bar u^\varepsilon(1)+L^\varepsilon(1)
 =V^\varepsilon(1^+),
\]
which shows that $V^\varepsilon$ is continuous at $t=1$.

The fluid controls are
\begin{equation}\label{4.23}
 G_u=
 \begin{cases}
 \xi^0+\varepsilon\xi^{1,\varepsilon}
       +\sqrt\delta(\xi^v)^\delta,&t<1,\\
 \varepsilon^2\bar f_u^\varepsilon,&t>1,
 \end{cases}
 \qquad
 \sigma^\varepsilon=
 \begin{cases}
 \sigma^0,&t<1,\\
 0,&t>1.
 \end{cases}
\end{equation}
As for $P_V$, the values are immaterial at $t=1$.

\paragraph{Boundary stress and its estimate.}
It follows from \eqref{4.2}, \eqref{4.10} and \eqref{4.13} that
\begin{align*}
 \NN(U^0+L^\varepsilon)
 ={}&\NN U^0-\mu\partial_\zeta v(t,x,0) +\sqrt\delta[2\mu D_x(v)(t,x,0)\nu+\beta v(t,x,0)]_{\tr}
       -\mu\sqrt\delta[\partial_\zeta w(t,x,0)]_{\tr}\\
 &+\delta\NN\nabla\theta^\delta
       +\delta[2\mu D_x(w)(t,x,0)\nu+\beta w(t,x,0)]_{\tr}\\
 ={}&\delta\NN\nabla\theta^\delta
       +\delta[2\mu D_x(w)(t,x,0)\nu+\beta w(t,x,0)]_{\tr}.
\end{align*}
Since $\NN\bar u^\varepsilon=0$, the remainder stress is
\begin{equation}\label{4.24}
 \begin{aligned}
 h^\varepsilon:=-\varepsilon^{-1}\NN V^\varepsilon
 ={}&-\mathbf1_{(0,1)}\NN z^\varepsilon-\mu\NN\nabla\theta^\delta -\mu[2\mu D_x(w)(t,x,0)\nu+\beta w(t,x,0)]_{\tr}.
 \end{aligned}
\end{equation}
According to the trace theorem and the Sobolev embedding, we have
\begin{align*}
 & \|\mathbf1_{(0,1)}\NN z^\varepsilon\|_{4,\partial\OO}
\le C\mathbf1_{(0,1)}, ~~~~
 \|\NN\nabla\theta^\delta\|_{4,\partial\OO}
\le C\left(\|D^2\theta^\delta\|_{H^1}
                    +\|\nabla\theta^\delta\|_{H^1}\right)
 \le C\delta^{-1/4}m(t),\\
 &\|[2\mu D_x(w)(t,\cdot,0)\nu+\beta w(t,\cdot,0)]_{\tr}\|_{4,\partial\OO}
 \le Cm(t),
\end{align*}
which implies 
\begin{equation}\label{4.25}
 \begin{aligned}
 \varepsilon\int_0^{S_\varepsilon}
                   \|h^\varepsilon\|_{4,\partial\OO}^2\dd t
 \le C\varepsilon\left(1+(1+\delta^{-1/2})
                                \int_0^\infty m(t)^2\dd t\right) \le C\varepsilon(1+\delta^{-1/2})\le C\varepsilon^{1/2}.
 \end{aligned}
\end{equation}

For convenience in the fluid residual calculations,
we introduce the following auxiliary velocity,
pressure and forcing profiles:
\[
 \begin{aligned}
 K^\varepsilon&=U^0+\delta u^1+\sqrt\delta\,v^\delta
                          +\delta w^\delta+\delta\nabla\theta^\delta,\\
 P_K&=p^0+\delta p^1+\delta\sigma^0+\delta Q^\delta
                  -\delta(\theta_t^\delta+U^0\cdot\nabla\theta^\delta),\\
 G_K&=\xi^0+\delta\xi^1+\sqrt\delta(\xi^v)^\delta, 
 \end{aligned}
\]
where $
 u^1=z^\varepsilon/\mu, 
 p^1=p^{1,\varepsilon}/\mu, 
 \xi^1=\xi^{1,\varepsilon}/\mu
$ as in \(\eqref{3.4}\). Notice that for $t>1$,  $U^0,u^1,p^0,p^1,\sigma^0,\xi^0,\xi^1,\xi^v,Q$ are zero. Moreover, $V^\varepsilon=K^\varepsilon$ for $t<1$, and
$V^\varepsilon=\varepsilon\bar u^\varepsilon+K^\varepsilon$ for $t>1$.  
Thus, define the fluid residual by
\begin{equation}\label{4.26}
 F_0^\varepsilon=\varepsilon^{-1}
 \bigl[G_K-K_t^\varepsilon-(K^\varepsilon\cdot\nabla)K^\varepsilon
                  +\delta\Delta K^\varepsilon-\nabla P_K\bigr].
\end{equation}
Here $F_0^\varepsilon$ is defined separately on $(0,1)$ and
$(1,S_\varepsilon)$, with classical time derivatives on each interval.
The auxiliary field $K^\varepsilon$ may be discontinuous
at $t=1$, but it is not the approximate velocity itself.
The actual profile satisfies
\[
 V^\varepsilon(1^-)=\varepsilon z^\varepsilon(1)+L^\varepsilon(1)
 =\varepsilon\bar u^\varepsilon(1)+L^\varepsilon(1)
 =V^\varepsilon(1^+).
\]
Thus its distributional time derivative has no Dirac mass at $t=1$.

Now, let us collect the estimates for the prepared fluid layer in the following lemma.
\begin{lemma}\label{lemma 4.2}
Let $U^0,z^\varepsilon$ be the return flow and
first-order correction of Lemma~\ref{lemma 3.1}. The exterior source $\xi^v$ in \eqref{4.3}
can be chosen such that the profiles \eqref{4.12}, \eqref{4.15}
and \eqref{4.18} satisfy
\begin{align*}
& v\cdot\nu=0,\qquad v(0)=0,\qquad
 \|v(t)\|_{H_x^7(H_\zeta^{6,6})}\le Cm(t),\\
 &\Div L^\varepsilon=0,\qquad L^\varepsilon\cdot\nu=0,\qquad
 L^\varepsilon(0)=0,\\
 & \|L^\varepsilon(t)\|_p
\le C\varepsilon^{1/2+1/(2p)}m(t)\qquad(p=2,4),\\
 & \|L^\varepsilon(t)\|_\infty
\le C\varepsilon^{1/2}m(t),\qquad
 \|\nabla L^\varepsilon(t)\|_\infty\le Cm(t),\\
 & \|L^\varepsilon\|_{L^2(0,1;L^2)}
\le C\varepsilon^{3/4},\qquad
 \int_0^{S_\varepsilon}\|\nabla L^\varepsilon\|_\infty\dd t\le C.
\end{align*}
The stress is exactly
\begin{equation*}
 \begin{aligned}
 &\NN(U^0+L^\varepsilon)
=\delta\NN\nabla\theta^\delta
    +\delta[2\mu D_x(w)(t,x,0)\nu+\beta w(t,x,0)]_{\tr},\\
 & h^\varepsilon=-\varepsilon^{-1}\NN V^\varepsilon
=-\mathbf1_{(0,1)}\NN z^\varepsilon-\mu\NN\nabla\theta^\delta
    -\mu[2\mu D_x(w)(t,x,0)\nu+\beta w(t,x,0)]_{\tr}.
 \end{aligned}
\end{equation*}
\end{lemma}

\subsection{Equations and estimates for the remainder}\label{subsection 4.3}
Substituting the approximate profiles \eqref{4.21}-\eqref{4.23} into the rescaled system \eqref{3.1}, we obtain the remainder system  \eqref{3.21} with \(U=V^\varepsilon+\varepsilon y\). The scalar and fluid defects are defined by \eqref{3.16}, \eqref{3.17} and \eqref{3.18}, while the boundary defect is given by \eqref{4.24}. To apply Proposition~\ref{prop 3.2}, it remains to verify the uniform profile bounds \eqref{3.23}  and the smallness condition \eqref{3.37}. We now estimate the additional terms arising from the fluid boundary layer.

\paragraph{Scalar residuals.}
Let $R_n^{\rm sp},R_c^{\rm sp}$ be the scalar residuals associated with
$V^\varepsilon-L^\varepsilon$ in Section~\ref*{section 3} and recall that
\[
 \|R_n^{\rm sp}(t)\|_4+\|R_c^{\rm sp}(t)\|_4
 \le C\varepsilon^{5/8}\mathbf1_{(0,1)}(t).
\]
Notice that only the transport terms change when $L^\varepsilon$ is added in this section, that is
\begin{equation*}
 R_n=R_n^{\rm sp}-L^\varepsilon\cdot\nabla n^{\rm app},
 \qquad R_c=R_c^{\rm sp}-L^\varepsilon\cdot\nabla c^{\rm app}.
\end{equation*}
Using \eqref{3.12} and \eqref{4.19}, one has
\[
 \|L^\varepsilon\cdot\nabla n^{\rm app}\|_4
 +\|L^\varepsilon\cdot\nabla c^{\rm app}\|_4
 \le\|L^\varepsilon\|_4
       (\|\nabla n^{\rm app}\|_\infty+\|\nabla c^{\rm app}\|_\infty)
 \le C\varepsilon^{5/8}m(t).
\]
In particular, since $n^{\rm app}=\bar n^\varepsilon$ and $c^{\rm app}=\bar c^\varepsilon$
for $1<t<S_\varepsilon$, we get
\begin{equation*}
 R_n=-L^\varepsilon\cdot\nabla\bar n^\varepsilon,
 \qquad R_c=-L^\varepsilon\cdot\nabla\bar c^\varepsilon
 \qquad(1<t<S_\varepsilon).
\end{equation*}
Thus, by \eqref{3.12} and \eqref{4.19}, $\|R_n(t)\|_4+\|R_c(t)\|_4\le C\varepsilon^{5/8}m(t)$, and
\begin{equation}\label{4.27}
\begin{aligned}
 &\|R_n\|_{L^1(0,S_\varepsilon;L^4)}
       +\|R_c\|_{L^2(0,S_\varepsilon;L^4)} \le C\varepsilon^{5/8}
 \left[\int_0^\infty m(t)\dd t
       +\left(\int_0^\infty m(t)^2\dd t\right)^{1/2}\right]
 \le C\varepsilon^{5/8}.
\end{aligned}
\end{equation}

\paragraph{Fluid residuals.}
We estimate $F_0^\varepsilon$ directly from the auxiliary profiles
$K^\varepsilon,P_K,G_K$ in \eqref{4.26}.
The calculation is performed separately on $(0,1)$ and
$(1,S_\varepsilon)$.

It follows from \eqref{3.2} and
\eqref{3.4} that the leading and first-order terms cancel. 
Then, we can use the following identities to cancel the viscous term \(-\delta\Delta U^0\), the order-\(\sqrt{\delta}\) terms generated by the boundary layer \(v^\delta\), and the time derivative and linear convection terms involving \(\delta\nabla\theta^\delta\).
\begin{align*}
 &-\delta\Delta U^0+\nabla(\delta\sigma^0)=0,  ~~~~\nabla(\delta Q^\delta)
       =\delta(\nabla_xQ)^\delta-\sqrt\delta(Q_\zeta)^\delta\nu,\\
 &v_t+(U^0\cdot\nabla_x)v+(v\cdot\nabla_x)U^0
       +U_b^0\zeta v_\zeta-v_{\zeta\zeta}-\xi^v
       =Q_\zeta\nu,\\
 &\partial_t\nabla\theta^\delta
 +(U^0\cdot\nabla)\nabla\theta^\delta
 +(\nabla\theta^\delta\cdot\nabla)U^0
       =\nabla(\theta_t^\delta+U^0\cdot\nabla\theta^\delta).
\end{align*}
However, the boundary layer pressure correction leaves the term \(\delta(\nabla_xQ)^\delta\), which is included in the fluid residual estimated below.

Expanding the remaining terms by
\eqref{4.7}-\eqref{4.9}, and using
$v\cdot\nu=0$ and $\delta/\varepsilon=\mu$, we obtain
\begin{equation*}
\begin{aligned}
 \mu^{-1}F_0^\varepsilon={}&
 \big[(\Delta \tilde{d})v_\zeta
       -2(\nu\cdot\nabla_x)v_\zeta+w_{\zeta\zeta}\big]^\delta +\sqrt\delta\big[\Delta_xv+(\Delta \tilde{d})w_\zeta
                         -2(\nu\cdot\nabla_x)w_\zeta\big]^\delta\\
 &+\delta(\Delta_xw)^\delta
       +\delta\Delta u^1+\delta\Delta\nabla\theta^\delta -\big[\big(v+\sqrt\delta(w+u^1+\nabla\theta^\delta)\big)
        \cdot\nabla_x
        \big(v+\sqrt\delta(w+u^1+\nabla\theta^\delta)\big)\big]^\delta\\
 &-\big[(U^0\cdot\nabla_x)w+(w\cdot\nabla)U^0
                              +U_b^0\zeta w_\zeta\big]^\delta +\big[\big((w+u^1+\nabla\theta^\delta)\cdot\nu\big)
                              (v_\zeta+\sqrt\delta w_\zeta)\big]^\delta -(\nabla_xQ)^\delta-(w_t)^\delta.
\end{aligned}
\end{equation*}
Now, we estimate this identity with $L^p$ norms for $p=2,4$. The smooth coefficients
$U^0,u^1,U_b^0$ and their required derivatives are bounded by
$C\mathbf1_{(0,1)}(t)\le Cm(t)$. The evaluation estimate
\eqref{4.14} controls the linear boundary layer terms, such as
\[
\begin{aligned}
 \|((\Delta \tilde{d})v_\zeta)^\delta\|_p
 +\|((\nu\cdot\nabla_x)v_\zeta)^\delta\|_p
 +\|(w_{\zeta\zeta})^\delta\|_p +\|(\nabla_xQ)^\delta\|_p+\|(w_t)^\delta\|_p
 \le C\delta^{1/(2p)}m(t).
\end{aligned}
\]
The remaining evaluated diffusion terms have additional factors
$\sqrt\delta$ or $\delta$, while
$\delta\|\Delta u^1\|_p\le C\delta m(t)$.
For the global potential correction, by \eqref{4.15}, one has
\begin{align*}
 \delta\Delta\nabla\theta^\delta
 =-\delta(\nabla_x\Div_xw)^\delta
       +\sqrt\delta\,\nu(\partial_\zeta\Div_xw)^\delta, 
~~~~ \|\delta\Delta\nabla\theta^\delta\|_p
 \le C\big(\delta^{1+1/(2p)}+\delta^{1/2+1/(2p)}\big)m(t).
\end{align*}
We also have for the terms with $U^0$ 
\[
 \|[(U^0\cdot\nabla_x)w+(w\cdot\nabla)U^0
                         +U_b^0\zeta w_\zeta]^\delta\|_p
 \le C\delta^{1/(2p)}m(t)^2.
\]
And using \eqref{4.14} and
\eqref{4.17}, it holds that
\begin{align*}
 &\big\|\big[\big(v+\sqrt\delta(w+u^1+\nabla\theta^\delta)\big)
       \cdot\nabla_x
       \big(v+\sqrt\delta(w+u^1+\nabla\theta^\delta)\big)\big]^\delta\big\|_p\\
 &\quad\le Cm(t)\big[\|(\nabla_xv)^\delta\|_p
       +\sqrt\delta\|(\nabla_xw)^\delta\|_p
       +\sqrt\delta\|\nabla u^1\|_p
       +\sqrt\delta\|D^2\theta^\delta\|_p\big]\\
 &\quad\le C\big(\delta^{1/(2p)}+\sqrt\delta\big)m(t)^2
 \le C\delta^{1/(2p)}m(t)^2,\\
 &\big\|\big[\big((w+u^1+\nabla\theta^\delta)\cdot\nu\big)
                      (v_\zeta+\sqrt\delta w_\zeta)\big]^\delta\big\|_p\\
 &\quad\le
 (\|w^\delta\|_\infty+\|u^1\|_\infty
                   +\|\nabla\theta^\delta\|_\infty)
 (\|(v_\zeta)^\delta\|_p+\sqrt\delta\|(w_\zeta)^\delta\|_p)\\
 &\quad\le C\delta^{1/(2p)}m(t)^2.
\end{align*}
It follows from \(\eqref{4.16}\) that \(\|D^2\theta^\delta(t)\|_p\le Cm(t)\).
Since $0<\delta\le1$, we conclude that
\begin{equation}\label{4.28}
\begin{aligned}
 \|F_0^\varepsilon(t)\|_p
 &\le C\delta^{1/(2p)}\big(m(t)+m(t)^2\big),\\
 \|F_0^\varepsilon\|_{L^1(0,S_\varepsilon;L^p)}
 &\le C\delta^{1/(2p)}\int_0^\infty(m(t)+m(t)^2)\dd t
 \le C\varepsilon^{1/(2p)},\qquad p=2,4.
\end{aligned}
\end{equation}
Thus, for $0<t<1$, using the fact
$V^\varepsilon=K^\varepsilon$,
$P_V=P_K$ and $G_u=G_K$,
by \eqref{3.18} and \eqref{4.26}, we obtain
\begin{align} \label{4.29}
F_u=F_0^\varepsilon+\varepsilon n^{\rm app}e_d,
\qquad
\int_0^1\|\varepsilon n^{\rm app}e_d\|_4\,dt
\le C\varepsilon.
\end{align}
From this and \eqref{4.28}, we deduce
\begin{align} \label{4.30}
\begin{aligned}
\|F_u\|_{L^1(0,1;L^4)}
\le \|F_0^\varepsilon\|_{L^1(0,1;L^4)}
+\varepsilon\int_0^1\|n^{\rm app}(t)\|_4\,dt\le C(\varepsilon^{1/8}+\varepsilon)
\le C\varepsilon^{1/8}.
\end{aligned}
\end{align}

Let us deal with the long-time fluid residuals.
For $1<t<S_\varepsilon$, we have
$V^\varepsilon=\varepsilon\bar u^\varepsilon+L^\varepsilon$.
The reference equation reads
\[
 (\varepsilon\bar u^\varepsilon)_t
 +(\varepsilon\bar u^\varepsilon\cdot\nabla)(\varepsilon\bar u^\varepsilon)
 -\delta\Delta(\varepsilon\bar u^\varepsilon)
 +\nabla(\varepsilon^2\bar p^\varepsilon)
 =\varepsilon^2(\bar n^\varepsilon e_d+\bar f_u^\varepsilon).
\]
Since $n^{\rm app}=\bar n^\varepsilon$, one has
\begin{equation*}
 F_u=F_0^\varepsilon
       -(\bar u^\varepsilon\cdot\nabla)L^\varepsilon
       -(L^\varepsilon\cdot\nabla)\bar u^\varepsilon,
\end{equation*}
where
\begin{equation}\label{4.31}
\begin{aligned}
 (\bar u^\varepsilon\cdot\nabla)L^\varepsilon={}&
 \sqrt\delta(\bar u^\varepsilon\cdot\nabla_xv)^\delta
 +(\bar u^\varepsilon\cdot\nabla \tilde{d})(v_\zeta)^\delta +\delta(\bar u^\varepsilon\cdot\nabla_xw)^\delta
 +\sqrt\delta(\bar u^\varepsilon\cdot\nabla \tilde{d})(w_\zeta)^\delta +\delta(\bar u^\varepsilon\cdot\nabla)\nabla\theta^\delta.
\end{aligned}
\end{equation}
From \(\bar u^\varepsilon\cdot\nu=0\), we have
$|\bar u^\varepsilon\cdot\nabla \tilde{d}|\le C\tilde{d}$, which implies
\begin{align*}
 \|(\bar u^\varepsilon\cdot\nabla \tilde{d})(v_\zeta)^\delta\|_4
 \le C\sqrt\delta\|(\zeta v_\zeta)^\delta\|_4
 \le C\delta^{5/8}m(t), ~~~~
 \|\sqrt\delta(\bar u^\varepsilon\cdot\nabla \tilde{d})(w_\zeta)^\delta\|_4
 \le C\delta\|(\zeta w_\zeta)^\delta\|_4
 \le C\delta^{9/8}m(t).
\end{align*}
The other terms in \eqref{4.31} are bounded by
\eqref{4.14} and \eqref{4.17}. Moreover,
\begin{equation*}
 \|(L^\varepsilon\cdot\nabla)\bar u^\varepsilon\|_4
 \le\|L^\varepsilon\|_4\|\nabla\bar u^\varepsilon\|_\infty
 \le C\delta^{5/8}m(t).
\end{equation*}
Therefore, for $1<t<S_\varepsilon$, we obtain
\[
 \|(\bar u^\varepsilon\cdot\nabla)L^\varepsilon\|_4
 +\|(L^\varepsilon\cdot\nabla)\bar u^\varepsilon\|_4
 \le C(\delta^{5/8}+\delta^{9/8}+\delta)m(t),
\]
which, together with \eqref{4.28}-\eqref{4.30}, yields
\begin{equation}\label{4.32}
\begin{aligned}
 \|F_u\|_{L^1(0,S_\varepsilon;L^4)}
 &\le C\varepsilon^{1/8}+C\varepsilon
 +C(\delta^{5/8}+\delta^{9/8}+\delta)\int_1^\infty m(t)\dd t\\
 &\le C\varepsilon^{1/8}.
\end{aligned}
\end{equation}

\paragraph{Uniform estimates and closure.}
The boundary estimate \eqref{4.25} already gives
\[
 \varepsilon\int_0^{S_\varepsilon}
          \|h^\varepsilon(t)\|_{4,\partial\OO}^2\dd t
 \le C\varepsilon^{1/2}.
\]
We next verify the coefficient bounds required by
Proposition~\ref{prop 3.2}. According to  \eqref{3.12}, \eqref{3.14}, \eqref{4.19}
and \eqref{4.20}, 
\begin{align*}
 \|n^{\rm app}\|_{W^{1,\infty}}+\|c^{\rm app}\|_{W^{1,\infty}}&\le C,\\
 \|\nabla V^\varepsilon\|_\infty+\|\Div V^\varepsilon\|_\infty
 +\|\alpha_n^\varepsilon\|_\infty
 +\|\alpha_c^\varepsilon\|_{W^{1,\infty}}
 &\le C\big(\mathbf1_{(0,1)}(t)+m(t)+\varepsilon\big).
\end{align*}
Since $\varepsilon S_\varepsilon=T_{\rm app}$, it follows that
\begin{equation*}
\begin{aligned}
 &\int_0^{S_\varepsilon}
 \big(\|\nabla V^\varepsilon\|_\infty+\|\Div V^\varepsilon\|_\infty
 +\|\alpha_n^\varepsilon\|_\infty
 +\|\alpha_c^\varepsilon\|_{W^{1,\infty}}\big)\dd t\\
 &\quad\le C\left(1+\int_0^\infty m(t)\dd t
                       +\varepsilon S_\varepsilon\right)
 \le C(1+T_{\rm app}).
\end{aligned}
\end{equation*}
Thus \eqref{3.23} holds. Substituting
\eqref{4.27}, \eqref{4.32} and \eqref{4.25}
into \eqref{3.37}, we obtain
\begin{equation}\label{4.33}
\begin{aligned}
 d_\varepsilon={}&
 \left(\int_0^{S_\varepsilon}(\|R_n\|_4+\|F_u\|_4)\dd t\right)^2
 +\varepsilon^{-1}\int_0^{S_\varepsilon}\|R_c\|_4^2\dd t +\varepsilon\int_0^{S_\varepsilon}
                     \|h^\varepsilon\|_{4,\partial\OO}^2\dd t\\
 \le{}&C\left[(\varepsilon^{5/8}+\varepsilon^{1/8})^2
            +\varepsilon^{-1}\varepsilon^{5/4}
            +\varepsilon^{1/2}\right]
 \le C\varepsilon^{1/4}.
\end{aligned}
\end{equation}
For sufficiently small $\varepsilon$, the estimates \(\eqref{3.23}\) and \(\eqref{4.33}\), together with the previously verified compatibility conditions, fulfill the hypotheses of Proposition~\(\ref{prop 3.2}\).
This, together with Lemma~\ref{lemma 3.3}, yields a strong solution on
$[0,S_\varepsilon]$ satisfying
\begin{equation*}
 \sup_{0\le t\le S_\varepsilon}
 \big(\|r(t)\|_4+\|q(t)\|_{W^{1,4}}+\|y(t)\|_4\big)
 \le C d_\varepsilon^{1/2}\le C\varepsilon^{1/8}.
\end{equation*}

At $t=S_\varepsilon=T_{\rm app}/\varepsilon$, the scalar profiles equal the
reference. With the help of \eqref{4.19}, one has
\[
 \varepsilon^{-1}\|L^\varepsilon(S_\varepsilon)\|_4
 \le C\varepsilon^{-3/8}(1+T_{\rm app}/\varepsilon)^{-2}
 \le C_{T_{\rm app}}\varepsilon^{13/8}.
\]
Transforming back to the original physical variables, we have
\[
\begin{aligned}
 n(s_0+T_{\rm app})-\bar n(s_0+T_{\rm app})&=r(S_\varepsilon),\\
 c(s_0+T_{\rm app})-\bar c(s_0+T_{\rm app})&=q(S_\varepsilon),\\
 u(s_0+T_{\rm app})-\bar u(s_0+T_{\rm app})
 &=y(S_\varepsilon)+\varepsilon^{-1}L^\varepsilon(S_\varepsilon).
\end{aligned}
\]
Therefore, we finally obtain
\begin{equation}\label{4.34}
 \|n-\bar n\|_4+\|c-\bar c\|_{W^{1,4}}+\|u-\bar u\|_4
 \le C\varepsilon^{1/8}\qquad\hbox{at }s=s_0+T_{\rm app}.
\end{equation}
The multiplicative controls \eqref{3.15} and
Lemma~\ref{lemma 2.2} ensure $n,c\ge0$.

\section{Local exact controllability to a positive trajectory}\label{section 5}
This section is devoted to proving the local exact controllability to the trajectories. 
The main difficulty is to handle the second order chemotactic
coupling together with the Navier slip-with-friction boundary
condition in the joint observability estimate.
We then obtain null controllability of the linearized system
by duality and prove the nonlinear result by a fixed point
argument.

Fix $[s_0,s_0+\tau]\subset[s_*,T]$, with $s_*>0$ as in Section~\ref{section 2}. All reference coefficients are evaluated at \((s_0+t,x)\), and their regularity bounds and positive lower bounds are independent of \(s_0\) for \(s_*\le s_0\le T-\tau\). Consider the state error system for $Z=(r,q,v)$
\begin{equation}\label{5.1}
\left\{\begin{aligned}
 &\begin{aligned}[b] &r_t-D_n\Delta r+\ub\cdot\nabla r+v\cdot\nabla\nb +\chi\Div(\nb\nabla q+r\nabla\cb)-\an r=F_n+h_n,\end{aligned} &&\text{in }Q_\tau,\\[2pt]
 &\begin{aligned}[b] &q_t-D_c\Delta q+\ub\cdot\nabla q+v\cdot\nabla\cb +(\nb-\ac)q+\cb r=F_c+h_c,\end{aligned} &&\text{in }Q_\tau,\\[2pt]
 &\begin{aligned}[b] &v_t-\mu\Delta v+(\ub\cdot\nabla)v+(v\cdot\nabla)\ub +\nabla\pi-r e_d=F_u+h_u,\end{aligned} ~~\Div v=0, &&\text{in }Q_\tau,\\[2pt]
 &\partial_\nu r=\partial_\nu q=0,\quad v\cdot\nu=0,\quad\NN v=0, &&\text{on }\Sigma_\tau,\\[2pt]
 &(r,q,v)|_{t=0}=(r_0,q_0,v_0), &&\text{in }\OO,
\end{aligned}\right.
\end{equation}
where the control $h=(h_n,h_c,h_u)$ is supported in $\omega\subset\subset\mathcal O\setminus\overline\Omega$. Let $F=(F_n,F_c,F_u)$ for a prescribed source and $Z_0=(r_0,q_0,v_0)$ for the initial error.

 Let $\phi,\psi,z$ and $\varpi$ be the adjoint states and pressure, with terminal data $(\phi_\tau,\psi_\tau,z_\tau)\in L^2\times L^2\times\Hh$ and $z_d=z\cdot e_d$. The homogeneous adjoint system is
\begin{equation}\label{5.2}
\left\{\begin{aligned}
 &\begin{aligned}[b] &-\phi_t-D_n\Delta\phi-\ub\cdot\nabla\phi -\chi\nabla\cb\cdot\nabla\phi-\an\phi+\cb\psi-z_d=0,\end{aligned} &&\text{in }Q_\tau,\\[2pt]
 &\begin{aligned}[b] &-\psi_t-D_c\Delta\psi-\ub\cdot\nabla\psi +(\nb-\ac)\psi+\chi\Div(\nb\nabla\phi)=0,\end{aligned} &&\text{in }Q_\tau,\\[2pt]
 &\begin{aligned}[b] &-z_t-\mu\Delta z-(\ub\cdot\nabla)z+(\nabla\ub)^{\mathsf T}z+\nabla\varpi+\phi\nabla\nb+\psi\nabla\cb=0,\end{aligned}~~~~\Div z=0, &&\text{in }Q_\tau,\\[2pt]
 &\partial_\nu\phi=\partial_\nu\psi=0,\quad z\cdot\nu=0,\quad\NN z=0, &&\text{on }\Sigma_\tau,\\[2pt]
 &(\phi,\psi,z)|_{t=\tau}=(\phi_\tau,\psi_\tau,z_\tau), &&\text{in }\OO.
\end{aligned}\right.
\end{equation}
According to the duality theory, the null controllability of \eqref{5.1} is equivalent to the observability of the adjoint system \eqref{5.2}.
\subsection{A Carleman inequality and observability}\label{subsection 5.1}
Choose connected smooth sets
$\omega_*\subset\subset\omega_2\subset\subset\omega_1\subset\subset\omega_0\subset\subset\omega$
and $\eta\in C^\infty(\overline\OO)$ satisfying
\[
 \eta>0\text{ in }\OO,\qquad \eta=0\text{ on }\partial\OO,
 \qquad |\nabla\eta|>0\text{ on }\overline\OO\setminus\omega_*.
\]
Let $\lambda,s>0$ be Carleman parameters and $\eta_{\max}=\max_{\overline\OO}\eta$. Define
\begin{align*}
 l(t)=t(\tau-t),~~~~
 \alpha=\frac{e^{2\lambda\eta_{\max}}-e^{\lambda\eta}}{l^4},~~~~
 \xi=\frac{e^{\lambda\eta}}{l^4},\\
 \alpha_*=\frac{e^{2\lambda\eta_{\max}}-1}{l^4},~~~~
 \widehat\alpha=\frac{e^{2\lambda\eta_{\max}}-e^{\lambda\eta_{\max}}}{l^4},~~~~
 b_{\rm C}=\frac{s}{l^4}.
\end{align*}
Let $w=e^{-2s\alpha}$, $w_*=e^{-2s\alpha_*}$ and $W=s\xi$. Then, one has
\begin{equation}\label{5.3}
 b_{\rm C}\le W\le e^{\lambda\eta_{\max}}b_{\rm C},\qquad w_*\le w\le e^{-2s\widehat\alpha}.
\end{equation}
For a scalar or vector field $f$, set
\begin{align*}
 I(f)&=\iint_{Q_\tau}w\big[W^{-1}(|f_t|^2+|\Delta f|^2)
             +\lambda^2W|\nabla f|^2+\lambda^4W^3|f|^2\big],\\
 I_1(f)&=\iint_{Q_\tau}w b_{\rm C}\big[W^{-1}(|f_t|^2+|\Delta f|^2)
             +\lambda^2W|\nabla f|^2+\lambda^4W^3|f|^2\big].
\end{align*}

The following classical estimates are available from \cite{CFKMS2023}, \cite{G2006} and \cite{IP2002}.
\begin{lemma}\label{lemma 5.1}
For sufficiently large $\lambda$ and $s\ge C_\lambda(\tau^4+\tau^8)$,
\begin{align}
 I(f)&\le C\iint_{Q_\tau}w|-f_t-\kappa\Delta f|^2
       +C\iint_{(0,\tau)\times\omega_2}w\lambda^4W^3|f|^2,
       \qquad\partial_\nu f=0,\label{5.4}\\
 I(z)&\le C\lambda\iint_{Q_\tau}w|-z_t-\mu\Delta z|^2
       +C\iint_{(0,\tau)\times\omega_2}w\lambda^4W^3|z|^2,
       \label{5.5}
\end{align}
where $\kappa>0$ is fixed. The second estimate
holds for $\Div z=0$, $z\cdot\nu=0$ and $\NN z=0$. If $\Delta p=\Div G$, then
\begin{align}
 \iint_{Q_\tau}w|\nabla p|^2
 \le{}&C\iint_{Q_\tau}wW|G|^2
 +C\int_0^\tau w_*b_{\rm C}^{1/2}\|p\|_{H^{1/2}(\partial\OO)}^2\dd t\nonumber\\
 &+C\iint_{(0,\tau)\times\omega_2}w(\lambda^2W^2|p|^2+|\nabla p|^2).
 \label{5.6}
\end{align}
The constants multiplying the volume terms in \eqref{5.5} and \eqref{5.6} are independent of sufficiently large $\lambda$.
\end{lemma}

We also recall the regularity estimates required for the pressure, which follow from \cite{G2006}. For a smooth bounded vector coefficient $b$, matrix coefficient $M$, and vector source $G$, the velocity $Z$ and pressure $P$ solving
\begin{equation*}
\left\{\begin{aligned}
 &-Z_t-\mu\Delta Z+b\cdot\nabla Z+MZ+\nabla P=G, &&\text{in }Q_\tau,\\[2pt]
 &\Div Z=0, &&\text{in }Q_\tau,\\[2pt]
 &Z\cdot\nu=0,\quad\NN Z=0, &&\text{on }\Sigma_\tau,\\[2pt]
 &Z(\tau,\cdot)=0, &&\text{in }\OO
\end{aligned}\right.
\end{equation*}
satisfy
\begin{equation}\label{5.7}
 \|Z_t\|_{L^2L^2}+\|Z\|_{L^2H^2}+\|P\|_{L^2H^1}
 \le C_T\|G\|_{L^2L^2}.
\end{equation}
The pressure is normalized by zero mean, either on $\OO$ or on $\omega_2$. By extending the coefficients and the forcing to the fixed interval \((0,T)\), with uniform bounds in the required norms, the constant can be chosen independently of \(0<\tau\le T\).

To control the pressure and higher order derivatives of the fluid component in Carleman estimates,
we need the following elliptic regularity result for the stationary Stokes system with Navier slip-with-friction boundary conditions.
\begin{lemma}\label{lemma 5.2}
For the symmetric matrix $\beta$, the stationary Navier problem
\begin{equation*}
\left\{\begin{aligned}
 &-\mu\Delta Z+\nabla P=G,\quad\Div Z=0, &&\text{in }\OO,\\[2pt]
 &Z\cdot\nu=0,\quad[2\mu D(Z)\nu+\beta Z]_{\tr}=h, &&\text{on }\partial\OO
\end{aligned}\right.
\end{equation*}
with $\int_\OO P\,dx=0$ and $h\cdot\nu=0$, has the regularity estimate, for $j=0,2$,
\begin{equation}\label{5.8}
 \|Z\|_{H^{j+2}}+\|P\|_{H^{j+1}}
 \le C\big(\|G\|_{H^j}+\|h\|_{H^{j+1/2}(\partial\OO)}+\|Z\|_2\big).
\end{equation}
\end{lemma}
\begin{proof}
Dividing the momentum equation by $\mu$ and adding the
zero-order term $Z$ gives the equivalent system
\begin{equation*}
\left\{\begin{aligned}
 &Z-\Delta Z+\nabla(P/\mu)=Z+G/\mu,\quad\Div Z=0,
 &&\text{in }\OO,\\
 &Z\cdot\nu=0,\quad
 [2D(Z)\nu]_{\tr}
 =\mu^{-1}(h-[\beta Z]_{\tr}),
 &&\text{on }\partial\OO.
\end{aligned}\right.
\end{equation*}
Applying the tangential-stress regularity theorem
\cite{Saito2004} in boundary charts yields
the corresponding local estimates. A finite covering of
$\partial\OO$, together with interior regularity, gives
the global estimate
\[
\begin{aligned}
 \|Z\|_{H^{k+2}}+\|P\|_{H^{k+1}}
 \le C\big(&\|G\|_{H^k}
 +\|h\|_{H^{k+1/2}(\partial\OO)} +\|Z\|_{H^{k+1}}+\|P\|_{H^k}\big),
 \qquad k=0,1,2.
\end{aligned}
\]
Indeed, localization with a cutoff $\zeta_{\rm loc}$ satisfying
$\partial_\nu\zeta_{\rm loc}=0$ produces
$\Div(\zeta_{\rm loc} Z)=\nabla\zeta_{\rm loc}\cdot Z$ and commutator terms
involving only lower-order derivatives of $Z$.
These terms are controlled by the lower-order norms
appearing above.

The standard pressure estimate, with $P$ normalized
to have zero mean, controls its lower order norm.
The interpolation inequality
\[
 \|Z\|_{H^{k+1}}
 \le \delta\|Z\|_{H^{k+2}}+C_\delta\|Z\|_2
\]
and the corresponding interpolation estimate for $P$
allow the lower-order terms to be absorbed. This proves
\eqref{5.8}. Changing the pressure
normalization adds only a spatial constant, which is
controlled by the normalization condition.
\end{proof}

Combining the Carleman estimates and the pressure
regularity estimate, the observability
inequality follows.

\begin{proposition}\label{prop 5.3}
There exist fixed $\lambda,K,S>0$ such that, for $s=S(\tau^4+\tau^8)$, every solution of \eqref{5.2} satisfies
\begin{equation}\label{5.9}
 K I_1(\phi)+I(\psi)+I(z)
 \le C\iint_{(0,\tau)\times\omega_0}\Theta(t)(|\phi|^2+|\psi|^2+|z|^2),
\end{equation}
where
$
 \Theta=e^{4\lambda\eta_{\max}}b_{\rm C}^{10}e^{-4s\widehat\alpha+2s\alpha_*}.
$
Furthermore,
\begin{equation}\label{5.10}
 \|\phi(0)\|_2^2+\|\psi(0)\|_2^2+\|z(0)\|_2^2
 \le C e^{C(1+\tau^{-4})}
 \iint_{(0,\tau)\times\omega_0}(|\phi|^2+|\psi|^2+|z|^2).
\end{equation}
Constants $C$ are uniform for $s_0\ge s_*$ and $s_0+\tau\le T$.
\end{proposition}
\begin{proof}
Let $R=\min b_{\rm C}=256s/\tau^8\ge256S$. Multiplying the first adjoint equation in
\eqref{5.2} by $b_{\rm C}^{1/2}$, one has
\[
\begin{aligned}
(-\partial_t-\kappa\Delta)
(b_{\rm C}^{1/2}\phi)
&=b_{\rm C}^{1/2}(-\phi_t-\kappa\Delta\phi)
-(b_{\rm C}^{1/2})'\phi.
\end{aligned}
\]
Using \eqref{5.4}, we obtain
\begin{equation}\label{5.11}
 I_1(f)\le C_\lambda\iint w b_{\rm C}|-f_t-\kappa\Delta f|^2
 +C_\lambda\iint_{(0,\tau)\times\omega_2}w b_{\rm C}W^3|f|^2,
\end{equation}
where we used $|(b_{\rm C}^{1/2})'|^2\le4\tau^2s^{-1/2}b_{\rm C}^{3/2}$ and
$\tau^2s^{-1/2}\le S^{-1/2}$.

The additional term arising from the time derivative
of $b_{\rm C}^{1/2}$ is bounded by
\[
C_\lambda S^{-1/2}R^{-5/2}I_1(f),
\]
and is absorbed into the left-hand side for sufficiently
large $S$. Applying \eqref{5.11} to the first
equation of \eqref{5.2} and \eqref{5.4}
to the second equation, we estimate the scalar components.
Since $\bar n$ and $\nabla\bar n$ are uniformly bounded,
$
|\Div(\bar n\nabla\phi)|^2
\le C\bigl(|\Delta\phi|^2+|\nabla\phi|^2\bigr).
$
Moreover, it follows from \eqref{5.3} that
\[
\iint_{Q_\tau}w|\Delta\phi|^2
\le e^{\lambda\eta_{\max}}I_1(\phi).
\]
Thus the second-order chemotactic coupling contributes
at most $C_\lambda I_1(\phi)$ to the estimate for
$I(\psi)$. Choosing $K$ sufficiently large allows this
term to be absorbed into $K I_1(\phi)$.

Next, applying \eqref{5.5} to the third
equation of \eqref{5.2} and absorbing the
remaining lower-order terms involving the fluid adjoint
variable by increasing $R$, we obtain
\begin{equation}\label{5.12}
\mathcal J
\le C_\lambda\mathcal I_{\rm obs}
+C_0\lambda\iint_{Q_\tau}w|\nabla\varpi|^2,
\end{equation}
where
\[
\mathcal J=K I_1(\phi)+I(\psi)+I(z),
\qquad
\mathcal I_{\rm obs}=
\iint_{(0,\tau)\times\omega_2}
w\bigl(b_{\rm C}W^3|\phi|^2
+W^3|\psi|^2+W^3|z|^2\bigr).
\]
The constant $C_0$ is independent of $K$ and 
sufficiently large $\lambda$, since the fluid Carleman
estimate is not multiplied by $K$.
To estimate the pressure term in \eqref{5.12},
let us set
\[
b=-\bar u,\qquad
M=(\nabla\bar u)^{\mathsf T},\qquad
F=-\phi\nabla\bar n-\psi\nabla\bar c.
\]
The fluid adjoint equation and its divergence yield
\begin{equation}\label{5.13}
\left\{\begin{aligned}
 &-z_t-\mu\Delta z+b\cdot\nabla z+Mz+\nabla\varpi=F,
 &&\text{in }Q_\tau,\\
 &\Delta\varpi=\Div(F-b\cdot\nabla z-Mz),
 &&\text{in }Q_\tau.
\end{aligned}\right.
\end{equation}
Applying the pressure estimate \eqref{5.6},
the first volume term is absorbed because
\[
C\lambda\iint_{Q_\tau}wW|\nabla z|^2
\le \frac{C}{\lambda}I(z).
\]
Moreover, since $W\geq R$, comparison with the
zero order terms in $\mathcal J$ gives
\[
\iint_{Q_\tau}wW
\bigl(|\phi|^2+|\psi|^2+|z|^2\bigr)
\le C_\lambda R^{-2}\mathcal J.
\]
We first fix $\lambda$ sufficiently large, then choose
$K$ to absorb the second-order chemotactic coupling,
and finally increase $S$ and $R$ as needed to absorb
the lower-order terms.

Normalize the pressure by
\[
\int_{\omega_2}\varpi(t,x)\,\dd x=0
\quad\text{for a.e. }t\in(0,\tau),
\]
and take
$
\rho_\Gamma=b_{\rm C}^{1/4}e^{-s\alpha_*}.
$
Applying \eqref{5.7}
to the first equation of \eqref{5.13}, after
multiplication by $\rho_\Gamma$, we have
\begin{align}
 \int_0^\tau w_*b_{\rm C}^{1/2}
       \|\varpi\|_{H^{1/2}(\partial\OO)}^2\,\dd t
 \le C_T\left(
       \|\rho_\Gamma F\|_{L^2(0,\tau;L^2)}^2
       +\|\rho_\Gamma' z\|_{L^2(0,\tau;L^2)}^2
       \right)\le C_\lambda
       (R^{-5/2}+S^{-1/2})\mathcal J,
       \label{5.14}
\end{align}
where
$
|\rho_\Gamma'|^2
\le C_\lambda\tau^2s^{-1/2}w_*b_{\rm C}^3.
$
Choosing $R$ and $S$ sufficiently large makes the
right-hand side of \eqref{5.14} small enough
to absorb the pressure boundary term in
\eqref{5.6}. Thus, we have
\begin{equation}\label{5.15}
\mathcal J\le C_\lambda(\mathcal I_{\rm obs}+\mathcal P_{\rm loc}),
\qquad
\mathcal P_{\rm loc}=
\iint_{(0,\tau)\times\omega_2}
w\bigl(\lambda^2W^2|\varpi|^2+|\nabla\varpi|^2\bigr).
\end{equation}

For the localized pressure term
$\mathcal P_{\rm loc}$, by the pressure normalization, Poincar\'e's
inequality on $\omega_2$, and the first equation of
\eqref{5.13}, we obtain
\begin{equation}\label{5.16}
\mathcal P_{\rm loc}\le C_\lambda
\iint_{(0,\tau)\times\omega_2}
P_0\bigl(
|z_t|^2+|\Delta z|^2+|\nabla z|^2
+|z|^2+|\phi|^2+|\psi|^2
\bigr),
\end{equation}
where $
P_0=e^{2\lambda\eta_{\max}}b_{\rm C}^2e^{-2s\widehat\alpha}.
$
To eliminate these local derivatives, set $\rho=b_{\rm C}^{-3}e^{-s\alpha_*}$ and
$d_S=R^{-5}+S^{-1/2}R^{-5/2}$. The weight comparisons give
\[
 \rho^2\le C_\lambda R^{-5}wW^{-1},\qquad
 |\rho'|^2\le C_\lambda S^{-1/2}R^{-5/2}wW^{-1}.
\]
Using the scalar Neumann estimate $\|f\|_{H^2}\le C(\|\Delta f\|_2+\|f\|_2)$ with $f=\phi, \psi$, we infer
\[
 \begin{aligned}
 &\|\rho\phi_t\|_{L^2L^2}^2+\|\rho\psi_t\|_{L^2L^2}^2
 +\|\rho\phi\|_{L^2H^2}^2+\|\rho\psi\|_{L^2H^2}^2 +\|\rho'z_t\|_{L^2L^2}^2+\|\rho z\|_{L^2H^1}^2
 \le C_\lambda d_S\mathcal J.
 \end{aligned}
\]
Differentiating \eqrow{5.13}{1} 
with respect to time, we find that the field
$Y=\rho z_t$ satisfies
\begin{equation*}
\left\{\begin{aligned}
 &\begin{aligned}[b] &-Y_t-\mu\Delta Y+b\cdot\nabla Y+MY+\nabla(\rho\varpi_t)=\rho F_t-\rho^\prime z_t-\rho b_t\cdot\nabla z-\rho M_tz,\end{aligned} &&\text{in }Q_\tau,\\[2pt]
 &\Div Y=0, &&\text{in }Q_\tau,\\[2pt]
 &Y\cdot\nu=0,\quad\NN Y=0, &&\text{on }\Sigma_\tau.
\end{aligned}\right.
\end{equation*}
Since $F_t$ contains only $\phi_t,\psi_t,\phi,\psi$ with bounded coefficients, \eqref{5.7} gives
\[
 \|\rho z_{tt}\|_{L^2L^2}^2+\|\rho z_t\|_{L^2H^2}^2
 \le C_\lambda d_S\mathcal J.
\]
Applying \eqref{5.8} with $j=2$ to \eqrow{5.13}{1} and interpolating $H^3$ between $L^2$ and $H^4$, we also have
\[
 \|\rho z\|_{L^2H^4}
 \le C\big(\|\rho z_t\|_{L^2H^2}+\|\rho\phi\|_{L^2H^2}
            +\|\rho\psi\|_{L^2H^2}+\|\rho z\|_{L^2L^2}\big).
\]
Therefore, we obtain
\begin{equation*}
 \|\rho z_{tt}\|_{L^2L^2}^2+\|\rho z_t\|_{L^2H^2}^2
 +\|\rho z\|_{L^2H^4}^2\le C_\lambda d_S\mathcal J.
\end{equation*}

By direct differentiation, one has
\[
 P_0+|P_0''|\le C_\lambda\Theta,\qquad \frac{P_0^2}{\rho^2}=\Theta,
 \qquad 2s\widehat\alpha-s\alpha_*=(e^{\lambda\eta_{\max}}-1)^2b_{\rm C}>0.
\]
Choosing $\zeta\in C_c^\infty(\omega_1)$ with $0\le\zeta\le1$ and $\zeta=1$ on $\omega_2$, and integrating by parts in time and twice in space, respectively, we have
\begin{align*}
 \iint_{(0,\tau)\times\omega_1}\zeta P_0|z_t|^2
 &=\frac12\iint_{(0,\tau)\times\omega_1}
       \zeta P_0''|z|^2
   -\iint_{(0,\tau)\times\omega_1}
       \zeta P_0 z\cdot z_{tt}\\
 &\le C\iint_{(0,\tau)\times\omega_1}
       \Theta|z|^2
   +C\|\rho z_{tt}\|_{L^2(0,\tau;L^2(\mathcal O))}^2,\\
 \iint_{(0,\tau)\times\omega_1}\zeta P_0|\Delta z|^2
 &=\iint_{(0,\tau)\times\omega_1}
       P_0z\cdot\Delta(\zeta\Delta z)\\
 &\le C\iint_{(0,\tau)\times\omega_1}
       \Theta|z|^2
   +C\|\rho z\|_{L^2(0,\tau;H^4(\mathcal O))}^2,\\
 \iint_{(0,\tau)\times\omega_1}\zeta P_0|\nabla z|^2
 &\le C\iint_{(0,\tau)\times\omega_1}
       \Theta|z|^2
   +C_\lambda d_S\mathcal J.
\end{align*}
Thus, we get
\[
 \mathcal P_{\rm loc}\le C_\lambda d_S\mathcal J
 +C_\lambda\iint_{(0,\tau)\times\omega_1}\Theta(|\phi|^2+|\psi|^2+|z|^2).
\]
Combining the bounds for the localized pressure and observation
terms with \eqref{5.15}, we obtain
\[
 \mathcal J
 \le C_\lambda d_S\mathcal J
 +C_\lambda
 \iint_{(0,\tau)\times\omega_0}
 \Theta(t)\bigl(|\phi|^2+|\psi|^2+|z|^2\bigr)\,dx\,dt.
\]
Since $R\ge256S$ and
$d_S=R^{-5}+S^{-1/2}R^{-5/2}$,
we choose $S$ sufficiently large so that
$C_\lambda d_S\le1/2$.
Absorbing the first term on the right-hand side into the
left-hand side yields \eqref{5.9}.
The calculations are first justified for smooth terminal
data by means of endpoint cutoffs. The cutoff errors
vanish because of the rapid decay of the Carleman weights
at $t=0$ and $t=\tau$. A density argument extends
the estimate to general $L^2$ terminal data.

It remains to recover the initial norm. Testing \eqref{5.2} by $K_{\rm adj}\phi,\psi,z$, choosing $K_{\rm adj} D_n>4\chi^2\|\nb\|_\infty^2/D_c$, and using \eqref{2.4}, we have
\begin{equation}\label{5.17}
 -\frac{\dd}{\dd t}(K_{\rm adj}\|\phi\|_2^2+\|\psi\|_2^2+\|z\|_2^2)
 \le C(K_{\rm adj}\|\phi\|_2^2+\|\psi\|_2^2+\|z\|_2^2),
\end{equation}
where we have used
\[
 \chi\left|\int\nb\nabla\phi\cdot\nabla\psi\right|
 \le\tfrac{D_c}{4}\|\nabla\psi\|_2^2
       +\frac{\chi^2\|\nb\|_\infty^2}{D_c}\|\nabla\phi\|_2^2.
\]
Combining 
\eqref{5.17} with integration over
$[\tau/3,2\tau/3]$ and the definitions of $I$ and
$I_1$, we obtain
\[
\|\phi(0)\|_2^2+\|\psi(0)\|_2^2+\|z(0)\|_2^2
\le
C\tau^{-1}(\tau^8/s)^3e^{Cs\tau^{-8}}\mathcal J
\le Ce^{C(1+\tau^{-4})}\mathcal J.
\]

Finally, $\Theta=e^{4\lambda\eta_{\max}}b_{\rm C}^{10}e^{-2(e^{\lambda\eta_{\max}}-1)^2b_{\rm C}}\le C_\lambda$. This proves \eqref{5.10} and the uniformity in $s_0$.
\end{proof}

\subsection{Null controllability of the linearized system}
We recall the spaces $\Y,\F,\Uc$ and $\X$ defined in \eqref{1.3}-\eqref{1.5}.
\begin{lemma}\label{lemma 5.4}
For $Z_0\in\Y$ and $F=0$, system
\eqref{5.1} is null controllable in time
$\tau$. That is, there exists a linear choice of controls
$h$, supported in $\omega$, such that the corresponding
solution satisfies
\[
Z(\tau)=0
\]
and
\begin{equation}\label{5.18}
 \|Z\|_{\X(0,\tau)}
 +\|h\|_{L^2(0,\tau;\Uc)}
 \le Ce^{C(1+\tau^{-4})}\|Z_0\|_\Y.
\end{equation}
Moreover, for zero initial data and zero control, the solution
with source $F$ satisfies
\begin{equation}\label{5.19}
 \|Z\|_{\X(0,\tau)}+\|Z(\tau)\|_\Y
 \le C_L\|F\|_{L^2(0,\tau;\F)},
 \qquad 0<\tau\le L.
\end{equation}
\end{lemma}

\begin{proof}
By \eqref{5.10}, the functional $-\langle Z_0,(\phi,\psi,z)(0)\rangle$ is bounded. Its Riesz representative defines a linear control $h\in L^2L^2$, supported in $\omega_0$, with
\[
 \begin{aligned}
 \|h\|_{L^2L^2}&\le Ce^{C(1+\tau^{-4})}\|Z_0\|_\Y,\\
 \langle Z(\tau),(\phi,\psi,z)(\tau)\rangle
 &=\langle Z_0,(\phi,\psi,z)(0)\rangle
   +\iint_{(0,\tau)\times\omega_0}h\cdot(\phi,\psi,z)=0.
 \end{aligned}
\]
Thus $Z(\tau)=0$. 

Testing \eqref{5.1} by $r,K(q-\Delta q),v$, taking $K$
sufficiently large and using \eqref{2.4}, we first obtain
\[
 \|r,v\|_{L^\infty L^2\cap L^2H^1}
 +\|q\|_{L^\infty H^1\cap L^2H^2}
 \le C_\tau(\|Z_0\|_\Y+\|h\|_{L^2L^2}).
\]
All reference coefficients are smooth and bounded. Neumann maximal
regularity in the chemical equation gives $q_t\in L^2L^2$.
Then $\Div(\nb\nabla q)$ and the other terms in the density equation
belong to $L^2L^2$, so that Neumann maximal regularity gives
$r\in H^1L^2\cap L^2H^2$. Stokes maximal regularity gives
$v\in\X_u$. Therefore, the trace theorem yields  
\begin{equation*}
 \|r\|_{H^1L^2\cap L^2H^2\cap CH^1}
 +\|q\|_{H^1L^2\cap L^2H^2\cap CH^1}
 +\|v\|_{\X_u}
 \le C_\tau(\|Z_0\|_\Y+\|h\|_{L^2L^2}).
\end{equation*}
We now improve the scalar control regularity by localization. With this fixed velocity $v$, let $(r^f,q^f)$ solve the equations \eqrow{5.1}{1,2} with zero source and zero scalar control, and the same initial data. Choose $\zeta\in C_c^\infty(\omega)$ equal to one near $\overline{\omega_0}$ and a time cutoff $\gamma$ equal to one near $0$ and zero near $\tau$. Let
\begin{equation*}
 \widetilde r=(1-\zeta)r+\zeta\gamma r^f,\qquad
 \widetilde q=(1-\zeta)q+\zeta\gamma q^f,\qquad \widetilde v=v.
\end{equation*}
The new scalar controls are
\begin{align*}
 \widetilde h_n={}&D_n\big(2\nabla\zeta\cdot\nabla(r-\gamma r^f)
                 +(\Delta\zeta)(r-\gamma r^f)\big)
 -(\ub+\chi\nabla\cb)\cdot\nabla\zeta\,(r-\gamma r^f)\\
 &-\chi\big(2\nb\nabla\zeta\cdot\nabla(q-\gamma q^f)
                +(q-\gamma q^f)\Div(\nb\nabla\zeta)\big)
 +\zeta\gamma' r^f+\zeta(1-\gamma)v\cdot\nabla\nb,\\
 \widetilde h_c={}&D_c\big(2\nabla\zeta\cdot\nabla(q-\gamma q^f)
                 +(\Delta\zeta)(q-\gamma q^f)\big)
 -\ub\cdot\nabla\zeta\,(q-\gamma q^f)
 +\zeta\gamma'q^f+\zeta(1-\gamma)v\cdot\nabla\cb, 
\end{align*}
which are supported in $\omega$. Since every term has at most one derivative of the scalar states, we have
\begin{equation*}
 \|\widetilde h_n,\widetilde h_c\|_{L^2H^1}
 \le C(1+\|\gamma'\|_\infty)
       (\|r,q,r^f,q^f\|_{L^2H^2}+\|v\|_{L^2H^1}).
\end{equation*}
The localized scalar controls satisfy
$\widetilde h_n,\widetilde h_c\in L^2(0,\tau;H^1)$,
and the fluid control is replaced by
\begin{equation*}
 \widetilde h_u=h_u+(r-\widetilde r)e_d,
 \qquad \widetilde r e_d+\widetilde h_u=r e_d+h_u.
\end{equation*}
By $r-\widetilde r=\zeta(r-\gamma r^f)$, this correction is supported
in $\omega$ and belongs to $L^2(0,\tau;L^2)^d$.
The initial and terminal states remain unchanged.
Since the chemical source belongs to
$L^2(0,\tau;H^1)$, by maximal regularity for the Neumann
equation, we first get $\widetilde q\in\mathcal X_c$.
The density equation then implies
$\widetilde r\in\mathcal X_n$, and Stokes maximal
regularity gives us $v\in\mathcal X_u$.
The localization procedure is linear. Moreover, the
derivatives of the spatial and temporal cutoff functions
introduce at most polynomial factors in $\tau^{-1}$.
Since the control cost in \eqref{5.18} grows exponentially
in $\tau^{-4}$, these factors can be absorbed into the
constant in \eqref{5.18}.

For \eqref{5.19}, it follows from the basic energy estimate
that $r\in\mathcal X_n$ and $v\in\mathcal X_u$.
Since $F_c\in L^2(0,\tau;H^1)$, the chemical equation
then yields $q\in\mathcal X_c$. Extending the source
to the fixed interval $(0,L)$ gives a constant
independent of $0<\tau\le L$ and the continuous
terminal trace in $\mathcal Y$.

\end{proof}

To treat nonlinear sources up to the terminal time, fix
$m=4$ and choose $\vartheta>1$ and
$\alpha_0,\beta_0>0$ such that
\begin{equation}\label{5.20}
\vartheta^8<2,\qquad
\alpha_0\vartheta^8<\beta_0<2\alpha_0,\qquad
\beta_0\vartheta^{-4}-\alpha_0\vartheta^4
>C_1(1-\vartheta^{-1})^{-4}+1,
\end{equation}
where $C_1$ exceeds the cost constant in
\eqref{5.18}. These inequalities can be satisfied
by choosing an admissible ratio $\beta_0/\alpha_0$
and then taking $\alpha_0,\beta_0$ sufficiently large.
For a fixed duration $L$, define
\[
\rho_0(t)=\exp\left(-\frac{\alpha_0}{(L-t)^4}\right),
\qquad
\rho_F(t)=\exp\left(-\frac{\beta_0}{(L-t)^4}\right),
\]
and
\begin{align*}
 \X_{\rho_0}&=\{Z:Z/\rho_0\in\X(0,L)\},&
 \F_{\rho_F}&=\{F:F/\rho_F\in L^2(0,L;\F)\},\\
 \Uc_{\rho_0}&=\{h:h/\rho_0\in L^2(0,L;\Uc),\ \supp h\subset[0,L]\times\overline\omega\},&
 \mathcal W&=\X_{\rho_0}\times\Uc_{\rho_0}.
\end{align*}
The weighted norms are $\|Z\|_{\X_{\rho_0}}=\|Z/\rho_0\|_\X$, $\|F\|_{\F_{\rho_F}}=\|F/\rho_F\|_{L^2\F}$ and $\|h\|_{\Uc_{\rho_0}}=\|h/\rho_0\|_{L^2\Uc}$.
\begin{proposition}\label{prop 5.5}
There exists a bounded linear operator that maps
each initial datum and source term to a state-control
pair for \eqref{5.1}, such that
\begin{equation}\label{5.21}
 \|Z\|_{\X_{\rho_0}}+\|h\|_{\Uc_{\rho_0}}
 \le C_L(\|Z_0\|_\Y+\|F\|_{\F_{\rho_F}}),\qquad Z(L)=0.
\end{equation}
\end{proposition}

\begin{proof}
Let
\[
t_k=L(1-\vartheta^{-k}),\qquad
I_k=[t_k,t_{k+1}],\qquad k\ge0.
\]
Then $t_k\to L$ as $k\to\infty$ and
$
|I_k|=(1-\vartheta^{-1})L\vartheta^{-k}.
$
We now construct the solution successively on these intervals.
First, use Lemma~\ref{lemma 5.4} to drive $Z_0$
to zero on $I_0$. For each $k\ge0$, solve the system
with zero initial data and source $F\mathbf1_{I_k}$
on $I_k$, and then use Lemma~\ref{lemma 5.4}
on $I_{k+1}$ to drive the resulting terminal state
to zero. Denote the resulting local state and control
by $(Z_k,h_k)$, extended by zero outside
$[t_k,t_{k+2}]$. The terminal states match at each
junction, and at most two local pairs are nonzero
on any interval $I_j$. Summing these pairs defines
the global state $Z$ and control $h$.

By \eqref{5.19}, the norm of the state to be
controlled on $I_j$ is bounded by
$
C_L\|F\|_{L^2(I_{j-1};\mathcal F)}.
$
Combining this bound with \eqref{5.18} and
the interval lengths above, the weighted control
cost is governed by the exponent
\[
\frac{\vartheta^{4j}}{L^4}
\left[
\alpha_0\vartheta^4
+C_1(1-\vartheta^{-1})^{-4}
-\beta_0\vartheta^{-4}
\right]
\le -\frac{\vartheta^{4j}}{L^4},
\]
by \eqref{5.20}. The directly forced terms
are also summable, since
$\beta_0>\alpha_0\vartheta^8>\alpha_0\vartheta^4$.
Moreover,
\[
\left|\frac{\rho_0'}{\rho_0}\right|
=4\alpha_0(L-t)^{-5},
\]
and the strict negative margin in the exponent
dominates this polynomial factor. Summing the squared
weighted estimates, we obtain
\[
\|Z\|_{\mathcal X_{\rho_0}}
+\|h\|_{\mathcal U_{\rho_0}}
\le C_L\left(
\|Z_0\|_{\mathcal Y}
+\|F\|_{\mathcal F_{\rho_F}}
\right).
\]

Finally, $Z/\rho_0\in C([0,L];\mathcal Y)$ and
$\rho_0(t)\to0$ as $t\to L^-$. Hence
\[
\lim_{t\to L^-}\|Z(t)\|_{\mathcal Y}=0.
\]
Defining $Z(L)=0$ gives the required continuous
terminal trace. Since every local control is chosen
linearly by Lemma~\ref{lemma 5.4}, the global
state-control map is linear and bounded. This proves
\eqref{5.21} and concludes the proof.
\end{proof}

\subsection{Local exact controllability of the nonlinear system}
For $(r,q,v,h_n,h_c,h_u)\in\mathcal W$, define
\begin{align}
 F_n&=-v\cdot\nabla r-\chi\nabla r\cdot\nabla q-\chi r\Delta q
                  +\nb^{-1}h_nr,\label{5.22}\\
 F_c&=-v\cdot\nabla q-rq+\cb^{-1}h_cq,\label{5.23}\\
 F_u&=-(v\cdot\nabla)v.\nonumber
\end{align}
These are the nonlinear sources corresponding to the actual controls
\begin{equation}\label{5.24}
 f_n=(\an+h_n/\nb)(\nb+r),\qquad
 f_c=(\ac+h_c/\cb)(\cb+q),\qquad f_u=\bar f_u+h_u.
\end{equation}
In particular, the control-state products are retained in \eqref{5.22}-\eqref{5.23}.

\begin{lemma} 
\label{lemma 5.6}
Let
$
Z=(r,q,v),\quad h=(h_n,h_c,h_u),\quad (Z,h)\in\mathcal W, 
R=\|(Z,h)\|_{\mathcal W}.
$
Then
\begin{equation}\label{5.25}
\begin{aligned}
\left\|\frac{F_n(Z,h)}{\rho_F}\right\|_{L^2(0,L;L^2)}
+\left\|\frac{F_c(Z,h)}{\rho_F}\right\|_{L^2(0,L;H^1)}
 +\left\|\frac{F_u(Z,h)}{\rho_F}\right\|_{L^2(0,L;L^2)}
\le CR^2.
\end{aligned}
\end{equation}
Moreover, for any two state-control pairs
$(Z_i,h_i)\in\mathcal W$,
$i=1,2$, with $R_i=\|(Z_i,h_i)\|_{\mathcal W}$, we have
\begin{equation}\label{5.26}
\begin{aligned}
&\left\|\frac{F_n(Z_1,h_1)-F_n(Z_2,h_2)}{\rho_F}
  \right\|_{L^2(0,L;L^2)}
\quad+
\left\|\frac{F_c(Z_1,h_1)-F_c(Z_2,h_2)}{\rho_F}
  \right\|_{L^2(0,L;H^1)}
\\ &\quad +
\left\|\frac{F_u(Z_1,h_1)-F_u(Z_2,h_2)}{\rho_F}
  \right\|_{L^2(0,L;L^2)}
\\
&\qquad\le
C(R_1+R_2)\|(Z_1-Z_2,h_1-h_2)\|_{\mathcal W}.
\end{aligned}
\end{equation}
\end{lemma}

\begin{proof}
In the nonlinear sources $F_n, F_c, F_u$, the pure state terms are estimated by
\eqref{3.44}. More precisely,
\[
\begin{aligned}
&\|v\cdot\nabla r\|_{L^2(0,L;L^2)}
+\|\nabla r\cdot\nabla q\|_{L^2(0,L;L^2)}
+\|r\Delta q\|_{L^2(0,L;L^2)}
 \le
C\bigl(\|v\|_{\mathcal X_u}\|r\|_{\mathcal X_n}
+\|r\|_{\mathcal X_n}\|q\|_{\mathcal X_c}\bigr),\\
&\|v\cdot\nabla q\|_{L^2(0,L;H^1)}
+\|rq\|_{L^2(0,L;H^1)}
 \le
C\bigl(\|v\|_{\mathcal X_u}\|q\|_{\mathcal X_c}
+\|r\|_{\mathcal X_n}\|q\|_{\mathcal X_c}\bigr),\\
&\|(v\cdot\nabla)v\|_{L^2(0,L;L^2)}
\le C\|v\|_{\mathcal X_u}^2.
\end{aligned}
\]
Here the constants may depend on the fixed interval length $L$.
Therefore, the six pure state terms are bounded by
$CR^2$.
For the control state products, thanks to Sobolev multiplication
and the definitions of $\mathcal X_n$ and $\mathcal X_c$, we have 
\[
\begin{aligned}
\|h_nr\|_{L^2(0,L;L^2)}
&\le C
\|h_n\|_{L^2(0,L;H^1)}
\|r\|_{L^\infty(0,L;H^1)},\\
\|h_cq\|_{L^2(0,L;H^1)}
&\le C
\|h_c\|_{L^2(0,L;H^1)}
\|q\|_{L^\infty(0,L;H^2)}.
\end{aligned}
\]
Since the reference trajectory is strictly positive and
has the required spatial regularity, multiplication by
$\bar n^{-1}$ and $\bar c^{-1}$ is bounded in these
spaces.

It remains to incorporate the time weights. By
\eqref{5.20}, $\beta_0<2\alpha_0$, and hence
\[
\frac{\rho_0(t)^2}{\rho_F(t)}
=\exp\left(
-\frac{2\alpha_0-\beta_0}{(L-t)^4}
\right)\le1,\qquad 0\le t<L.
\]
Combining this inequality with the product estimates
above, we obtain
\[
\begin{aligned}
&\left\|\frac{F_n(Z,h)}{\rho_F}\right\|_{L^2(0,L;L^2)}
+\left\|\frac{F_c(Z,h)}{\rho_F}\right\|_{L^2(0,L;H^1)}
 +
\left\|\frac{F_u(Z,h)}{\rho_F}\right\|_{L^2(0,L;L^2)}
\le C\|(Z,h)\|_{\mathcal W}^2
=CR^2,
\end{aligned}
\]
which implies \eqref{5.25}.

For two state control pairs $(Z_1,h_1),(Z_2,h_2)$,  
applying the same product and weight estimates to these
differences of the nonlinear source terms
$F_n$, $F_c$ and $F_u$ yields \eqref{5.26}.
\end{proof}

\begin{theorem} 
\label{thm 5.7}
For fixed $s_*\in(0,T)$ and $0<L\le T-s_*$,
there exists $\eta_L>0$, uniform for
$s_*\le s_0\le T-L$, such that every
$Z_0=(r_0,q_0,v_0)\in\mathcal Y$ satisfying
\[
\|Z_0\|_{\mathcal Y}\le\eta_L,\qquad
\bar n(s_0)+r_0\ge0,\qquad
\bar c(s_0)+q_0\ge0
\]
can be driven to the reference trajectory at time
$s_0+L$ by the exterior controls
\eqref{5.24}. More precisely, the
corresponding solution $(r,q,v)$ belongs to
\eqref{1.5}, satisfies
\[
r(L)=q(L)=v(L)=0
\qquad\text{in the relative time }t=s-s_0,
\]
and the biological states remain nonnegative:
\[
\bar n+r\ge0,\qquad \bar c+q\ge0
\quad\text{a.e. in }\mathcal O,\qquad 0\le t\le L.
\]
\end{theorem}

\begin{proof}
Starting from the zero state control pair, apply Proposition~\ref{prop 5.5} at step $k+1$ with initial state $Z_0$ and sources \eqref{5.22}-\eqref{5.23} evaluated at step $k$. If $C_R$ is the norm of the linear construction and $R_k$ is the norm of the $k$th pair, then
\[
 R_{k+1}\le C_R\|Z_0\|_\Y+C_RC R_k^2.
\]
Choose $R=2C_R\|Z_0\|_\Y$ and $\eta_L$ such that $2C_RC R\le1/2$. By Lemma~\ref{lemma 5.6}, we have
\[
 \|(Z^{k+1},h^{k+1})-(Z^k,h^k)\|_{\mathcal W}
 \le\tfrac12\|(Z^k,h^k)-(Z^{k-1},h^{k-1})\|_{\mathcal W}.
\]
The sequence converges in $\mathcal W$ to a solution of the nonlinear system, and \eqref{5.21} gives $Z(L)=0$. Finally,
\[
 \alpha_n=\an+h_n/\nb,\quad \alpha_c=\ac+h_c/\cb\in L^2H^1\subset L^2L^3,
 \qquad \nabla c\in L^\infty L^6\subset L^4L^6.
\]
By Lemma~\ref{lemma 2.2}, the limiting solution $n,c\ge0$.  
\end{proof}

 \section{The global arguments}\label{section 6}
We conclude by combining the preparation, approximate controllability, regularization and local controllability results.

\subsection{Regularization near a nonconstant reference}\label{subsection 6.1}
We return to physical time and choose the exterior sources $f_n=\an n$, $f_c=\ac c$ and $f_u=\bar f_u$. 

Fix $s_*>0$, $s_0\ge s_*$ and $L>0$ with $s_0+L\le T$. Let $t=s-s_0$ and 
\begin{equation*}
 r(t)=n(s_0+t)-\nb(s_0+t),\quad
 q(t)=c(s_0+t)-\cb(s_0+t),\quad
 v(t)=u(s_0+t)-\ub(s_0+t).
\end{equation*}
Let $\pi(t)=p(s_0+t)-\pb(s_0+t)$ and $(r_0,q_0,v_0)=(r,q,v)(0)$. Subtracting the reference equations, one can see
\begin{equation}\label{6.1}
\left\{\begin{aligned}
 &\begin{aligned}[b] &r_t-D_n\Delta r+\ub\cdot\nabla r+v\cdot\nabla\nb +\chi\Div(\nb\nabla q+r\nabla\cb)-\an r=-v\cdot\nabla r-\chi\Div(r\nabla q),\end{aligned} &&\text{in }Q_L,\\[2pt]
 &\begin{aligned}[b] &q_t-D_c\Delta q+\ub\cdot\nabla q+v\cdot\nabla\cb +(\nb-\ac)q+\cb r=-v\cdot\nabla q-rq,\end{aligned} &&\text{in }Q_L,\\[2pt]
 &\begin{aligned}[b] &v_t-\mu\Delta v+(\ub\cdot\nabla)v+(v\cdot\nabla)\ub +\nabla\pi-r e_d=-(v\cdot\nabla)v,\end{aligned} ~~\Div v=0, &&\text{in }Q_L,\\[2pt]
 &\partial_\nu r=\partial_\nu q=0,\quad v\cdot\nu=0,\quad\NN v=0, &&\text{on }\Sigma_L,\\[2pt]
 &(r,q,v)|_{t=0}=(r_0,q_0,v_0), &&\text{in }\OO.
\end{aligned}\right.
\end{equation}

\begin{lemma}
\label{lemma 6.1}
Fix $s_*\in(0,T)$ and $0<L\le T-s_*$. There are $\delta_L,C_L>0$, uniform for $s_*\le s_0\le T-L$. Assume that $(r_0,q_0,v_0)\in\Y$ satisfies
\begin{equation*}
 \|r_0\|_4+\|q_0\|_{W^{1,4}}+\|v_0\|_4\le\delta_L.
\end{equation*}
Then system \eqref{6.1} has a strong solution on $[0,L]$ and
\begin{align}
 \sup_{0\le t\le L}(\|r(t)\|_4+\|q(t)\|_{W^{1,4}}+\|v(t)\|_4)
 &\le C_L(\|r_0\|_4+\|q_0\|_{W^{1,4}}+\|v_0\|_4),\label{6.2}\\
 \sup_{0\le t\le L}\|(r,q,v)(t)\|_{\mathcal X_0}^2
 +\int_0^L\|(r,q,v)(t)\|_{\mathcal Y}^2\dd t
 &\le C_L\|(r_0,q_0,v_0)\|_{\mathcal X_0}^2.\label{6.3}
\end{align}
At some $t_*\in[L/4,L/2]$,
\begin{equation}\label{6.4}
 \|(r,q,v)(t_*)\|_{\mathcal Y}\le C_L\|(r_0,q_0,v_0)\|_{\mathcal X_0}.
\end{equation}
Nonnegative actual biological states remain nonnegative.  
\end{lemma}
\begin{proof}

Applying \eqref{3.29} and
\eqref{3.36} with $\varepsilon=1$,
$n^{\rm app}=\bar n$, $c^{\rm app}=\bar c$, $V^\varepsilon=\bar u$,
$\alpha_n^\varepsilon=\bar\alpha_n$, $\alpha_c^\varepsilon=\bar\alpha_c$,
and $R_n=R_c=F_u=h^\varepsilon=0$, on any interval
where
\[
E(t):=
\bigl(\|r(t)\|_4^4+\|q(t)\|_{W^{1,4}}^4\bigr)^{1/2}
+\|v(t)\|_4^2\le1,
\]
we obtain
\[
E'(t)\le CE(t).
\]
Then, using Gronwall's inequality, we get 
\[
E(t)\le e^{CL}E(0),\qquad 0\le t\le L.
\]
Choose the initial data such that
\[
e^{CL}\left[
\bigl(\|r_0\|_4^4+\|q_0\|_{W^{1,4}}^4\bigr)^{1/2}
+\|v_0\|_4^2
\right]<\frac12,
\]
which implies $E(t)<1/2$ throughout the interval of existence.
By continuity, $E(t)$ cannot reach $1$, and hence
\eqref{6.2} holds. Thus, Lemma~\ref{lemma 3.3}
 extends the strong solution to the whole interval
$[0,L]$.

We next derive an energy estimate that allows us to
select a time $t_*\in[L/4,L/2]$ at which the solution
belongs to the stronger space $\mathcal Y$. Testing the equations by $r$, $K(q-\Delta q)$ and $v$, the linear terms satisfy
\begin{align}
 &\chi\left|\int\bar n\nabla q\cdot\nabla r\right|
\le\tfrac18D_n\|\nabla r\|_2^2+C\|\nabla q\|_2^2,\label{6.5}\\
 &\chi\left|\int r\nabla\bar c\cdot\nabla r\right|
\le\tfrac18D_n\|\nabla r\|_2^2+C\|r\|_2^2,\label{6.6}\\
 &|(v\cdot\nabla\bar n,r)|+| (\bar\alpha_nr,r)|\le C\big[\|r\|_2^2+K(\|q\|_2^2+\|\nabla q\|_2^2)+\|v\|_2^2\big],\label{6.7}\\
& |(v\cdot\nabla\bar c+(\bar n-\bar\alpha_c)q+\bar c r,q-\Delta q)|
\le\eta\|\Delta q\|_2^2+C_\eta \big[\|r\|_2^2+K(\|q\|_2^2+\|\nabla q\|_2^2)+\|v\|_2^2\big],\label{6.8}\\
 &|((v\cdot\nabla)\bar u,v)|+|(r e_d,v)|\le C\big[\|r\|_2^2+K(\|q\|_2^2+\|\nabla q\|_2^2)+\|v\|_2^2\big].\notag
\end{align}
Denote $b_* =\sup_{[0,L]}(\|r\|_4+\|q\|_{W^{1,4}}+\|v\|_4)$. For the remaining nonlinear terms, we also have
\begin{align}
 \chi\left|\int r\nabla q\cdot\nabla r\right|
 &\le C\|r\|_3\|\nabla q\|_6\|\nabla r\|_2
 \le C b_*(\|\Delta q\|_2+\|q\|_2)\|\nabla r\|_2 \le Cb_*(\|\nabla r\|_2^2+\|\Delta q\|_2^2+\|q\|_2^2),\label{6.9}\\
 |(rq,q-\Delta q)|
 &\le\|q\|_\infty\|r\|_2(\|q\|_2+\|\Delta q\|_2) \le Cb_*(\|r\|_2^2+\|q\|_2^2+\|\Delta q\|_2^2),\label{6.10}\\
 |(v\cdot\nabla q,\Delta q)|
 &\le\|v\|_3\|\nabla q\|_6\|\Delta q\|_2
 \le Cb_*(\|q\|_2^2+\|\Delta q\|_2^2).\label{6.11}
\end{align}
If we take $K$ sufficiently large, the terms
involving $\|\nabla q\|_2^2$ in
\eqref{6.5}-\eqref{6.7} can be controlled
by the chemical energy. Then, choose $\eta>0$
sufficiently small to absorb the term
$\eta\|\Delta q\|_2^2$ in \eqref{6.8}.
Finally, using the smallness of $b_*$, the nonlinear
terms in \eqref{6.9}-\eqref{6.11}
are absorbed into the diffusion terms.

Thus, by \eqref{2.4}, we obtain
\[
 \begin{aligned}
 &\frac{\dd}{\dd t}\big[\|r\|_2^2+K(\|q\|_2^2+\|\nabla q\|_2^2)+\|v\|_2^2\big] +\gamma_{\rm E}(\|\nabla r\|_2^2+\|\Delta q\|_2^2+\|\nabla q\|_2^2+\|v\|_{H^1}^2)\\
 &\le C\big[\|r\|_2^2+K(\|q\|_2^2+\|\nabla q\|_2^2)+\|v\|_2^2\big] +Cb_*(\|\nabla r\|_2^2+\|\Delta q\|_2^2+\|r\|_2^2+\|q\|_2^2).
 \end{aligned}
\]
From \eqref{6.2} and choosing $Cb_*\le \gamma_{\rm E}/2$, we have
\[
 \begin{aligned}
 &\frac{\dd}{\dd t}\big[\|r\|_2^2+K(\|q\|_2^2+\|\nabla q\|_2^2)+\|v\|_2^2\big] +\tfrac{\gamma_{\rm E}}2(\|\nabla r\|_2^2+\|\Delta q\|_2^2+\|\nabla q\|_2^2+\|v\|_{H^1}^2)\\
 &\le C_L\big[\|r\|_2^2+K(\|q\|_2^2+\|\nabla q\|_2^2)+\|v\|_2^2\big],\\
 &\sup_{0\le t\le L}(\|r(t)\|_2^2+\|q(t)\|_{H^1}^2+\|v(t)\|_2^2) +\int_0^L(\|\nabla r\|_2^2+\|\Delta q\|_2^2+\|\nabla q\|_2^2+\|v\|_{H^1}^2)\dd t\\
 &\le C_L(\|r_0\|_2^2+\|q_0\|_{H^1}^2+\|v_0\|_2^2),
 \end{aligned}
\]
which implies
\eqref{6.3}. In particular,
\[
 \frac4L\int_{L/4}^{L/2}
      (\|r\|_{H^1}^2+\|q\|_{H^2}^2+\|v\|_{H^1}^2)\dd t
 \le\frac{C_L}{L}(\|r_0\|_2^2+\|q_0\|_{H^1}^2+\|v_0\|_2^2).
\]
Therefore, \eqref{6.4} holds at some $t_*\in[L/4,L/2]$.  Finally, applying Lemma~\ref{lemma 2.2} to the
first and second equations of \eqref{2.1} yields $n\ge0$ and $c\ge0$.
\end{proof}

\subsection{Global approximation}
By Proposition~\ref{prop 2.1}, we extend the 
initial data to $\mathcal O$ such that
\eqref{2.2} holds. Starting from these
extended data, we apply Lemma~\ref{lemma 2.4} with all
exterior sources set to zero. The resulting solution is
nonnegative and becomes smooth for every positive time.
Choose $T_r\in(0,T/16)$ within the existence interval.
Then $(n,c,u)(T_r)$ is smooth and satisfies the boundary
compatibility conditions, and its scalar
components are nonnegative.
Then, starting from this state at $T_r$,
we apply Lemma~\ref{lemma 2.5} with the prescribed
time upper bound $T/2$. Let $T_p\in(0,T/16)$ be the duration
of this preparation and set
\begin{equation}\label{6.12}
T_0=T_r+T_p<T/8,\qquad s_*=T_0.
\end{equation}
By Lemma~\ref{lemma 2.5}, the prepared state
satisfies
\[
\min_{\overline{\mathcal O}}n(T_0)>0,\qquad
\min_{\overline{\mathcal O}}c(T_0)>0.
\]
The solution on the initial interval retains the
prescribed physical initial trace $(n_0,c_0,u_0)$ and
satisfies \eqref{2.6}. Moreover, the reference
trajectory and its associated gains are smooth on
$[T_0,T]$, and there exists $m_{T_0}>0$ such that
\[
\bar n,\bar c\ge m_{T_0}
\quad\text{on }[T_0,T]\times\overline{\mathcal O}.
\]

We first note that the piecewise solutions can be joined without loss of time regularity. Let $f_j$ denote
one component of the solution on $(t_j,t_{j+1})$.
If the traces agree at every junction $t_j$, then
\begin{equation}\label{6.13}
\begin{aligned}
\partial_s f
&=\sum_j\mathbf1_{(t_j,t_{j+1})}\partial_s f_j
+\sum_{j\ge1}\bigl(f_j(t_j)-f_{j-1}(t_j)\bigr)\delta_{t_j} =\sum_j\mathbf1_{(t_j,t_{j+1})}\partial_s f_j
\end{aligned}
\end{equation}
in the sense of distributions. Thus, the distributional time derivative contains no Dirac mass at any joining time. Since only finitely
many pieces are joined on each $[\delta,T]$, the strong regularity estimates on the individual
segments combine to give the corresponding
piecewise strong regularity on $[\delta,T]$.
For the velocity, the return flow $U^0$ vanishes near
the endpoints of its construction interval by
\eqref{3.2}. Hence its divergence
$\sigma^0=\Div U^0$ also vanishes there, which implies the
prescribed divergence is compatible with the adjacent
stages. On the return stage, the actual velocity is
\[
 u(s)=\varepsilon^{-1}V^\varepsilon(t)+y(t),\qquad t=(s-s_0)/\varepsilon.
\]
It is the remainder $y$, after subtraction of its smooth boundary stress lifting, that belongs to $\X_u$. The prescribed field $V^\varepsilon$ is piecewise smooth with matching state traces. Consequently, for fixed $\varepsilon>0$, the joined actual velocity belongs to
\[
 H^1(\delta,T;L^2(\OO)^d)\cap L^2(\delta,T;H^2(\OO)^d)
 \cap C([\delta,T];H^1(\OO)^d),
\]
while its divergence is the prescribed exterior source. The scalar components have the corresponding strong regularity in \eqref{1.5}.
Combining this regularity with \eqref{2.6}
and the Sobolev embeddings, we obtain \eqref{1.9}. Since all exterior
sources vanish in $\Omega$ and the extended divergence
source is zero there, restriction to $\Omega$ gives
the physical system \eqref{1.1}.

Define the boundary controls by
\[
\left\{\begin{aligned}
h_n^*&=\partial_\nu n|_{\Gamma_c},&
h_c^*&=\partial_\nu c|_{\Gamma_c},\\
h_\nu^*&=(u\cdot\nu)|_{\Gamma_c},&
h_\tau^*&=(\mathcal N_\beta u)|_{\Gamma_c}.
\end{aligned}\right.
\]
By virtue of the strong regularity on $[\delta,T]$, we have
\[
n,u\in L^2(\delta,T;H^2(\Omega)),\qquad
c\in L^2(\delta,T;H^3(\Omega)).
\]
Then, by the Sobolev trace theorem, \eqref{1.10} holds.
 The homogeneous boundary
conditions on $\Gamma_0$ follow from the extended
boundary conditions. Moreover, we have
\[
\int_{\Gamma_c}h_\nu^*\,dS
=\int_{\partial\Omega}u\cdot\nu\,dS
=\int_\Omega\Div u\,dx=0.
\]
Finally, nonnegativity of $n$ and $c$ follows from
Lemma~\ref{lemma 2.2} and Lemma~\ref{lemma 2.4}.

\begin{proof}[Proof of Theorem~\ref{thm 1.3}]
Let $T_0$ be given by \eqref{6.12},
and set
\[
T_{\rm app}=T-T_0,\qquad
s_0=T_0,\qquad
(n_{\rm in},c_{\rm in},u_{\rm in})=(n,c,u)(T_0).
\]
By \eqref{6.12}, the initial state
is smooth, compatible and strictly positive. Applying
the construction of Section~\ref{section 4} on
$[T_0,T]$, there is an extended controlled solution whose
initial trace is exactly $(n_{\rm in},c_{\rm in},u_{\rm in})$ and whose
terminal error satisfies \eqref{4.34}.

Restricting this solution to $\Omega$, we obtain a solution
of the original system \eqref{1.1}, since the
exterior sources vanish in $\Omega$. \eqref{4.34} yields \eqref{1.8}. 
The state and control regularity
in \eqref{1.9} and \eqref{1.10}
follows from the strong estimates on $[\delta,T]$,
the Sobolev trace theorem, and \eqref{6.13}. Nonnegativity is
preserved on each construction interval and at every
junction.

According to  \eqref{1.8}, the terminal error is bounded by
$C\varepsilon^{1/8}$. Hence, for any prescribed
$\eta>0$, we can take 
\[
0<\varepsilon<
\min\left\{\varepsilon_0,\,
\left(\frac{\eta}{C}\right)^8,\,
\frac{T_{\rm app}}{2}\right\}
\]
such that $C\varepsilon^{1/8} \leq \eta$.
This proves the global approximate controllability.
\end{proof}

\subsection{Global exact controllability}
\begin{proof}[Proof of Theorem~\ref{thm 1.4}]
By \eqref{6.12}, we set
\[
 T_{\rm app}=T/4,\qquad L_s=T/4,\qquad L_c=T/4,\qquad s_1=T_0+T_{\rm app}.
\]
With the help of the approximation construction on $[T_0,s_1]$, we have
\[
 (r_1,q_1,v_1)=(n-\nb,c-\cb,u-\ub)(s_1).
\]
It follows from \eqref{4.34} that
\[
 \|r_1\|_4+\|q_1\|_{W^{1,4}}+\|v_1\|_4
 \le C_{\rm app}\varepsilon^{1/8}.
\]
Let $\delta_{L_s},C_s>0$ be the constants given by
Lemma~\ref{lemma 6.1} with $C_s$ enlarged, if necessary, to absorb
the constant in the embedding
$L^4\times W^{1,4}\times L^4\hookrightarrow\mathcal X_0$.
Let $\eta_{L_c}>0$ denote the uniform radius provided
by Theorem~\ref{thm 5.7} for $s_*=T_0$. Choose $\varepsilon>0$ such that
\begin{equation}\label{6.14}
 C_{\rm app}\varepsilon^{1/8}\le\delta_{L_s},\qquad
 C_sC_{\rm app}\varepsilon^{1/8}\le\eta_{L_c}.
\end{equation}
From $s_1$, we use the exterior sources $\an n$, $\ac c$, $\bar f_u$ and
$\sigma=0$. Thus, by Lemma~\ref{lemma 6.1}, there is a time
$t_*\in[L_s/4,L_s/2]$ at which, with $s_2=s_1+t_*$,
\[
 \|(n-\nb,c-\cb,u-\ub)(s_2)\|_{\Y}
 \le C_sC_{\rm app}\varepsilon^{1/8}\le\eta_{L_c}.
\]

Since the state at $s_2$ satisfies the smallness and
nonnegativity conditions of Theorem~\ref{thm 5.7},
we apply Theorem~\ref{thm 5.7} on $[s_2,s_2+L_c]$ to obtain a controlled solution satisfying
\[
(n,c,u)(s_3)=(\bar n,\bar c,\bar u)(s_3),
\qquad s_3=s_2+L_c.
\]

Notice that all reference dependent constants in \eqref{6.14} are
fixed on $[T_0,T]$, independently of $\varepsilon$, and
\[
 s_3\le T_0+T/4+T/8+T/4
       =T_0+5T/8<3T/4<T.
\]
For $s\in[s_3,T]$, we set
$(n,c,u)=(\bar n,\bar c,\bar u)$ and use the
corresponding exterior sources
\[
f_n=\bar\alpha_n\bar n=\bar f_n,\qquad
f_c=\bar\alpha_c\bar c=\bar f_c,\qquad
f_u=\bar f_u.
\]
By virtue of \eqref{6.13}, we obtain the required regularity and nonnegativity  across
all junctions. Restriction to $\Omega$ proves \eqref{1.11}, which completes the proof of Theorem~\ref{thm 1.4}.
\end{proof}

\section*{Acknowledgement}
Tao was partially supported by the National Natural Science Foundation of China (12471423),
by the Guangdong Basic and Applied Basic Research Foundation (2025A1515010485), and by the
Shenzhen Science and Technology Program (RCJC20221008092753082).
Yang
was partially supported by the Shenzhen Science and Technology Program (JCYJ20241202124209011),
and by the Shenzhen Natural Science Fund (the Stable Support Plan Program 20231122110251002).


\vskip8mm


\begin{thebibliography}{99}

\bibitem{AS2006}
A.A. Agrachev, A.V. Sarychev,
Controllability of 2D Euler and Navier--Stokes equations by degenerate forcing,
Comm. Math. Phys. 265 (2006) 673--697.

\bibitem{AEG2018}
C. Amrouche, M. Escobedo, A. Ghosh,
Semigroup theory for the Stokes operator with Navier boundary condition on $L^p$ spaces,
in: T. Bodn\'ar, G.P. Galdi, \v S. Ne\v casov\'a (Eds.), Waves in Flows,
Birkh\"auser, Cham, 2021, pp.~1--51.
 

\bibitem{CarG2013}
N. Carre\~no, S. Guerrero,
Local null controllability of the $N$-dimensional Navier--Stokes system with $N-1$ scalar controls in an arbitrary control domain,
J. Math. Fluid Mech. 15 (2013) 139--153.

\bibitem{CFKMS2023}
F.W. Chaves-Silva, E. Fern\'andez-Cara, K. Le Balc'h, J.L.F. Machado, D.A. Souza,
Global controllability of the Boussinesq system with Navier-slip-with-friction and Robin boundary conditions,
SIAM J. Control Optim. 61 (2023) 484--510.
 
\bibitem{CG2015}
F.W. Chaves-Silva, S. Guerrero,
A uniform controllability result for the Keller--Segel system,
Asymptot. Anal. 92 (2015) 313--338.

\bibitem{CG2017}
F.W. Chaves-Silva, S. Guerrero,
A controllability result for a chemotaxis--fluid model,
J. Differential Equations 262 (2017) 4863--4905.

\bibitem{CL2014}
J.-M. Coron, P. Lissy,
Local null controllability of the three-dimensional Navier--Stokes system with a distributed control having two vanishing components,
Invent. Math. 198 (2014) 833--880.

\bibitem{CMS2020}
J.-M. Coron, F. Marbach, F. Sueur,
Small-time global exact controllability of the Navier--Stokes equation with Navier slip-with-friction boundary conditions,
J. Eur. Math. Soc. 22 (2020) 1625--1673.


\bibitem{CMSZ2019}
J.-M. Coron, F. Marbach, F. Sueur, P. Zhang,
Controllability of the Navier--Stokes equation in a
rectangle with a little help of a distributed phantom force,
Ann. PDE 5 (2019), Article 17.

\bibitem{FGIP2004}
E. Fern\'andez-Cara, S. Guerrero, O.Yu. Imanuvilov, J.-P. Puel,
Local exact controllability of the Navier--Stokes system,
J. Math. Pures Appl. 83 (2004) 1501--1542.

\bibitem{G2006}
S. Guerrero,
Local exact controllability to the trajectories of the Navier--Stokes system with nonlinear Navier-slip boundary conditions,
ESAIM Control Optim. Calc. Var. 12 (2006) 484--544.

\bibitem{GZ2014}
B.  Guo, L. Zhang,
Local null controllability for a chemotaxis system of parabolic--elliptic type,
Systems Control Lett. 65 (2014) 106--111.

\bibitem{GZ2016}
B.  Guo, L. Zhang,
Local exact controllability to positive trajectory for parabolic system of chemotaxis,
Math. Control Relat. Fields 6 (2016) 143--165.

\bibitem{I2001}
O.Yu. Imanuvilov,
Remarks on exact controllability for the Navier--Stokes equations,
ESAIM Control Optim. Calc. Var. 6 (2001) 39--72.

\bibitem{IP2002}
O.Yu. Imanuvilov, J.-P. Puel,
Global Carleman estimates for weak solutions of elliptic nonhomogeneous Dirichlet problems,
C. R. Math. Acad. Sci. Paris 335 (2002) 33--38.

\bibitem{KS1970} E.F. Keller, L.A. Segel, Initiation of slime mold aggregation viewed as an instability, J. Theoret. Biol. 26 (1970) 399–415.

\bibitem{KS1971}  E.F. Keller, L.A. Segel, Model for chemotaxis, J. Theoret. Biol. 30 (1971) 225–234.

\bibitem{K2025}
P.C. Kunstmann,
$H^\infty$-calculus for the Stokes operator with Hodge, Navier, and Robin boundary conditions on unbounded domains,
\href{https://arxiv.org/abs/2504.18895v1}{arXiv:2504.18895v1} (2025).

\bibitem{LSZ2022}
J. Liao, F. Sueur, P. Zhang,
Smooth controllability of the Navier--Stokes equation with Navier conditions: application to Lagrangian controllability,
Arch. Ration. Mech. Anal. 243 (2022) 869--941.

\bibitem{LTZ2017}
J. Loh\'eac, E. Tr\'elat, E. Zuazua,
Minimal controllability time for the heat equation
under unilateral state or control constraints,
Math. Models Methods Appl. Sci. 27 (2017) 1587--1644.

\bibitem{NRB2025}
V. Nersesyan, M. Rissel,
Global controllability of Boussinesq flows by using only a temperature control,
Arch. Ration. Mech. Anal. 249 (2025), Article 58.

\bibitem{NRNS2025}
V. Nersesyan, M. Rissel,
Localized and degenerate controls for the incompressible Navier--Stokes system,
Comm. Pure Appl. Math. 78 (2025) 1285--1319.

\bibitem{PZ2018}
D. Pighin, E. Zuazua,
Controllability under positivity constraints of
semilinear heat equations,
Math. Control Relat. Fields 8 (2018) 935--964.

\bibitem{ProtterWeinberger1984}
M.H. Protter, H.F. Weinberger,
\emph{Maximum Principles in Differential Equations},
Springer-Verlag, New York, 1984, Chapter~III, pp.~159--194.

\bibitem{Saito2004}
N. Saito,
On the Stokes equation with the leak and slip boundary conditions of friction type: regularity of solutions,
Publ. Res. Inst. Math. Sci. 40 (2004) 345--383.

\bibitem{TYY2025}
Q. Tao, Z. Yao, X. Yin,
Global controllability of the micropolar fluids with Navier slip-with-friction and Dirichlet boundary conditions,
J. Differential Equations 440 (2025) 113437.

 

\bibitem{TZ2025}
Q. Tao, M. Zhang,
Exact controllability to nonnegative trajectory for a chemotaxis system,
Appl. Math. Optim. 91 (2025), Article 21.

\bibitem{TCDWKG2005} I. Tuval, L. Cisneros, C. Dombrowski, C.W. Wolgemuth, J.O. Kessler, R.E. Goldstein, Bacterial swimming and oxygen transport near contact lines, Proc. Natl. Acad. Sci. USA 102 (2005) 2277–2282.

\end{thebibliography}
\end{document}